%% file: 00-main.tex
\documentclass{article}
\usepackage[a4paper,margin=1in]{geometry}

\usepackage{amsmath}
\usepackage{amssymb}
\usepackage{amsthm}
\usepackage{boldtensors}
\usepackage{geometry}
\usepackage{graphicx}
\usepackage{lineno}
\usepackage{mathtools}
\usepackage{microtype}
\usepackage{multirow}
\usepackage{url} %this package should fix any errors with URLs in refs.
\usepackage{xcolor}
\theoremstyle{plain}
\usepackage{algpseudocode}
\usepackage{booktabs}
\usepackage{makecell}
\usepackage{tikz}
\usepackage{xcolor}
\usepackage[table,xcdraw]{xcolor}
\usepackage[normalem]{ulem}
\usepackage[backend=bibtex,maxcitenames=2]{biblatex}
\usepackage[font=normalsize]{caption}
\usepackage{subcaption}
\usepackage{xparse}

\usepackage{amsthm}
\usepackage[ruled,linesnumbered]{algorithm2e}
\newtheorem{remark}{Remark}

\newcommand{\lr}[1]{\left(#1\right)}

\input{commands}

\begin{document}

% \title{Domain decomposition-based coupling of non-intrusive operator inference reduced order models via the Schwarz alternating method}
\title{Hybrid coupling with operator inference and the overlapping Schwarz alternating method}

\author{Irina Tezaur$^1$\thanks{Email: ikalash@sandia.gov}, Eric Parish$^{1}$, Anthony Gruber$^1$, Ian Moore$^{1,2}$, \\ Christopher R. Wentland$^1$, Alejandro Mota$^1$    
  \\
  \\
  \small $^1$Sandia National Laboratories, Livermore, CA and Albuquerque, NM, USA\\
  \small $^2$Rice University, Houston, TX, USA\\
}

\date{}

\maketitle

%\ikt{Author list is currently alphabetical, with Irina last.  We can decide on the author order later on once the paper is fleshed out more.  Also, title should probably be shortened.}

% \adg{I took a stab at the title and left a review comment about my place in the author order.}

%\ikt{What do folks think about adding ``hybrid models" to the title?} \adg{Something like ``Coupled hybrid modeling with operator inference and the overlapping Schwarz alternating method'' would work, but then we lose the subdomain-local piece.  Or maybe ``Hybrid coupling with...''?}

%\todo{Throughout paper, switch $^T$ to $^{\intercal}$ to denote transpose.  I did this in a lot of places but might have missed some. \\}

%\todo{Irina.  State explicitly that we are considering physically-motivated DDs and not something like Ramin and Youngsoo where you're tiling together hundreds of small subdomains.\\}

%\ikt{I introduce the acronym O-SAM but then mostly use SAM throughout the paper.  We could cut the ``O" or try to use O-SAM in more places.  Thoughts?  Maybe it is OK as it.} \adg{I think O-SAM is cute and worth using, but it's probably okay as it is. \\} 

%\ikt{I don't know if we explicitly say anywhere that we do multiplicative Schwarz.  Should we?\\}

%\ikt{In terms of where to submit, here is my list, in order of preference: IJNME, CMAME, Computational Mechanics, FINEL.  Any thoughts / alternate suggestions are welcome.} \adg{I agree with this list. \\}

% Abstract
%\begin{center}
\noindent {\bf Abstract.  } This paper presents a 
%novel
hybrid approach for coupling subdomain-local, non-intrusive Operator Inference (OpInf) reduced order models (ROMs) with each other and with subdomain-local, high-fidelity full order models (FOMs) using the overlapping Schwarz alternating method (O-SAM). The proposed methodology addresses significant challenges in multiscale modeling and simulation, particularly the long runtime and complex mesh generation requirements associated with traditional high-fidelity simulations. By leveraging the flexibility of O-SAM, we enable the seamless integration of disparate models, meshes, and time integration schemes, enhancing computational efficiency while maintaining high accuracy. Our approach is demonstrated through a series of numerical experiments on increasingly complex three-dimensional solid dynamics problems, showcasing speedups of up to 106 times compared to conventional FOM-FOM couplings. %The results suggest that the spatial localization and FOM coupling enabled by SAM integration of OpInf allows for a non-intrusive coupling strategy that significantly reduces implementation burdens and enhances the predictive capabilities of ROMs. 
This work paves the way for more efficient simulation workflows in engineering applications, with potential extensions to a wide range of partial differential equations. \\

%\end{center}

\input{01-intro}

% \todo{Anthony.}

%\ikt{We could give here the model problem, e.g., if we decide to show just solid mechanics results from Norma.  Describe linear and quadratic OpInf models.}

\input{02-solid-mechanics}
\input{03-opinf}

% Schwarz method

%\noindent \ikt{Optional topics that could potentially be added (don't write up yet): 
%\begin{itemize}
%    \item Centering approach for BCs - only include if we feel it's needed.
%    \item If non-overlapping works and is not much different, can include this in the paper.  Evaluate in mid-August after Cam is done with his internship.  OpInf-FOM with D-N BCs might be low-hanging fruit that should work as is. 
%\end{itemize}
%}

%\ikt{10/1/25 comment on above: I would say lets skip both centering and NO-SAM.  NO-SAM will be its own paper.}

\input{04-schwarz}

%\input{05-cost-analysis}

\input{06-results}

%\input{results}

% Conclusions

\input{07-conclusion}

\input{08-appendix}

\input{reproducibility}

\input{acknowledgments}

% Appendices...
%\input{AppendixA}

\printbibliography

\end{document}

%% file: commands.tex
\newcommand{\numTime}{\tau}

\newcommand{\timeInit}{\ensuremath{t_0}}
\newcommand{\timeFinal}{\ensuremath{T}}

\newcommand{\domain}{\Omega}
\newcommand{\domainOneArg}[1]{\domain_{#1}}

\newcommand{\domainBoundaryOneArg}[1]{\partial \Omega_{#1}}
\newcommand{\domainBoundaryDirichlet}{\partial_{\varphi} \Omega }
\newcommand{\domainBoundaryDirichletOneArg}[1]{\partial_{\varphi} \Omega_{#1} }
\newcommand{\domainBoundaryNeumann}{\partial_{T} \Omega }
\newcommand{\domainBoundaryNeumannOneArg}[1]{\partial_{T} \Omega_{#1}}
\newcommand{\domainBoundarySchwarzOneArg}[1]{\partial_{S} \Omega_{#1}}
\newcommand{\mass}{\boldsymbol{M}}
\newcommand{\stiffness}{\boldsymbol{K}}

\newcommand{\displacement}{\boldsymbol u}
\newcommand{\velocity}{\dot{\displacement}}
\newcommand{\acceleration}{\ddot{\displacement}}

\newcommand{\bounadryDirichletOneArg}[1]{\partial_{\chi} \Omega}
\newcommand{\bounadryTractionOneArg}[1]{\partial_{T} \Omega}
\newcommand{\schwarzIt}{n}
\newcommand{\numControllerSteps}{C}
\newcommand{\numTimeStepsOneArg}[1]{n_{#1}}
\newcommand{\numTotalTimeStepsOneArg}[1]{\tau_{#1}}

\newcommand{\timeStepOneArg}[1]{\Delta t_{#1}}
\newcommand{\accelerationThreeArg}[3]{\acceleration_{#1}^{#2,#3}}
\newcommand{\accelerationTwoArg}[2]{\acceleration_{#1}^{#2}}
\newcommand{\displacementThreeArg}[3]{\displacement_{#1}^{#2,#3}}
\newcommand{\displacementTwoArg}[2]{\displacement_{#1}^{#2}}
\newcommand{\displacementOneArg}[1]{\displacement_{#1}}
\newcommand{\velocityTwoArg}[2]{\velocity_{#1}^{#2}}
\newcommand{\velocityThreeArg}[3]{\velocity_{#1}^{#2,#3}}
\newcommand{\schwarzBoundaryBCTwoArg}[2]{\boldsymbol s_{#1}^{#2}}
\newcommand{\schwarzBoundaryBCThreeArg}[3]{\boldsymbol s_{#1}^{#2,#3}}
\newcommand{\numFreeDofsOneArg}[1]{N_{#1}}
\newcommand{\numDirichletBoundaryDofsOneArg}[1]{N_{#1}^B}
\newcommand{\numSchwarzDirichletBoundaryDofsOneArg}[1]{N_{#1}^S}
\newcommand{\spatialSchwarzProjectionOperator}{~\Pi}
\newcommand{\temporalSchwarzProjectionOperator}{~\Xi}
\newcommand{\trialVarCap}{\Phi}

\newcommand{\trialBasis}{\ensuremath{\mathbf{\trialVarCap}}}

\newcommand{\displacementSnapshots}{\boldsymbol U}
\newcommand{\dirichletBoundarySnapshots}{\boldsymbol G}
\newcommand{\reducedDirichletBoundarySnapshots}{\overline{\boldsymbol G}}
\newcommand{\regularizationParameter}{\lambda}
\newcommand{\reducedAcceleration}{\ddot{\reducedDisplacement}}
\newcommand{\reducedVelocity}{\dot{\reducedDisplacement}}
\newcommand{\reducedDisplacement}{\hat{\displacement}}
\newcommand{\reducedBoundaryBC}{\hat{\boldsymbol g}}
\newcommand{\reducedStiffness}{\bar{\boldsymbol{K}}}
\newcommand{\reducedStiffnessQuadratic}{\bar{\boldsymbol{H}}}
\newcommand{\reducedStiffnessCubic}{\bar{\boldsymbol{C}}}
\newcommand{\reducedBoundaryOperator}{\bar{\boldsymbol{B}}}

\newcommand{\romDimOneArg}[1]{r_{#1}}
\newcommand{\reducedBoundaryBCOneArg}[1]{\reducedBoundaryBC_{#1}}

\newcommand{\reducedBoundaryBCThreeArg}[3]{\reducedBoundaryBC_{#1}^{#2,#3}}
\newcommand{\reducedAccelerationOneArg}[1]{\reducedAcceleration_{#1}}
\newcommand{\reducedAccelerationTwoArg}[2]{\reducedAcceleration_{#1}^{#2}}
\newcommand{\reducedAccelerationThreeArg}[3]{\reducedAcceleration_{#1}^{#2,#3}}

\newcommand{\reducedVelocityTwoArg}[2]{\reducedVelocity_{#1}^{#2}}
\newcommand{\reducedVelocityThreeArg}[3]{\reducedVelocity_{#1}^{#2,#3}}
\newcommand{\reducedDisplacementOneArg}[1]{\reducedDisplacement_{#1}}
\newcommand{\reducedDisplacementTwoArg}[2]{\reducedDisplacement_{#1}^{#2}}
\newcommand{\reducedDisplacementThreeArg}[3]{\reducedDisplacement_{#1}^{#2,#3}}
\newcommand{\trialBasisDirichlet}{\trialBasis_{\domainBoundaryDirichletOneArg{1}}}
\newcommand{\trialBasisDirichletTwo}{\trialBasis_{\domainBoundaryDirichletOneArg{2}}}
\newcommand{\trialBasisSchwarz}{\trialBasis_{\domainBoundarySchwarzOneArg{1}}}
\newcommand{\trialBasisSchwarzTwo}{\trialBasis_{\domainBoundarySchwarzOneArg{2}}}
\newcommand{\numReducedDirichletBoundaryDofsOneArg}[1]{r_{#1}^{\mathrm{B}}}
\newcommand{\numReducedSchwarzDirichletBoundaryDofsOneArg}[1]{r_{#1}^{\mathrm{S}}}
\newcommand{\schwarzBoundarySnapshots}{\boldsymbol S}
\newcommand{\reducedSchwarzBoundarySnapshots}{\overline{  \schwarzBoundarySnapshots}}
\newcommand{\approximateSchwarzBoundaryBCTwoArg}[2]{\tilde{\boldsymbol s}_{#1}^{#2}}
\newcommand{\approximateSchwarzBoundaryBCThreeArg}[3]{\tilde{\boldsymbol s}_{#1}^{#2,#3}}

\newcommand{\RR}[1]{\mathbb{R}^{#1}}
\newcommand{\RRStar}[2]{\mathbb{V}_{#2}\left( \RR{#1} \right)}
\newcommand{\defEq}{\ensuremath{:=}}
\DeclareMathOperator*{\argmin}{arg\,min}

%% file: 01-intro.tex
% Introduction
\section{Introduction} \label{sec:intro}

Multiscale modeling and simulation (mod/sim) is crucial in engineering, as it enables the design, optimization and qualification of complex engineered components and systems across a variety of applications,
enhancing the understanding of sophisticated multiscale and multiphysics processes.  % without requiring
%costly experiments or field tests. 
Unfortunately, analysts running traditional high-fidelity simulation
codes often face significant delays due to not only long runtime requirements, but also to the mesh
generation step of the mod/sim workflow: creating a high-quality mesh for a single multiscale component
can take weeks, making it ``the single biggest bottleneck in [mod/sim analyses]" \cite{sandialabnews}. During the past two
decades, projection-based reduced order models (ROMs) have emerged as a promising data-driven tool 
with potential to reduce the online computational complexity of numerical simulations, especially
in multi-query applications, such as design optimization or uncertainty quantification (UQ). However, these
models can face their own shortcomings, including a lack of systematic refinement mechanisms, a lack of
stability and/or accuracy in the predictive regime, and lengthy implementation time requirements. While
recently-proposed physics-informed neural networks (PINNs) \cite{Raissi:2019} and Physics-Informed Deep Operator Networks (PI-DeepONets) \cite{Mandl:2025} have been advertised as mesh-free
methods that can bypass meshing by relying %instead 
on a set of scattered collocation points at which
the governing partial differential equations (PDEs) are evaluated, these models lack interpretability, can
be 
%incredibly 
costly to train, and can suffer from convergence issues. This is because they effectively replace a set of linear
algebraic solves arising from a traditional (e.g., finite element or finite difference) discretization with a
complex nonlinear optimization problem involving the neural network’s (NN’s) parameters.

Alternatively, this paper develops a novel, hybrid, domain decomposition- (DD-)based approach that has the potential to mitigate both the meshing and long runtime requirement issues described above by ``gluing together," in a plug-and-play fashion, arbitrary combinations of subdomain-local, high-fidelity full order models (FOMs) with subdomain-local ROMs using the Schwarz alternating method (SAM)~\cite{Schwarz:1870}. The key idea behind SAM is to first decompose the physical domain on which a given PDE is posed into smaller subdomains, and then to solve a sequence of problems on each subdomain, while exchanging boundary condition information across subdomain interfaces to ensure solution compatibility. The approach presented herein builds on 
our
past work in developing SAM as a means to achieve concurrent multiscale coupling of high-fidelity finite element models in solid mechanics~\cite{Mota:2017, Mota:2022}, and on intrusive projection-based ROMs in hyperbolic fluid problems discretized using the cell-centered finite volume method~\cite{wentland2024Schwarz}. 
% It was demonstrated in these manuscripts that
SAM 
has a number of advantages demonstrated in these manuscripts: (i) it is minimally intrusive to implement in existing HPC software frameworks; (ii) it is capable of coupling regions with different mesh resolutions, different element types, different time integration schemes (e.g., implicit and explicit), and even different models (e.g., FOM and ROM), all without introducing any artifacts exhibited by alternative coupling methods; and (iii) it possesses rigorous convergence properties/guarantees~\cite{Mota:2017, Mota:2022}. Additionally, 
%we argue in 
\cite{wentland2024Schwarz} argues that DD-based FOM-ROM and ROM-ROM couplings %such as those proposed herein 
have the potential to improve the predictive viability of ROMs by enabling their spatial localization via domain decomposition, as well as the online integration of high-fidelity information 
%into these models 
via FOM coupling. 

The aim of the present work is to extend SAM further to enable the DD-based coupling of non-intrusive ROMs constructed via an approach known as operator inference, or OpInf, pioneered by Peherstorfer and Willcox~\cite{willcox2016opinf}. 
%\ikt{I really like the previous sentence and stand by it as a main motivation for this work.  However, one valid criticism from a reviewer could be that we don't show that the localization/coupling is improving accuracy relative to a monolithic ROM.  I did not show those results for a few reasons including the length of the paper and difficulty getting a monolithic mesh (for the bolted joint problem) that would run with Norma, but also the fact that some of the monolithic ROMs I tried were actually pretty accurate.  I do believe it will not always be the case.  I am OK keeping the sentence and seeing what the reviewers say, but I wanted to share this concern with the team to see if there are any thoughts/comments.}
In developing ROM-ROM and FOM-ROM couplings, our goal is to reduce the implementational burden in combining SAM with model order reduction (MOR), making it accessible to a wider range of applications and codes. Traditional intrusive MOR requires access to the underlying FOM code in order to project the governing PDE(s) onto a reduced subspace, posing a high development cost. OpInf instead assumes a functional form (usually globally linear or quadratic~\cite{Geelen:2023}) for the ROM in terms of reduced operators which are learned offline from data.
% This is possible with OpInf since, unlike traditional intrusive MOR, which requires access to the underlying FOM code in order to project the governing PDE(s) onto a reduced subspace, OpInf works by assuming a functional form (usually linear or quadratic~\cite{Geelen:2023}) for the ROM in terms of to-be-learned reduced operators within a latent space spanned by a reduced (e.g., a linear Proper Orthogonal Decomposition or a quadratic) basis, and solving an optimization problem offline for these operators.
Importantly, this procedure can be implemented entirely outside the FOM code, significantly reducing both the development time and the time-to-impact.

We focus our attention herein on the simplest version of SAM, namely overlapping SAM (O-SAM), in which the physical geometry is decomposed into overlapping subdomains and coupling is achieved via Dirichlet transmission boundary conditions on the subdomain boundaries; further extensions to non-overlapping domain decompositions have been explored in the preliminary work~\cite{Rodriguez:2025} and will be the subject of a future publication. Since we target an application space in which system geometries can be too complex to mesh using a monolithic meshing scheme, both \textit{a priori} ROM training and \textit{a posteriori} accuracy and runtime assessments are performed with respect to a SAM-based coupled simulation involving finite element FOMs. The assumption that a FOM-FOM coupled solution on the full physical domain of interest is possible to obtain via SAM enables us to utilize a ``top-down" (vs. a ``bottom-up"~\cite{Chung:2024}) training approach, which does not require ROM subdomains to be simulated independently, and generally gives rise to more expressive reduced bases and more accurate ROMs. While 
%we develop our method
the present method is developed around a solid dynamics exemplar, our SAM-based coupling strategy described herein can be applied to any set of PDEs. We note that the proposed approach 
%can be applied 
is also applicable in the case ``bottom-up" training is utilized to build the ROMs being coupled.

% We describe the
The specific contributions and differentiating features
of our approach are described in Section~\ref{sec:contribs}, after surveying the literature for related past work in Section~\ref{sec:past_work}.

\subsection{Overview of related past work} \label{sec:past_work}

The method proposed herein is an addition to a growing body of literature on DD-based couplings involving a variety of ROM formulations. For the sake of brevity, we will focus our literature review on Schwarz and Schwarz-like coupling methods for FOM-ROM and ROM-ROM couplings involving non-intrusive ROMs, as this class of methods is most closely related to the approach developed in this paper. A more detailed overview that includes a comprehensive survey of a variety of coupling methods involving intrusive projection-based ROMs can be found in 
%, the interested reader is referred to
\cite{wentland2024Schwarz} and \cite{Ruan2026}.%\ikt{I am trying to keep the lit review in this paper shorter relative to past papers. I am thinking to avoid discussing intrusive PROM couplings, especially non-Schwarz ones, and to just cite our earlier paper with Chris, which goes into this ad nauseam. Do people think this is bad? I know the literature really well for PROM and other types of coupling (LM, OBC) and can add it, but it will lengthen the paper.}

Some of the earliest Schwarz-based DD approaches for ROM-ROM and ROM-FOM coupling are based on non-intrusive ROMs constructed via Galerkin-free Proper Orthogonal Decomposition (POD). Galerkin-free POD ROMs are a class of reduced order modeling techniques that aim to efficiently approximate the behavior of complex dynamical systems without relying on the traditional Galerkin projection method. Instead of projecting the governing equations onto a reduced POD basis, as done in traditional intrusive projection-based MOR, Galerkin-free POD ROMs predict the solution by applying interpolation techniques to a POD representation of the primary solution field. Three methods to perform Galerkin-free ROM-FOM and ROM-ROM coupling are presented by Buffoni \textit{et al.}~\cite{Buffoni:2009}: (i) a Schur iteration where the solution of the PROM is obtained by a projection step in the space spanned by the POD modes, (ii) a Dirichlet–Dirichlet iteration in the frame of a classical Schwarz method, and (iii) an approach obtained by minimizing the residual norm of the canonical approximation in the space spanned by the POD modes. Another Galerkin-free ROM-FOM and ROM-ROM coupling approach that makes use of SAM is the work of Cinquegrana \textit{et al.}~\cite{Cinquegrana:2011}, but, unlike our approach, requires matching meshes at the subdomain interfaces. 
%are required. 
A third Galerkin-free POD approach, termed zonal Galerkin-free POD~\cite{Bergmann:2018}, defines an optimization problem that minimizes the difference between the POD reconstruction and its corresponding FOM solution in the overlapping region between a ROM and a FOM domain.

Schwarz-based methods have also been explored as a means to couple subdomain-local ROMs obtained non-intrusively via the proper generalized decomposition (PGD). In the recent pre-print~\cite{discacciati2025pgdbasedlocalsurrogatemodels}, Discacciati \textit{et al.} develop an algebraic SAM for coupling together subdomain-local non-intrusive PGD ROMs following a domain decomposition of the underlying spatial domain into overlapping subdomains. The proposed method is %effectively 
an improvement of the method originally presented in reference~\cite{Discacciati:2024}. %However, 
The coupling strategies proposed in~\cite{Discacciati:2024, discacciati2025pgdbasedlocalsurrogatemodels} are currently limited to linear elliptic PDEs and ROM-ROM (vs. ROM-FOM) couplings.

Recent years have seen the development of coupling methods for various types of NN-based ROMs. In~\cite{Wang:2022}, Wang \textit{et al.} introduce a genomic flow network (GFNet), in which pre-trained PINNs or NNs are trained locally on small ``genomes” (subdomains) and stitched together via a Schwarz-like iteration to represent the solution on an arbitrary domain comprised of the subdomain-local genomes. While the approach is presented as a non-overlapping method, convergence requires the introduction of ``auxiliary” genomes, which effectively introduce overlap into the domain decomposition. There are several other methods that utilize Schwarz for PINN coupling~\cite{LiD3M, LiDeepDDM, Snyder:2023}; however, the primary purpose of these approaches is to accelerate the training stage of the PINN construction process, rather than to glue together pre-trained subdomain-local models, as we attempt to do herein. The work by Discaciatti and Hesthaven~\cite{Discacciati:2024p2} similarly proposes a Schwarz-like method leveraging NNs, where the authors assume that local models assigned to non-overlapping subdomains are glued together using the relaxed Dirichlet-Neumann SAM, and construct a surrogate for the Dirichlet and Neumann maps defining the coupling conditions using ROMs based on kernel interpolation methods and NNs. One disadvantage of this method, as presented in~\cite{Discacciati:2024p2}, is that it does not automatically provide the local solution in specific subdomains, as its main purpose is boundary map representation for Schwarz transmission conditions. In~\cite{Goswami:2025}, Wang \textit{et al.} develop a hybrid framework that integrates PI-DeepONets~\cite{Mandl:2025} with finite element method (FEM) models via domain decomposition and O-SAM, focusing on a two-dimensional (2D) solid dynamics exemplar. The basic domain decomposition strategy is to assign PI-DeepONets to computationally demanding regions, while resolving the remainder of the computational domain via the classical FEM. Extensive numerical results are presented demonstrating that error margins below 1\% are achievable; however, the subdomain-local PI-DeepONets are extremely challenging to train, requiring up to $O(10^{6})$ epochs to converge.  The approach in \cite{Goswami:2025} was recently extended to non-overlapping domain decompositions, which employs displacement–traction information exchange across the subdomain interface \cite{wang2026nonoverlapping}.  For time-dependent problems, the approach employs a time-marching neural operator that incorporates Newmark time integration and uses the solution state from the previous time step together with interface data to predict the subdomain response.

In addition to DD-based approaches, localized reduced order modeling strategies have been developed in which multiple local reduced models are constructed over different regions of the solution or parameter space. Such methods seek to improve approximation accuracy by replacing a single global reduced basis with a collection of local models tailored to distinct dynamical regimes, as done in \cite{Geelen2022, Buhr2019}. In contrast, the present work employs a geometric decomposition of the computational domain and couples independently constructed reduced models through a Schwarz iteration framework.

The two most closely related works to the present manuscript are~\cite{Farcas:2023} and~\cite{Gkimisis:2025}, both of which propose Schwarz-like couplings involving non-intrusive OpInf ROMs~\cite{willcox2016opinf}. In~\cite{Farcas:2023}, Farcas \textit{et al.} develop an overlapping DD-based methodology for coupling subdomain-local OpInf ROMs by learning appropriate reduced operators responsible for the coupling, and demonstrates the method on a challenging three-dimensional (3D) combustion example. While this coupling framework is not presented in the context of SAM, it is effectively equivalent to applying a single iteration of the parallel (additive)~\cite{Gander:2008, wentland2024Schwarz} version of O-SAM, rather than iterating to convergence as done within the classical Schwarz framework considered herein. Another recent related work is that of Gkimisis \textit{et al.}~\cite{Gkimisis:2025}, which presents a fully non-intrusive hybrid ROM-FOM coupling framework, in which OpInf ROMs are coupled to ``sparse FOMs" (sFOMs) via an overlapping Schwarz iteration procedure. The sFOMs are inferred offline using training data collected from an existing high-fidelity simulation code by learning the numerical stencil of an adjacency-based FOM corresponding to the high-fidelity system dynamics. The primary motivation for this is to enable the creation of ROM-FOM couplings without access to the FOM code. 

Finally, the recent pre-print~\cite{Rodriguez:2025} from our group presents the first application, to the authors’ knowledge, of the non-overlapping SAM for the coupling of subdomain-local, non-intrusive OpInf ROMs with each other and with subdomain-local FOMs. Both Robin-Robin and alternating Dirichlet-Neumann transmission conditions are explored within this coupling framework, with the former yielding improved convergence and accuracy relative to the latter.  In our other recent pre-print \cite{sambataro2026rolerelaxationaccelerationnonoverlapping}, we demonstrate that convergence of the Dirichlet-Neumann non-overlapping SAM can be improved significantly through the use of Aitken \cite{Deparis2004} or Anderson acceleration \cite{Walker:2011}.

\subsection{Contributions, differentiating features and organization} \label{sec:contribs}

While the methodology described herein shares some similarities with the work of Farcas \textit{et al.}~\cite{Farcas:2023} and Gkimisis \textit{et al.}~\cite{Gkimisis:2025}, both of which propose overlapping Schwarz-like formulations for the DD-based coupling subdomain-local OpInf models, there are some important distinctions. Whereas~\cite{Farcas:2023} and~\cite{Gkimisis:2025} limit their attention to linear and quadratic OpInf ROMs, we also consider a higher-order cubic OpInf ROM, or ``COpInf." Moreover, neither~\cite{Farcas:2023} nor~\cite{Gkimisis:2025} describe and demonstrate a use case in which different time-integrators with different time steps are coupled together using SAM, and the former work does not consider FOM-OpInf couplings. As explained in Section~\ref{sec:past_work}, the method developed in Farcas \textit{et al.}~\cite{Farcas:2023} can be viewed as an additive Schwarz-like coupling framework in which each set of subdomain problems is solved once, rather than by iterating to convergence. We demonstrate in Section~\ref{sec:tension-specimen} that such an approach is inadequate for certain classes of nonlinear problems. Whereas the hybrid FOM-OpInf couplings developed herein require access to the FOM code in a minimally-intrusive way, the work of Gkimisis \textit{et al.}~\cite{Gkimisis:2025} takes the approach of using overlapping SAM to couple subdomain-local OpInf ROMs with sFOMs, sparse FOMs learned offline from available simulation data, as described in more detail in Section~\ref{sec:past_work}. One advantage of this method is that it yields a coupling framework that is fully non-intrusive even when performing OpInf-FOM coupling. However, learning the sFOM requires the numerical solution of a complex optimization problem, which can be infeasible when the sFOM is posed on a complex, unstructured 3D geometry. This issue is compounded for multi-query scenarios in which the geometry, mesh and/or material parameters may change between simulations, requiring re-learning of the sFOM. Additionally, the learned sFOM in~\cite{Gkimisis:2025} can introduce an additional and difficult to control source of error into the overall approximation pipeline. 

It is worth noting that the use case targeted by our SAM-based hybrid approach is similar to that of~\cite{Gkimisis:2025}: in a complex multiscale simulation, our strategy is to assign the subdomain containing the more complex, nonlinear dynamics to a FOM
%, as done in~\cite{Gkimisis:2025}, 
so as to maximize the accuracy of the coupled model. In contrast, other works (e.g., \cite{Goswami:2025, parish2024embedded}) take the approach of assigning more 
%refined 
dynamically complex subdomains to a (NN-based) ROM, towards maximizing improvements in online CPU time. Further, the OpInf ROMs proposed herein are substantially easier to train than the PI-DeepONets considered in~\cite{Goswami:2025}.

It is also worth mentioning that the domain decompositions we employ are physically motivated and mimic practical problems of interest in solid mechanics.   It is thus common in our applications to have 2--3 subdomain, and we do not envision having to couple more than 5--10 subdomains in general.  
%and in practical problems of interest to us are always physically-motivated; 
It follows that we do not target the use case where a complex domain is ``tiled" with hundreds of pre-trained subdomains, as done in~\cite{Wang:2022, Chung:2024, choi2025definingfoundationmodelscomputational}.

While the present work is an extension of some of our earlier research on SAM for FOM-FOM coupling~\cite{Mota:2017, Mota:2022} and contact~\cite{Mota:2025} in solid mechanics, intrusive ROM-ROM/ROM-FOM coupling~\cite{wentland2024Schwarz, Barnett:2022Schwarz}, and non-intrusive OpInf-OpInf/OpInf-FOM coupling~\cite{Moore:2024, Rodriguez:2025}, it contains a number of new contributions. This paper is not only one of the earliest to investigate SAM-based coupling of subdomain-local OpInf ROMs, it is also the first to apply this approach to realistic 3D nonlinear problems in solid mechanics. We feature problems posed on  geometries discretized with varying mesh resolutions and element types, and explore the integration of disparate time integration schemes, that may have different time-steps, into the coupling scheme. We demonstrate that it is possible to obtain accurate and efficient hybrid models for \textit{fully nonlinear} solid dynamics problems by coupling together subdomain-local OpInf ROMs based on \textit{global polynomial approximations} of the governing PDEs. Additionally, we propose several novel strategies for improving the efficiency and robustness of O-SAM, including boundary POD bases to reduce the size of boundary operators and methods for optimizing regularization parameters in subdomain-local OpInf ROMs within the Schwarz coupling framework. %\todo{Say something about cost analysis, if we have that in the paper?}

The remainder of this paper is organized as follows. In Section~\ref{sec:solid_mechanics}, we describe our model solid dynamics problem and its monolithic finite element discretization. Section~\ref{sec:opinf} overviews the OpInf approach to building non-intrusive ROMs, assuming a global cubic functional form for the learned model. In Section~\ref{sec:schwarz}, we describe our overlapping Schwarz algorithm and how it can be used to create hybrid FOM-OpInf and OpInf-OpInf models for generic solid dynamics applications. %\todo{If we include some variant of Section~\ref{sec:analysis}, describe it here.} 
Section~\ref{sec:results} evaluates the accuracy and efficiency of the proposed approach on four solid mechanics test cases, in which subdomains are discretized using disparate meshes/element types and/or advanced forward in time using different time-integrators with disparate time steps. Conclusions are offered in Section~\ref{sec:conc}.

%% file: 02-solid-mechanics.tex
\section{Solid mechanics problem formulation and monolithic finite element discretization} \label{sec:solid_mechanics}

%\ikt{I think the FOM description should be its own section, rather than a part of the intro.  I view the purpose of the intro as primarily reviewing past work and outlining contributions.}

%\todo{Consider changing presentation to start with variational form rather than strong form.  Alejandro likes, but not necessary.}

Consider the Euler-Lagrange equations for a general solid mechanics problem:
\begin{equation}
    \label{eq:dynamic_elasticity_pde}
    \nabla\cdot ~P + \rho ~R = \rho \ddot{~\varphi} \quad \text{in } I\times\Omega.
\end{equation}
Here, $\Omega \in \mathbb{R}^d$ for $d \in \{1, 2, 3\}$ is an open and bounded domain, $I \coloneqq \{t \in [\timeInit, \timeFinal]\}$ is a closed interval with $0 \leq \timeInit < \timeFinal$, and $~x = ~\varphi(t, ~X): I \times\Omega \rightarrow \mathbb{R}^d$ describes the mapping from the reference configuration $~X \in \Omega$ to the current configuration. The ``double-dot" notation $\ddot{~\varphi}$ indicates a second derivative in time. The symbol $\nabla\cdot$ represents the row-wise divergence, $~P:=~P(~\varphi)$ denotes the first Piola-Kirchoff stress tensor, $~R$ is the specific body force, and $\rho$ is the mass density defined in the reference configuration. We note that $~P$ is described by a constitutive model which is generally nonlinear in the current state $~\varphi$. Specifically, $~P$ is defined as the derivative of the Helmholtz free-energy density, denoted by $A(~F)$, with respect to the deformation gradient $~F$, that is, $~P := \frac{\partial A}{\partial ~F}$. In the present work, we restrict attention to hyperelastic material models, in which $A$ is not path-dependent (i.e., it does not depend on a collection of internal variables). For a concrete example of $A(~F)$ and $~P$, the reader is referred to Appendix A, where expressions for these variables are given for the linear elastic, Saint Venant--Kirchhoff \cite{holzapfel2000} and Neohookean \cite{Mota:2011} material models.

To ensure well-posedness, the problem \eqref{eq:dynamic_elasticity_pde} is subject to initial and boundary conditions (BCs),
% \begin{equation}
%     \label{eqn:dynamic_elasticity_bcs}
%     \begin{aligned}
%     ~{\varphi}(t_0,~X) &= ~X_0 && \text{in } \Omega \\
%     \dot{~{\varphi}}(t_0, ~X) &= ~v_0 && \text{in } \Omega \\
%     ~{\varphi}(t, ~X) &= \boldsymbol{\chi} && \text{on } \partial_{\varphi} \Omega \times I \\
%     ~P ~N &= ~T && \text{on } \partial_{T} \Omega \times I
%     \end{aligned}
% \end{equation}
%\begin{align}\label{eq:dynamic_elasticity_bcs}
%    ~{\varphi}(t_0,~X) &= ~X_0 \quad \text{in } \Omega, \qquad\qquad \dot{~{\varphi}}(t_0, ~X) = ~v_0 \quad \text{in } \Omega \\
%    ~{\varphi}(t, ~X) &= \boldsymbol{\chi} \quad \text{on } I \times \partial_{\varphi}\Omega \qquad\qquad 
%    ~P ~n = ~T \quad \text{on } I \times \partial_{T} \Omega
%\end{align}
\begin{equation} \label{eq:dynamic_elasticity_bcs}
    \begin{array}{ll}
    ~{\varphi}(t_0,~X) = ~x_0 \quad \text{ in } \Omega, & \quad\dot{~{\varphi}}(t_0, ~X) = ~v_0 \quad \text{ in } \Omega, \\
      ~{\varphi}(t, ~X) = \boldsymbol{\chi} \quad \text{ on } I \times \partial_{\varphi}\Omega, &\quad ~P ~n = ~T  \quad\text{ on } I \times \partial_{T} \Omega, 
      \end{array}
\end{equation}
where the boundary $\partial \Omega = \domainBoundaryDirichlet   \cup \domainBoundaryNeumann$ is the union of Dirichlet $\domainBoundaryDirichlet$ and Neumann $\domainBoundaryNeumann$ parts, which satisfy the non-overlapping condition $\domainBoundaryDirichlet \cap \domainBoundaryNeumann= \emptyset$, and $~n$ denotes the unit vector normal to $\domainBoundaryNeumann$.  

In most practical applications, we wish to discretize the PDEs above in space using the classical Galerkin FEM. To this end, consider the function space $Q = C^2\big(I; W^2_1(\Omega)\cap C^0(\partial\Omega)\big)$. Let $\mathcal{V} \coloneqq \{~\Theta \in Q: ~\Theta = ~\chi \text{ on } I\times\partial_{\varphi}\Omega\}$ and $\mathcal{V}_0 \coloneqq \{ ~{\zeta} \in Q: ~\zeta=~0 \text{ on } (\partial_{\varphi} \Omega \times I) \cup (\Omega \times \partial I)$ denote the test and trial function spaces, respectively, with $~\varphi\in \mathcal{V}$ and $~\xi\in \mathcal{V}_0$. It is now straightforward to derive the weak variational form of \eqref{eq:dynamic_elasticity_pde}--\eqref{eq:dynamic_elasticity_bcs}: 
%The present work is concerned with the weak variational form of \eqref{eq:dynamic_elasticity_pde},\eqref{eq:dynamic_elasticity_bcs},
%\begin{equation}\label{eq:weak_form}
%    \int_I \left[ \int_{\Omega} \Big( \rho \ddot{~{\varphi}}\cdot ~{\xi} + ~P : \nabla ~\xi - \rho ~R\cdot ~{\xi} \Big) \, d \Omega  - \int_{\partial_T \Omega} (~P~n)\cdot ~{\xi}  \, d S  \right] \, dt = 0.
%\end{equation}
\begin{equation}\label{eq:weak_form}
     \int_{\Omega} \Big( \rho \ddot{~{\varphi}}\cdot ~{\xi} + ~P : \nabla ~\xi - \rho ~R\cdot ~{\xi} \Big) \, d \Omega  - \int_{\partial_T \Omega} ~T \cdot ~{\xi}  \, d S   = 0.
\end{equation}
%The Neumann boundary condition from \eqref{eq:dynamic_elasticity_bcs} 
%is implemented by inserting $~P ~n = ~T$ into the second integral in \eqref{eq:weak_form}. \adg{Is it not better to simply replace $~P~n$ with $~T$ in (3) and avoid this comment?}
%\ikt{Yes, I have fixed it.}
%where the solution $~\varphi\in \mathcal{V} \coloneqq \{~\Theta \in C^2\big(I; W^2_1(\Omega)\big): ~\Theta = ~\chi \text{ on } I\otimes\partial_{\varphi}\Omega\}$ is sought in a space satisfying the Neumann BC and $~\xi\in \mathcal{V}_0 \coloneqq \{ ~{\zeta} \in C^2\big(I; W_1^2(\Omega)\big): ~\zeta=~0 \text{ on } (\partial_{\varphi} \Omega \times I) \cup (\Omega \times \partial I)$ is a test function in the space $\mathcal{V}$ vanishing identically on the Neumann boundary $\partial_{\varphi}\Omega$ and on the interior at times $t \in \{t_0,t_f\}$. 
%\adg{Why is this integrated in time?  Are we doing weak-form in time as well?}  \ikt{I think this is just how they commonly write the weak form in solid mechanics, with the implication being that you will handle the time-discretization in some way.  We are NOT doing space-time -- that is actually a really important detail.  I would be fine getting rid of the integral over time, as I agree it could be confusing.  We can check if Alejandro has any objections to this.}
Semi-discretization of \eqref{eq:weak_form} in space using standard isoparametric finite elements yields a system of the form 
\begin{equation}\label{eq:matrix_vector_form}
    ~M \ddot{~\varphi}_h + ~K(~\varphi_h) = ~f(~T).
\end{equation}
%\ejp{f on the RHS will need to be a function of the BCs in the case where $\chi$ or $T$ changes.}
Here, $~\varphi_h \in \mathbb{R}^M$ is the discretized solution vector, the mass matrix $~M\in\mathbb{R}^{M\times M}$ represents the usual discretization of the $\rho$-weighted $L^2$ inner product $(\cdot,\cdot)$, the nonlinear function $~K(~\varphi_h)\in\mathbb{R}^M$ discretizes the weak-form divergence of $~P$, and  $~f(~T) \in\mathbb{R}^{M}$ captures the contribution from the external force term $\rho~R$ along with the traction BC $~T$. The non-negative integer $M \in \mathbb{N}$ denotes the number of degrees of freedom (DoFs) of the semi-discrete system \eqref{eq:matrix_vector_form}.

To complete the discrete problem formulation, \eqref{eq:matrix_vector_form} must be augmented with appropriate Dirichlet BCs.  The standard way to impose Dirichlet BCs in the FEM is to split the solution vector $~u$ into constrained and unconstrained DoFs, modifying \eqref{eq:matrix_vector_form} such that only the unconstrained DoFs are solved for. 
%Toward this effect, let us 
To see this, partition the solution vector as $~\varphi_h := \left( \begin{array}{cc} ~u^{\intercal} & ~u_c^{\intercal} \end{array} \right)^{\intercal} := \left( \begin{array}{cc} ~u^{\intercal} & ~\chi_h^{\intercal} \end{array} \right)^{\intercal}$, where $~u \in \mathbb{R}^{N}$ is the vector of unconstrained DoFs to be solved for, and $~u_c = ~\chi_h \in \mathbb{R}^{K}$ are discrete values of the Dirichlet data $~\chi$, where $N, K \in \mathbb{N}$ and $N + K = M$. Substituting this decomposition into \eqref{eq:matrix_vector_form}, a system of the form
\begin{equation}\label{eq:matrix_vector_form_with_dbcs}
    \tilde{~M} \ddot{~u} + \tilde{~K}(~u, ~\chi_h) = \tilde{~f}(~T, ~\chi_h, \ddot{~\chi}_h)
\end{equation}
is obtained for the unconstrained DoFs $~u$, where  $\tilde{~M} \in \mathbb{R}^{N \times N}$ and $\tilde{~K}, \tilde{~f} \in \mathbb{R}^{N}$.  

To create a fully-discrete model, \eqref{eq:matrix_vector_form_with_dbcs} must be discretized in time.  Herein, the semi-discrete system is time-advanced using the standard Newmark-$\beta$ time-integration scheme \cite{Mota2003}, a popular choice in solid mechanics due to its flexibility, rigorous convergence guarantees, and preservation of the symplectic structure of the underlying PDEs. More details about this scheme can be found in Section \ref{sec:schwarz}.
%\ikt{I made $M$ the total number of DoFs, $N$ the number of unconstrained DoFs and $K$ the number of constrainted DoFs.  I think we are not using $M$ and $K$ later in the paper, but please correct me if I'm mistaken.}

The functional forms of $\tilde{~K}(~u, ~\chi_h)$ and $\tilde{~f}(~T, ~\chi_h, \ddot{~\chi}_h)$ in \eqref{eq:matrix_vector_form_with_dbcs} depend on the nonlinearities present in the vector of internal forces, $~K(\cdot)$.  Consider the example of a simple linear elastic material (Appendix A.1), for which $~K(~\varphi_h):= ~K ~\varphi_h$. 
Observe that the semi-discrete system \eqref{eq:matrix_vector_form_with_dbcs} can be written as 
\begin{equation}\label{eq:matrix_vector_form_linear}
    \tilde{~M} \ddot{~u} + \tilde{~K}~u = ~B ~g, 
\end{equation}
where $\tilde{~M} :=  ~M_{[1:N, 1:N]}$, $\tilde{~K}:= ~K_{[1:N, 1:N]}$, and
\begin{equation} \label{eq:B}
~B:= \left( \begin{array}{ccc}~I_u, & -~M_{uc}, & -~K_{uc} \end{array} \right) \in \mathbb{R}^{N \times (N + 2K)}, 
\end{equation}
\begin{equation}\label{eq:g}
    ~g := \left( \begin{array}{ccc} ~f_{u}^{\intercal}, & \ddot{\boldsymbol{\chi}}_h^{\intercal}, & \boldsymbol{\chi}_h^{\intercal} \end{array}\right)^{\intercal} \in \mathbb{R}^{N + 2K}, 
\end{equation}
with $~M_{uc}:= ~M_{[1:N, N+1:N]}$, $~K_{uc}:= ~K_{[1:N, N+1:M]}$, $~f_u:=~f_{[1:N]}$, $~I_u$ denoting the $N \times N$ identity matrix. Here, the notation $~M_{[i,j]}$ denots the $(i,j)^{th}$ component of the matrix $~M$ (and similarly for $~K$), and $~f_{[i]}$ denotes the $i^{th}$ component of the vector $~f$. %\ikt{Should subscript bracket notation be defined here?} 
%\adg{I think yes.  Also, I 
% \adg{I would recommend putting the brackets in-line instead of subscripted, for improved readability.} \ikt{I had that originally but thought it looks kind of ugly.}

In the present work, we focus our attention on developing ROMs for approximations of \eqref{eq:matrix_vector_form} which have a globally polynomial (at most, cubic) structure, i.e., systems of the form 
\begin{equation}
\label{eq:quadratic_model}
    ~M\ddot{~\varphi}_h + ~K ~\varphi_h + ~H~\varphi_h^{\otimes 2} + ~C~\varphi_h^{\otimes 3} = ~f(~T).
\end{equation}
Here, $~\varphi_h^{\otimes{k}}$ is a Kronecker product (i.e., matricized tensor product) of $~\varphi_h$ with itself $k$ times, and the matrices $~M,~K\in\mathbb{R}^{M \times M}$, $~H\in\mathbb{R}^{M\times M^2}$, and $~C\in\mathbb{R}^{M\times M^3}$ arise via the discretization of differential operators in the governing equations\footnote{Since $~\varphi_h^{\otimes k}$ is symmetric, it is computed in practice using a compressed vector representation of length $(1/k!)M(M+1)...(M+k-1)$.  Similarly, the symmetric operators $~H\in\mathbb{R}^{M\times M^2}$ and $~C\in\mathbb{R}^{M \times M^3}$ are parameterized in terms of only their independent degrees of freedom. 
% \ikt{Consider changing $N_k$ to something else since Eric uses $N_k$ in the Schwarz section.}
}, or by performing a Taylor expansion of \eqref{eq:matrix_vector_form}. 
% \crw{Is it necessary to have additional $N_2$ and $N_3$ terms? The fact that these don't strict need to be $N^2$ and $N^3$ in practical implementations doesn't seem particularly important, and just adds additional dimensions.}  \ikt{I agree with Chris.  Anthony, if you agree, can I ask you to please simplify?} \adg{The problem with not saying anything about this is that the OpInf really should be over only $N_i$ DoFs, not $N^i$ DoFs, otherwise you are not respecting symmetry.  I have tried to communicate this without the additional notation.}
For example, in the case of linear elasticity, $~M,~K$ represent the usual mass and stiffness matrices discretizing the $L^2$ inner product $(\cdot,\cdot)$ and weak-form Laplace operator $(\nabla\cdot,\nabla\cdot)$.  As shown in Appendix A.2, the PDEs \eqref{eq:dynamic_elasticity_pde}--\eqref{eq:dynamic_elasticity_bcs} reduce to a cubic functional form \eqref{eq:quadratic_model} when a Saint Venant--Kichhoff material model \cite{holzapfel2000} is specified. While, for a general hyperelastic material, the governing PDEs will not reduce to \eqref{eq:quadratic_model}, we show in Section \ref{sec:results} that the system \eqref{eq:quadratic_model} can be effective as a surrogate model when combined with Schwarz coupling even when the governing equations under consideration are fully nonlinear.  %This is partly attributed to the fact that \eqref{eq:quadratic_model} can be seen as arising from a Taylor expansion of \eqref{eq:matrix_vector_form}.  

As for the linear case, to apply Dirichlet BCs within \eqref{eq:quadratic_model}, we must substitute the decomposition  $~\varphi_h := \left( \begin{array}{cc} ~u^{\intercal}, & ~\chi_h^{\intercal} \end{array} \right)^{\intercal}$ into \eqref{eq:quadratic_model} and bring the terms corresponding to the Dirichlet boundary condition DoFs to the right-hand side.  Doing so gives rise to a system of the form 
\begin{equation}
\label{eq:quadratic_model_with_bcs}
    \tilde{~M}\ddot{~u} + \tilde{~K}(~\chi_h) ~u + \tilde{~H}(~\chi_h)~u^{\otimes 2} + \tilde{~C}~u^{\otimes 3} = ~f_u - ~M_{uc} \ddot{~\chi}_h - ~K_{uc} ~\chi_h - ~h(~H, ~C, ~\chi_h),
\end{equation}
for matrix $\tilde{~K}(~\chi_h) \in \mathbb{R}^{N \times N}$, vector $~h\in \mathbb{R}^{N}$, and appropriately-sized tensors $\tilde{~H}(~\chi_h)$ and $\tilde{~C}$. The linear and quadratic operators are functions of the boundary data $~\chi_h$, as they include contributions from cross terms of the form $~u \otimes ~\chi_h$ and $~u \otimes ~\chi_h \otimes ~\chi_h$ that arise when substituting the decomposed $~\varphi_h$ vector into \eqref{eq:quadratic_model}.   
%The tensors $~H_{ucc}$, $~G_{uuuu}$ and  $~G_{uccc}$ are defined analogously to $~K_{uu}$ and $~K_{uc}$, whereas $\tilde{~H}(~u_c)$ and $\tilde{~G}(~u_c)$ include also the cross terms $~u_u ~u_c$, $~u_u ~u_u ~u_c$, etc. arising from computing Kronecker products involving $~u_u$ and $~u_c$. 
%\ikt{Is the above too imprecise?  I don't want to give expressions for each of the matrices/tensors b/c they will be too complicated.  Maybe some of the multiple subscripts should be dropped also...}

In defining the functional form of our OpInf ROMs in Sections \ref{sec:opinf} and \ref{sec:schwarz}, we will approximate the right-hand-side of \eqref{eq:quadratic_model_with_bcs} in several ways. First, because none of the numerical examples considered herein have a body source term or inhomogeneous traction boundary conditions, we will assume from this point forward that $~f_u = ~0$ in \eqref{eq:quadratic_model_with_bcs}. We will additionally neglect the $\ddot{~\chi}_h$ term in \eqref{eq:g}; this term drops out if the mass matrix is lumped\footnote{Many production codes, such as Sandia's {\tt SIERRA/SM} code \cite{sierrasm}, utilize a lumped matrix for both implicit and explicit solves.} and is expected to be small in the general case. Finally, we will account for the influence of the $~h(~H, ~C, ~\chi_h)$ term in \eqref{eq:quadratic_model_with_bcs}
through a new matrix $\tilde{~B}\in\mathbb{R}^{N\times K}$, discussed in more detail in Remark \ref{remark:bc_approx}. %\adg{Instead of ``neglect the ... term'', we should say that we account for its influence through $\tilde{~B}$ and reference the remark.}. \ikt{I agree with this and have changed the text.} 
With these assumptions, our approximate semi-discrete FOM model on which the algebraic structure of our OpInf ROMs is based (see Sections \ref{sec:opinf} and \ref{sec:schwarz}) takes the form: 
\begin{equation}
\label{eq:quadratic_model_with_bcs_approx}
     \tilde{~M}\ddot{~u} + \tilde{~K}(~\chi_h) ~u + \tilde{~H}(~\chi_h)~u^{\otimes 2} + \tilde{~C}~u^{\otimes 3} =  \tilde{~B} \tilde{~g},
\end{equation}
where %\adg{I would completely cut (12). It undermines our argument because $\tilde{~B}$ should catch other contributions too.} \ikt{Good catch, I have removed this.}
\begin{equation}\label{eq:gtilde}
    \tilde{~g} :=  \boldsymbol{\chi}_h \in \mathbb{R}^{K}.   
\end{equation}
As we show in Section \ref{sec:results}, subdomain-local OpInf ROMs based on \eqref{eq:quadratic_model_with_bcs_approx}, when coupled to other subdomain-local OpInf ROMs and/or subdomain-local FOMs, can yield tremendously accurate approximations despite the fact that a number of assumptions and simplifications have been made to go from \eqref{eq:matrix_vector_form} to \eqref{eq:quadratic_model_with_bcs_approx}.  

\begin{remark} \label{remark:bc_approx}
  The approximation of \eqref{eq:quadratic_model_with_bcs} by  \eqref{eq:quadratic_model_with_bcs_approx} should not be interpreted as neglecting all nonlinear dependence on the interface data. After omitting the first two terms on the right-hand side of \eqref{eq:quadratic_model_with_bcs}, the remaining contribution consists of a term that is linear in the interface data together with the nonlinear interaction term $~h(~H,~C,\chi_h)$. As discussed in more detail in Section \ref{sec:opinf}, the ROM ansatz in \eqref{eq:quadratic_model_with_bcs_approx} replaces the combined effect of these contributions with a learned linear boundary operator $\tilde{~B}\chi_h$. One may view this step as approximating
\begin{equation}
 ~h(~H,~C,\chi_h)\approx ~A\chi_h,
\end{equation}
for some operator $~A$, yielding $\tilde{~B}=-~K_{uc}-~A$.

More broadly, it is important to note that the polynomial OpInf ROMs considered herein are themselves approximations of the underlying FOMs, which are based on fully nonlinear equations. Hence, for the fully nonlinear problems considered in this work, neither the reduced operators nor the boundary representations are expected to exactly reproduce the governing equations. Rather, the objective is to identify a compact reduced model that accurately reproduces the observed system dynamics. 
Because $\tilde{~B}$ 
% is
will be inferred simultaneously with the remaining reduced operators, the learned boundary operator may absorb part of the effective influence of the neglected interaction terms over the training trajectories.

More expressive boundary representations could also be considered, for example, by augmenting the reduced model with quadratic or higher-order functions of the interface data. However, the numerical results presented in Section \ref{sec:results} suggest that the linear boundary representation adopted here is sufficiently expressive for the problems considered, as evidenced by the small solution errors and robust Schwarz convergence observed across all test cases. We therefore do not consider the investigation of enriched boundary operators in this work.
\end{remark}

\begin{remark} \label{remark:bc_approx_penalty}
There is an alternate justification for the approximation of \eqref{eq:quadratic_model_with_bcs} by \eqref{eq:quadratic_model_with_bcs_approx}.  The boundary forcing term appearing in \eqref{eq:quadratic_model_with_bcs_approx} is also consistent with formulations in which Dirichlet boundary conditions are imposed weakly using a penalty method or Nitsche's method \cite{hansbo2002unfitted}. For example, a penalty formulation augments the variational statement with a term of the form
\begin{equation}
    \frac{\gamma}{h}\int_{\Gamma}(~u-~g)~v,dS,
\end{equation}
where $\Gamma$ is a given boundary, $~g$ denotes prescribed boundary data, $h$ is a characteristic mesh size, $~u$ and $~v$ are test and trial functions, respectively, and $\gamma$ is a penalty parameter. After spatial discretization and projection onto a reduced basis, the resulting reduced equations contain additional forcing terms that depend linearly on the boundary data. Consequently, the linear boundary operator appearing in \eqref{eq:quadratic_model_with_bcs_approx} may also be interpreted as a reduced representation of boundary forcing effects arising from weak enforcement of interface conditions.

We emphasize that the present work does not explicitly derive the ROM from a penalty or Nitsche formulation. Rather, this observation provides additional intuition for why a linear dependence on the interface data constitutes a reasonable and commonly encountered reduced-order modeling ansatz.

\end{remark}

\begin{remark} \label{remark:nonlinear_fom}
The polynomial approximation \eqref{eq:quadratic_model_with_bcs_approx} to \eqref{eq:dynamic_elasticity_sam_bcs} is introduced only for the purpose of motivating the cubic OpInf ROMs described in Section \ref{sec:opinf} and considered herein.  As will be made clear in Section \ref{sec:opinf-fom}, our FOM-OpInf couplings create hybrid models in which fully nonlinear subdomain-local FOMs based on \eqref{eq:dynamic_elasticity_sam_bcs} are coupled to polynomial subdomain-local ROMs based on \eqref{eq:quadratic_model_with_bcs_approx}.       
\end{remark}

%% file: 03-opinf.tex
\section{Model reduction via non-intrusive operator inference (OpInf)} \label{sec:opinf}

Many effective model reduction strategies are intrusive to the underlying FOM discretization. Galerkin ROMs built with a data-driven technique such as POD~\cite{Holmes:1996, Kunisch2001POD,Sirovich:1987}, for example, require the ability to project the discrete operators governing the FOM onto the span of the POD basis, which %in turn 
requires the user to directly modify the FOM implementation. Conversely, many FOMs corresponding to practical applications are proprietary, built on legacy software, or require specialized expertise to simulate, in which case the underlying FOM must be treated as inaccessible. This has motivated the development of non-intrusive model reduction strategies such as OpInf, which will now be described. 

OpInf is a projection-inspired method for model reduction, originally introduced by Peherstorfer and Willcox~\cite{willcox2016opinf} as a non-intrusive alternative to standard Galerkin ROMs. The algorithm proceeds in two major stages. First, a basis for the solution is built from FOM snapshot data similarly to the intrusive Galerkin case. However, instead of projecting the FOM onto the span of this basis, the second stage infers a data-driven approximation to the necessary projected operators from (projections of) the same FOM snapshot data, allowing for a fully non-intrusive ROM. The remainder of this section describes the particular OpInf strategy that is employed in the present work.

% \adg{the following paragraph is (almost) lifted from Ian's paper. Not sure it belongs here.}

% Continuing the example of the heat equation from~\eqref{eq:heat_pde}, a spatially discretized monolithic FOM for the heat equation typically appears in the following form after a boundary lift:
% \begin{equation}
% \label{eq:heat_eq_discretized}
%     ~M\dot{~x} = ~K ~x + ~B ~g,
% \end{equation}
% where $\mathbf{x} \in \mathbb{R}^{N}$ is a discretized vector corresponding to the unconstrained state degrees of freedom (DoFs), and $~g \in \mathbb{R}^{m}$ discretizes the Dirichlet boundary condition. The matrices $~M,~K \in \mathbb{R}^{N\times N}$ come from discretizing the $L^2$ inner product $(\cdot,\cdot)$ and weak-form Laplace operator $(\nabla\cdot, \nabla\cdot)$, while $~B \in \mathbb{R}^{N \times m}$ deals with boundary-to-interior coupling. The full state representation $~u \in \mathbb{R}^{N + m}$ for all DoFs is obtained by augmenting the unconstrained solution $~x$ with the known boundary condition $~g$.

% \subsection{Operator Inference Preliminaries} \label{sec:opinf}

 % Similarly to other Galerkin projection ROMs, the method begins by constructing a reduced basis from data. OpInf then diverges from the standard method in the construction of the reduced operators, which are estimated using data through a regression problem. 

\subsection{Stage 1: proper orthogonal decomposition} \label{sec:pod}

% \ikt{I don't think we want to use $~x$ for the state because that is typically used to denote the position in solid mechanics. The displacement is usually called $~u$. We could alternatively use $~q$.}

%\ikt{We have a notation problem here. $~\Psi$ was used earlier in the solid mechanics formulation. I will ask Alejandro if there is another reasonable variable to use there, as I'd like to stick with $~\Phi$ and $~\Psi$ for POD/ROM.}

The first stage of OpInf constructs a reduced basis, commonly with the standard POD algorithm~\cite{Holmes:1996, Kunisch2001POD,Sirovich:1987}.
% , though one might also use other methods, e.g., the reduced basis method with a greedy algorithm~\cite{RozzaGreedy2014}. To perform POD, we require snapshots from some data source or FOM. %\irinanote{New comment 8/7: maybe define $\mathbf{u}$ mathematically?} \irm{sufficient?}
To illustrate this procedure, suppose that the monolithic problem~\eqref{eq:quadratic_model} has been solved in time for $\numTime+1$ separate states $0 = T_0 < T_1 < \dots < T_{\tau} = \timeFinal$, yielding a collection of snapshots of the unconstrained state, 
\begin{equation} 
~U = [~u(T_0), ~u(T_1), \dots, \textcolor{black}{~u(T_r)}] \in \mathbb{R}^{N 
\times (\numTime+1)}.
 \label{eq:snapshots}
\end{equation}
%\ikt{There is another notation problem: we use $d$ earlier to denote the number of spatial dimensions. (BTW, what is the syntax used for lower vs. upper case letters denoting natural number?). Also we need to make the time notation consistent with the solid mechanics formulation section. There, time goes from $t_0$ to $t_f$. I think $\tau$ should be renamed $f$. }
% \irinanote{I think the size of $\mathbf{X}$ should be $n \times \tau$, based on how you index the time steps, no? Also, why are $n$ and $N_{\tau}$ bold in the above equation?}
The governing assumption in POD is that the solution space to the FOM \eqref{eq:quadratic_model} has low rank relative to the DoF space $\mathbb{R}^N$. Correspondingly, the snapshot matrix $~U$ possesses a low-rank structure, i.e., rank$(~U) \ll N$, which can be exploited for model reduction. Considering the (thin) singular value decomposition (SVD) $~U = ~{\Phi} ~{\Sigma} ~\Psi^{\intercal}$, the column-orthonormal matrix $~{\Phi} \in \RRStar{N}{\text{rank}(~U)}$ then spans the range of the discrete dynamics exhibited in $~U$; here $\RRStar{N}{K}$ denotes the space of $N \times K$ orthonormal matrices (i.e., $\RRStar{N}{K} \defEq \{ \mathbf{X} \in \RR{N \times K} \ | \ \mathbf{X}^T \mathbf{X} = \mathbf{I} \})$. %\crw{$m$ is used later as an index. Is it useful to define a separate rank of the solution space as well as a ROM basis dimension? Seems like the former doesn't add much.} \ikt{I got rid of it.}
Moreover, the first $r\leq \text{rank}(~U)$ columns of $~{\Phi}$ form a variance-maximizing basis, denoted $~{\Phi}_r$ which is $\ell^2$-optimal in the sense that 
% \crw{Is this really stating optimality? Shouldn't this instead read that $~\Phi_r$ is the solution to an argmin of the norm?} \adg{fair enough.  I changed it slightly}
% \begin{equation} \label{eq:pod_cutoff}
%     \min_{~V\in\mathbb{R}^{N\times r}} \left\|~U-~V~V^\intercal~U\right\|_F^2 = \left\|~U-~{\Phi}_r~{\Phi}_r^\intercal~U\right\|_F^2 = \sum_{t=0}^{\numTime} \left\|~u(t_i)-~{\Phi}_r~{\Phi}_r^\intercal~u(t_i)\right\|_2^2 = \sum_{i=r+1}^{\text{rank}(~U)}\sigma_i^2,
% \end{equation}
\begin{equation} \label{eq:pod_cutoff}
    \min_{~V\in\mathbb{R}^{N\times r}} \left\|~U-~V~V^\intercal~U\right\|_F^2 = \left\|~U-~{\Phi}_r~{\Phi}_r^\intercal~U\right\|_F^2 = \sum_{i=r+1}^{\text{rank}(~U)}\sigma_i^2,
\end{equation}
where $\sigma_i$ denotes the $i^{\mathrm{th}}$ diagonal entry in $~{\Sigma}$. The basis $~{\Phi}_r$ is employed in all present cases and will be referred to as the POD basis of rank $r$.
%\adg{We can talk about orthonormality in other IPs here, but I don't think we ever use it for the experiments, so I would recommend leaving it out.}
%\ikt{I think it's fine to skip discussion of other IPs.}

%\adg{I swapped $\Phi$ and $\Psi$, because the other way is nonstandard and weird to me alphabetically. We can fight about it if you want but I am not likely to back down.}
%\ikt{I'm fine with this, but note my commenta above about duplicate notation.}

\begin{remark} \label{remark1}
  The discussion above suggests that only the primary solution field, in this case, the displacement field $~u$, is used to build the POD basis $~\Phi_r$. While this is the present approach, we note that it is possible to augment the snapshot set $~U$ in \eqref{eq:snapshots} with other data, e.g., snapshots of the velocity and/or acceleration fields, as discussed in~\cite{Barnett:2022}.
\end{remark}

\begin{remark} \label{remark2}
    While we limit this discussion to affine POD bases, we point out that the ROM latent space can also be comprised of nonlinear representations, e.g., quadratic bases \cite{Barnett:2022, Geelen:2023}, convolutional autoencoders \cite{Lee:2020}, or other neural networks \cite{Barnett:2023}.  In particular, our SAM-based coupling framework can be used to couple OpInf ROMs with latent spaces spanned by either linear or nonlinear bases.  
\end{remark}

\subsection{Stage 2: least-squares regression}
\label{subsec:Opinf_Details}

%\ikt{In my mind, there is some novelty in what we are doing with the BC enforcement. I don't think that comes out here. Is there a way to emphasize it more either here or in the Schwarz section?}

%\ikt{We totally don't mention regularization when this is a **critical** part of the approach and Eric has a novel-ish approach to find the optimal regularization parameter within the Schwarz context. This needs to be discussed here and/or in the Schwarz section.}

Following the calculation of the POD basis $~{\Phi}_r$, it is necessary to construct a non-intrusive approximation to the usual projection-based Galerkin ROM, which forms the core of the OpInf routine. 
Since we know that the algebraic structure of a polynomial FOM is preserved under projection, the POD ROM corresponding to \eqref{eq:quadratic_model_with_bcs_approx} is also cubic and takes the form\footnote{Postulation of \eqref{eq:quadratic_ROM} assumes that the Dirichlet degrees of freedom have been explicitly removed, unlike the FOM analogues of these matrices in \eqref{eq:quadratic_model_with_bcs_approx}.}
% \footnote{Note that, in postulating \eqref{eq:quadratic_ROM}, we are assuming $\hat{K}$ and $\hat{H}$ are independent of the Dirichlet data, unlike the FOM analogs of these matrices in \eqref{eq:quadratic_model_with_bcs_approx}. This is not unreasonable to assume, as $\hat{K}$ and $\hat{H}$ are generic operators to-be-learned.  \ikt{Anthony, please check this / reword if needed.} } 
%Recall that the intrusive POD-ROM corresponding to the discrete equation \eqref{eq:quadratic_model_with_bcs_approx} is given by \adg{This is not quite correct in the sense that the K-hat and H-hat operators are a function of the boundary data in the derivation of the FOM analog of the functional form below.}  \ikt{Maybe add this as a remark?} 
\begin{equation}\label{eq:quadratic_ROM}
\hat{~M}\ddot{\hat{~u}} + \hat{~K}\hat{~u} + \hat{~H}\hat{~u}^{\otimes 2} + \hat{~C}\hat{~u}^{\otimes 3} = \hat{~B}\tilde{~g},
\end{equation}
where $\hat{~M}=~{\Phi}_r^\intercal~M~{\Phi}_r$, $\hat{~K}=~{\Phi}_r^\intercal~K~{\Phi}_r$, $\hat{~H} = ~\Phi^\intercal_r~H\lr{~\Phi_r\otimes~\Phi_r}$, $\hat{~C}= ~\Phi_r^\intercal~C\lr{~\Phi_r\otimes~\Phi_r\otimes~\Phi_r}$, and $\hat{~B}=~{\Phi}_r^\intercal ~B$ are projections of the FOM operators onto the span of the reduced basis $~{\Phi}_r$. The goal of OpInf is to infer approximations of these objects from the collected snapshot data $~U$. Assuming an $~M$-orthonormal basis $~{\Phi}_r$, or, non-equivalently, left multiplying by $\hat{~M}^{-1}$, the alternative model form referred to as the monolithic OpInf formulation can be written as
\begin{equation} \label{eq:opinf_rom}
    \ddot{\hat{~u}} + \bar{~K} \hat{~u} + \bar{~H}\hat{~u}^{\otimes 2} + \bar{~C}\hat{~u}^{\otimes 3} = \bar{~B} \tilde{~g}.
\end{equation}
Here, the low-dimensional operators $\bar{~K},\bar{~H},\bar{~B}$ are data-driven surrogates for their intrusive counterparts in \eqref{eq:quadratic_ROM}, which are assumed inaccessible\footnote{As before, the operators $\bar{~H}$ and $\bar{~C}$ are parameterized in terms of their compressed representation, so that a symmetric operator is guaranteed after inference.}.  These operators can be computed as the solution to a convex learning problem with a least-squares objective. More precisely, given projected state data $\bar{~U} = ~{\Phi}_r^\intercal~U$ and projected velocity information $D^2_t(\bar{~U}) \approx ~{\Phi}_r^\intercal\ddot{~U}$ computed with, e.g., a finite difference operator $D_t$, the OpInf learning problem is given as a minimization of the approximate residual: 
% \adg{can we go back to bars here instead of overlines?  It looks like the quantities KU, HU etc are single objects, which is not the case.}
% \begin{equation}\label{eq:quadratic_opinf}
% \argmin_{\overline{~K},\overline{~B}} \left\|D^2_t(\overline{~U})-\overline{~K}\overline{~U} - \overline{~H}\lr{\overline{~U}\ast\overline{~U}} - \overline{~G}\lr{\overline{~U}\ast\overline{~U}\ast\overline{~U}} - \overline{~B}~g~1^\intercal\right\|_F^2,
% \end{equation}
%\irm{I've added the norms here for consistency with what is currently (7) below.} \ejp{adding regularizaiton parameters} 
% \begin{equation}\label{eq:quadratic_opinf}
% \argmin_{\overline{~K},\overline{~H},\overline{~G},\overline{~B}} \left\|D^2_t(\overline{~U}) +\overline{~K}\overline{~U} + \overline{~H}\overline{~U}^{\ast 2} + \overline{~G}\overline{~U}^{\ast 3} - \overline{~B}~g~1^\intercal\right\|_F^2 + \regularizationParameter \sum_{\mathcal{O} \in \{\overline{~K},\overline{~H},\overline{~G},\overline{~B} \}  } \| \mathcal{O} \|_F 
% \end{equation}
\begin{equation}\label{eq:quadratic_opinf}
\argmin_{\bar{~K},\bar{~H},\bar{~C},\bar{~B}} \left\|D^2_t(\bar{~U}) +\bar{~K}\bar{~U} + \bar{~H}\bar{~U}^{\ast 2} + \bar{~C}\bar{~U}^{\ast 3} - \bar{~B}~{G}\right\|_F^2 + \regularizationParameter \big( \|\bar{~K}\|_F^2 + \|\bar{~H}\|_F^2 + \|\bar{~C}\|_F^2 + \|\bar{~B}\|_F^2 \big)
\end{equation}
where $\bar{~U}^{\ast k}\in\mathbb{R}^{r^k\times \numTime}$ is the column-wise Kronecker (i.e., Khatri--Rao) power of $~U$, $~{G}$ is a matrix whose $i^{\mathrm{th}}$ column is $\tilde{~g}(t_i)$, and $\regularizationParameter \in \mathbb{R}^+$ is a (scalar-valued) regularization parameter. The most basic, linear form of OpInf neglects the $\bar{~H}$ and $\bar{~C}$ terms, while quadratic OpInf (or QOpInf) only neglects the $\bar{~C}$ term.  Cubic OpInf (COpInf) retain the all terms in \eqref{eq:quadratic_opinf}.

\begin{remark} \label{remark3}
    We note that various OpInf works, e.g., \cite{McQuarrie2021combustion, qian:2023nonlinearopinf}, examine vector-valued regularization parameters such that each operator has its own regularization parameter. We do not consider this here as the resulting grid search for the optimal parameter combination scales exponentially with the number of parameters to infer.
\end{remark}

The advantage of this procedure is that the simulation of the approximate ROM \eqref{eq:opinf_rom} is non-intrusive, since the surrogate operators $\bar{~K},\bar{~H},\bar{~C},\bar{~B}$ are learned directly from data.  Moreover, under the assumptions that (i) the time-integration of the FOM \eqref{eq:quadratic_model} is convergent with decreasing step-size $\Delta t\to 0$, (ii) the data $\bar{~U}$ (and its Khatri--Rao powers) have full column rank,
% \irm{I think some parentheses would make this read easier, e.g. $(Bg)1^T$ or $Bg(1^T)$?} \adg{I don't think it matters since all the ops are associative. But if we must put something I would choose $B(g1^T)$.} \adg{Double check that this is actually all that is needed as far as rank requirements.}
% $~X^\intercal\otimes_k~I + ~1~g^\intercal\otimes_k~I$ are column-wise linearly independent,
and (iii) the time derivative approximations $D^2_t(~U)$ converge to their instantaneous equivalents with decreasing difference parameter, it follows that for $\lambda = 0$ and any $\varepsilon>0$ %\ejp{for $\regularizationParameter = 0$} 
, the learned approximations satisfy (c.f.~\cite{willcox2016opinf}) 
\begin{equation}
    \| \bar{~K} - \hat{~K} \|_F < \varepsilon, \qquad \| \bar{~H} - \hat{~H} \|_F < \varepsilon, \qquad
    \| \bar{~B} - \hat{~B} \|_F < \varepsilon, \qquad \| \bar{~C}-\hat{~C} \|_{F} < \varepsilon,
\end{equation}
for some $r \leq N$ and $\Delta t>0$, both depending on $\varepsilon$. Therefore, convergence of the learned operators to their intrusive counterparts is guaranteed in appropriate limits. Conversely, it is clear from the lack of commutativity between Galerkin projection and the flow maps corresponding to~\eqref{eq:quadratic_model} and \eqref{eq:quadratic_ROM} (see, e.g.,~\cite{gruber2025variationally}) that solutions to the OpInf problem \eqref{eq:quadratic_opinf} will never pre-asymptotically recover the corresponding intrusive operators appearing in \eqref{eq:quadratic_ROM}. This has motivated semi-intrusive methods such as re-projection~\cite{Peherstorfer:2020} and rank-aware snapshot collection~\cite{rosenberger2025exact} which are guaranteed to eliminate this troublesome ``closure error'' that limits ROM performance. Even without these techniques, OpInf remains a powerful method for building and deploying non-intrusive surrogate models in a variety of settings. The next section will discuss its use within the Schwarz coupling framework.

% While this overview has focused on the heat equation, we again stress that operator inference is general to low-order polynomial non-linearity, with results previously shown for a polynomial non-linearity of degree three in a one-dimensional (1D) nuclear reactor model in~\cite{willcox2016opinf}, for a 2D single injector combustion model in~\cite{McQuarrie2021combustion}, and for a 3D rotating detonation rocket engine in~\cite{Farcas:2023}, among others.

%% file: 04-schwarz.tex
\section{The Schwarz alternating method for heterogeneous model coupling} \label{sec:schwarz}

%\ikt{Describe the following:
%\begin{itemize}
%\item Schwarz algorithm.
%\item Describe OpInf-OpInf and OpInf-FOM coupling (in particular, how BCs are handled in OpInf ROM).
%\item Describe regularization parameter tuning.
%\item Describe expansion of BC data in its own basis to reduce cost.
%\item Stick to multiplicative Schwarz.  Compare to corresponding FOM-FOM coupling rather than single-domain FOM. \\
%\item \ejp{If it pans out, I think adding an error estimate based on the overlap region could be nice.} % \ikt{I agree.  I added an optional Section 4 for this, but it could go here as well.  The only other consideration is we may want to add results that have different overlap regions in this case.}
%\end{itemize}}

%\ikt{The following ingredients should be discussed:
%\begin{itemize}
% \item 	Projection operators
%\item 	How we handle time-stepping with different dts and mesh resolutions, along with figures.
%\item	Where the B operator term ``comes from” relative to the underlying finite element discretization.
%\item 	How we do training (top-down vs. bottom-up, though bottom-up is possible within the Schwarz framework)
%\item Reduction of boundary condition operator in the case there are a lot of DBC nodes.
%\item 	How we optimize the relaxation parameters
%\item Describe expansion of BC data in its own basis to reduce cost
%\item 	Highlight along the way what is novel.
%\item  Add algorithm explaining offline/online workflows.
%\end{itemize}
%}

The second ingredient necessary for the present approach is the Schwarz alternating method. SAM is a ``divide-and-conquer'' strategy which first decomposes the physical domain on which a given set of PDEs is posed into smaller subdomains.  Following this, the governing PDEs are solved by iterating between subdomain-local problems, with boundary information exchanged to ensure compatibility across subdomain interfaces (see Algorithm \ref{algo:sam-fom}). After overviewing O-SAM for FOM-FOM coupling in solid dynamics (Section \ref{sec:schwarz_fom}), we describe our extension of the approach to non-intrusive OpInf ROMs (Section \ref{sec:schwarz_rom}). We note that this section proposes several novel methodologies for improving the efficiency and robustness of O-SAM for couplings involving OpInf ROMs, in particular: (i) the introduction of boundary POD bases to reduce the size of the boundary operators $\bar{~B}_i$ in the case there are a large number of Dirichlet boundary nodes (Section \ref{sec:schwarz_opinf_offline}), and (ii) an efficient and robust algorithm for optimally selecting the regularization parameter in the subdomain-local OpInf least squares optimization problems (Section \ref{sec:regularization}).

\subsection{O-SAM for FOM-FOM coupling } \label{sec:schwarz_fom}

\begin{figure}[ht!]
    \centering
    \includegraphics[width=0.6\linewidth]{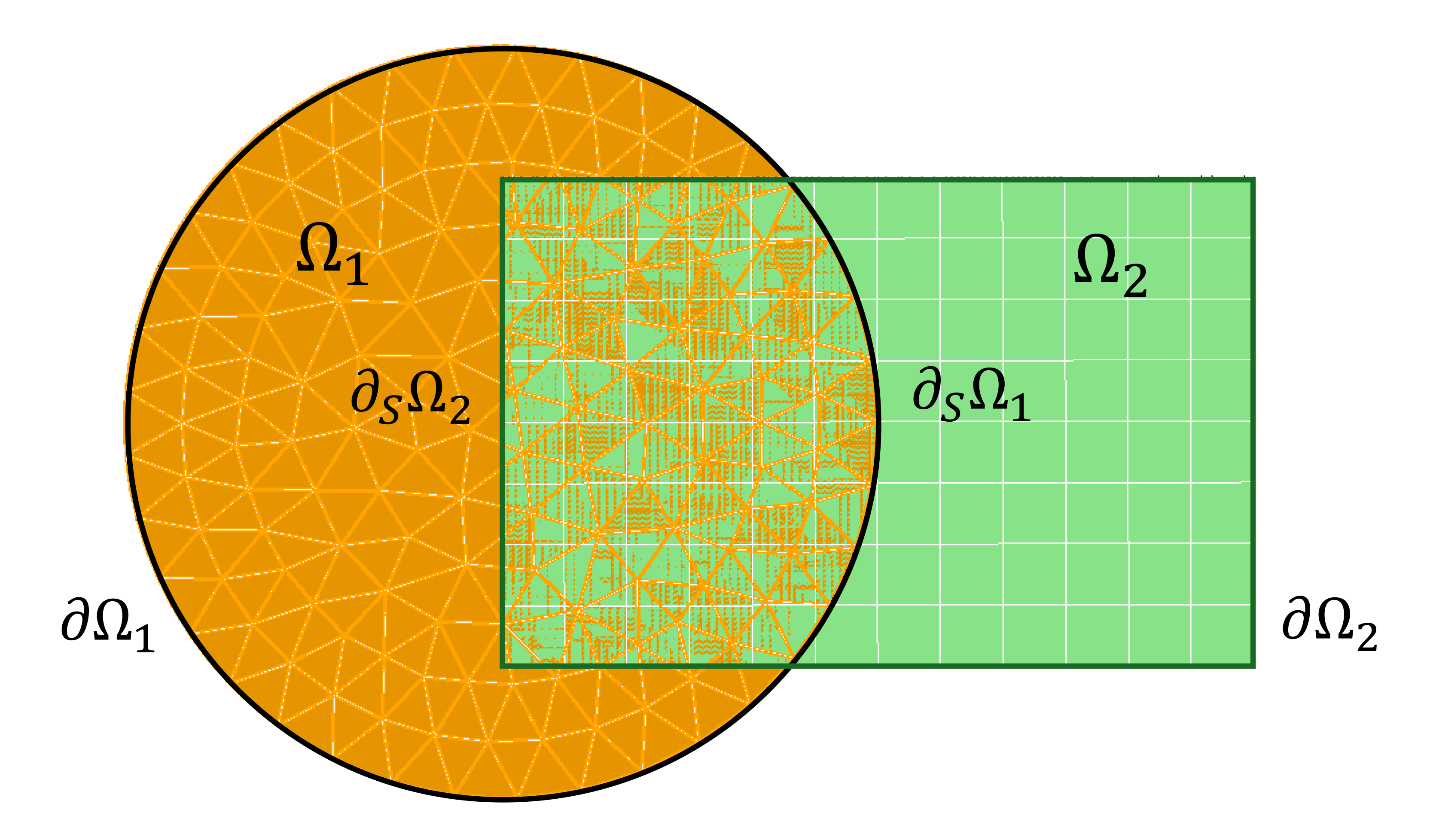}    
    \caption{Illustration showing an domain decomposition of a 2D domain $\Omega$ into two overlapping subdomains, $\Omega_1$ and $\Omega_2$, for the application of O-SAM.  } 
    \label{fig:schwarz_overlapping}
\end{figure}

Without loss of generality, consider the decomposition of $\Omega$ into two overlapping subdomains, $\Omega_1$ and $\Omega_2$, with boundaries $\partial \Omega_1$ and $\partial \Omega_2$, respectively, such that $\Omega_1 \cap \Omega_2 \ne \emptyset$, as shown in Figure \ref{fig:schwarz_overlapping}. The boundary of each subdomain $\domainBoundaryOneArg{i} = \domainBoundaryDirichletOneArg{i} \cup \domainBoundaryNeumannOneArg{i} \cup \domainBoundarySchwarzOneArg{i}$ is the union of the (exterior) Dirichlet $\domainBoundaryDirichletOneArg{i}$, the (exterior) Neumann $\domainBoundaryNeumannOneArg{i}$, and (interior) interface $\domainBoundarySchwarzOneArg{i}$ parts, which satisfy the non-overlapping condition $\domainBoundaryDirichletOneArg{i} \cap \domainBoundaryNeumannOneArg{i}\cap \domainBoundarySchwarzOneArg{i} =\emptyset$. The interior interface parts are defined as $\domainBoundarySchwarzOneArg{1} = \domainBoundaryOneArg{1} \cap \domainOneArg{2}$ and $\domainBoundarySchwarzOneArg{2} = \domainBoundaryOneArg{2} \cap \domainOneArg{1}.$ We refer the reader to Figure~\ref{fig:schwarz_overlapping} for a depiction.

%\input{schwarz_tikz_figure}

%Without loss of generality, consider the decomposition of
%$\Omega$ into two overlapping subdomains, $\Omega_1$ and $\Omega_2$, with
%boundaries $\partial \Omega_1$ and $\partial \Omega_2$, respectively, such that
%$\Omega_1 \cup \Omega_2 = \Omega$, $\Omega_1 \cap \Omega_2 \ne \emptyset$.
%The domain boundary is decomposed as $\partial \Omega = \Gamma_1 \cup
%\Gamma_2$ where $\Gamma_i = \Gamma_i^D \cup \Gamma_i^{T}$, $i=1,2$, is composed of
%a Dirichlet boundary $\Gamma_i^D$ and a Neumann boundary $\Gamma_i^{T}$. 
%\ikt{We should make the notation of the system boundaries consistent with Section \ref{sec:solid_mechanics} (either change notation there or here).}
%The domain boundaries on each
%subdomain are then given as $\partial \Omega_1 = \Gamma_1 \cup \Gamma_{12}$ and
%$\partial \Omega_2 = \Gamma_2 \cup \Gamma_{21}$, where $%\Gamma_{12} = \partial
%\Omega_1 \cap \Omega_2$ and $\Gamma_{21} = \partial \Omega_2 \cap \Omega_1$.
%We refer the reader to Figure~\ref{fig:} \todo{Irina - make/add figure} for a depiction.
Furthermore, let the time domain be decomposed into $\numControllerSteps$ time intervals (referred to as ``controller time-steps" in our earlier work \cite{Mota:2022}), $I_k = [t_{k-1},t_{k}]$, $k=1,\ldots,\numControllerSteps$ with $t_k < t_{k+1}$, $t_{\numControllerSteps} = \timeFinal$.
%{\color{red} This is conflicting with the definition of $t_i$ above Eq. 14. Are we ok with this? We could introduce a time step in 14 and use notation similar to 24 to avoid this.}
SAM operates by solving the following sequence of problems for $k=1,\ldots,\numControllerSteps$, with inner Schwarz iterations $\schwarzIt=1,2,\ldots$,
%{\color{red} I think we need to drop the projection operators since this is defined at the continuous level; the comments on different spatial discretizations doesn't make sense yet. Any arguments?} 
\begin{equation}\label{eq:dynamic_elasticity_sam_bcs} 
    \begin{aligned} \left\{
    \begin{aligned}
        \nabla\cdot ~P_1^{(\schwarzIt)} + \rho ~R_1 &= \rho \ddot{~\varphi}_1^{(\schwarzIt)} &\text{ in } I_k\times \domainOneArg{1}, \\
        {~\varphi}_1^{(\schwarzIt)}(t, ~X) &= \boldsymbol{\chi} &\text{ on } I_k \times \domainBoundaryDirichletOneArg{1} \\ 
        ~P_1^{(\schwarzIt)} ~n &= ~T  &\text{ on } I_k \times \domainBoundaryNeumannOneArg{1}, \\
        {~\varphi}_1^{(\schwarzIt)}(t, ~X) &= {~\varphi}_2^{(n-1)}(t, ~X) &\text{ on } I_k \times  \domainBoundarySchwarzOneArg{1}, 
%\dot{\varphi}_1^{(\schwarzIt)}(t, ~X) &= \dot{{\varphi}}_2^{(n-1)}(t, ~X) \quad \text{ on }
%I_k \times \Gamma_{12}, \\ 
%\ddot{{\varphi}}_1^{(\schwarzIt)}(t, ~X) &=
%\ddot{{\varphi}}_2^{(n-1)}(t, ~X) \quad \text{ on } I_k \times \Gamma_{12},
    \end{aligned} \right. 
    \end{aligned} \quad
    \begin{aligned} \left\{
    \begin{aligned}
        \nabla\cdot ~P_2^{(\schwarzIt)} + \rho ~R_2 &= \rho \ddot{~\varphi}_2^{(\schwarzIt)} &\text{ in } I_k\times\domainOneArg{2}, \\ 
        {~\varphi}_2^{(\schwarzIt)}(t, ~X) &= \boldsymbol{\chi} &\text{ on } I_k \times \domainBoundaryDirichletOneArg{2}, \\
        ~P_2^{(\schwarzIt)} ~n &= ~T  &\text{ on } I_k \times \domainBoundaryNeumannOneArg{2}, \\
        {~\varphi}_2^{(\schwarzIt)}(t, ~X) &= {~\varphi}_1^{(\schwarzIt)}(t, ~X) &\text{ on } I_k \times \domainBoundarySchwarzOneArg{2}, 
%\dot{\varphi}_2^{(\schwarzIt)}(t, ~X) &= \dot{{\varphi}}_1^{(\schwarzIt)}(t, ~X) \quad \text{ on } I_k
%\times \Gamma_{21}, \\ 
%\ddot{{\varphi}}_2^{(\schwarzIt)}(t, ~X) &=
%\ddot{{\varphi}}_1^{(\schwarzIt)}(t, ~X) \quad \text{ on } I_k \times \Gamma_{21}.
        \end{aligned} \right.
        \end{aligned} 
\end{equation} 
given an initial guess for $~\varphi_2^{(0)}$. %In \eqref{eq:dynamic_elasticity_sam_bcs}, $\Pi_1$ denotes the projection of the solution $~\varphi_2$ onto $\partial_S \Omega_1$ and $\Pi_2$ denotes the projection of the solution $~\varphi_1$ onto $\partial_S \Omega_1$.  
These operators are needed since the subdomains $\Omega_1$ and $\Omega_2$ can be discretized by different meshes having different resolutions and/or element types, as shown in Figure \ref{fig:schwarz_overlapping}.  
%\ikt{I added $\Pi_i$, continuous projection operators.  I think it makes \eqref{eq:dynamic_elasticity_sam_bcs} more precise and I use it in Algorithm \ref{algo:sam-fom}.  Thoughts/concerns?}
%\ejp{When I wrote it, I felt it wasn't necessary since the BC definition is only on the boundary. The explanation of different meshes/element types doesn't make sense to me here since this is defined at the PDE level. That being said, I don't think it hurts to have it.}
The system~\eqref{eq:dynamic_elasticity_sam_bcs} is supplemented with initial conditions, as given in~\eqref{eq:dynamic_elasticity_bcs}, restricted to each subdomain. In essence, SAM solves the PDE on $\domainOneArg{1}$ using an initial guess to the solution on $\domainOneArg{2}$ as a Dirichlet boundary condition on the interface $\domainBoundarySchwarzOneArg{1}$, and then proceeds to solve the PDE on $\domainOneArg{2}$ using the updated solution from $\domainOneArg{1}$ as a Dirichlet boundary on the interface $\domainBoundarySchwarzOneArg{2}$. This iteration is continued until convergence, discussed in more detail below. We emphasize that the Schwarz iteration is not performed in a global space-time setting, but rather over each controller time-step. This is equivalent to a global space-time implementation (due to causality in time), but far more efficient and easier to implement, as discussed in~\cite{Mota:2022}. %\ejp{Irina/Alejandro, can you confirm my statement is correct?} \ikt{Yes!}

\begin{remark} \label{remark:many_sds}
    Although the discussion and numerical results presented herein assume the coupling of just two subdomains, we emphasize that SAM is capable of coupling an arbitrary number of subdomains.  Indeed, we have considered SAM-based couplings of three or more subdomains in several of our past works, e.g.,~\cite{Mota:2017, Snyder:2023, Moore:2024, wentland2024Schwarz}. However, since we consider herein DDs that are physically motivated, we are not targeting problems having greater than 10 subdomains; in fact, it is rare to encounter a problem requiring more than 4--5 subdomains.
\end{remark}

\IncMargin{1.5em}
\begin{algorithm}[h]
	\SetKwInOut{Input}{\hspace{-0.5cm}Input}
	\SetKwInOut{Output}{\hspace{-0.5cm}Output}
	%\Input{$\bm{\mathcal{D}} = \{(\mathbf{x}_{n}, \mathbf{y}_{n})\}_{n=1}^{N}$ - the set of observations of two r.v.s;\\ $\lambda$ - parameter of independence}
	%\Output{The causal direction}
	%\BlankLine
%\ikt{This is my initial attempt at writing an algorithm to describe SAM. Still need to check everything.  Feedback welcome.} \ejp{This looks great to me}\;
\Input{Overlapping subdomains $\Omega_1$ and $\Omega_2$ each with a mesh,   initial condition ($~\varphi_1^{(0)}(t_0, ~X)$ and $~\varphi_2^{(0)}(t_0, ~X)$, respectively), time integrator, uniform time-step ($\Delta t_1$ and $\Delta t_2$, respectively). }
\Output{Converged solutions $~\varphi_1(t, ~X)$ and $~\varphi_2(t,~X)$ to the subproblems in \eqref{eq:dynamic_elasticity_sam_bcs} in $\Omega_1$ and $\Omega_2$, respectively.}
	%Create a domain decomposition of $\Omega$ into two overlapping subdomains, $\Omega_1$ and $\Omega_2$ (Figure \ref{fig:schwarz_overlapping}) \;
    %Mesh $\Omega_1$ and $\Omega_2$ \;  
    %Select time-integration schemes and time-steps, $\Delta t_1$ and $\Delta t_2$, respectively, to use in each subdomain \; 
    Divide the global time $[0,T]$ into $C$ time intervals $I_k=[t_{k-1},t_k]$, for $k=1, ..., C$\;
    %Specify initial conditions $~\varphi_1^{(1)}(t_0)$ and $~\varphi_2^{(1)}(t_0)$ for the solutions in $\Omega_1$ and $\Omega_2$, respectively\;
  
    \For{$k=1$ to $C$}{
      Initialize $n=1$ \;
      Initialize \textit{SAM converged} to false \; 
      \While{SAM converged is false }{
        Project $~\varphi_2^{(n-1)}(t, ~X)$ onto $\partial_S \Omega_1$ in space and in time to obtain $~\varphi_2^{(n-1)}(t, ~X)$ on $\partial_S \Omega_1$ for $t = t_{k-1}, t_{k-1} + \Delta t_1, t_{k-1} + 2\Delta t_1, ..., t_{k}$ \;  
        Advance the $\Omega_1$ subproblem in \eqref{eq:dynamic_elasticity_sam_bcs} from time $t_{k-1}$ to $t_k$  to obtain $~\varphi_1^{(n)}(t, ~X)$ in $\Omega_1$ for $t = t_{k-1}, t_{k-1} + \Delta t_1, t_{k-1} + 2\Delta t_1, ..., t_k$\; 
        Project $~\varphi_1^{(n)}(t, ~X)$ onto $\partial_S \Omega_2$ in space and in time to obtain $~\varphi_1^{(n)}(t, ~X)$ on $\partial_S \Omega_2$ for $t = t_{k-1}, t_{k-1} + \Delta t_2, t_{k-1} + 2\Delta t_2, ..., t_{k}$ \;   
        Advance the $\Omega_2$ subproblem in  \eqref{eq:dynamic_elasticity_sam_bcs} from time $t_{k-1}$ to $t_k$  to obtain $~\varphi_2^{(n)}(t, ~X)$ in $\Omega_2$ for $t = t_{k-1}, t_{k-1} + \Delta t_2, t_{k-1} + 2\Delta t_2, ..., t_k$\; 
        Check SAM convergence criteria \eqref{eq:conv_criterion} at time $t_k$ \; 
        \uIf{SAM convergence criteria satisfied}{
          Set \textit{SAM converged} to true
        }
	    \Else{
          Increment $n = n+1$ \; 
        }
      }
    }
	\caption{O-SAM for FOM-FOM coupling in solid dynamics and the specific case of two overlapping subdomains.}\label{algo:sam-fom}
\end{algorithm}
\DecMargin{1.5em}

\begin{figure}[ht!]
    \centering
    \begin{subfigure}{0.48\linewidth}
        \includegraphics[width=1.0\linewidth]{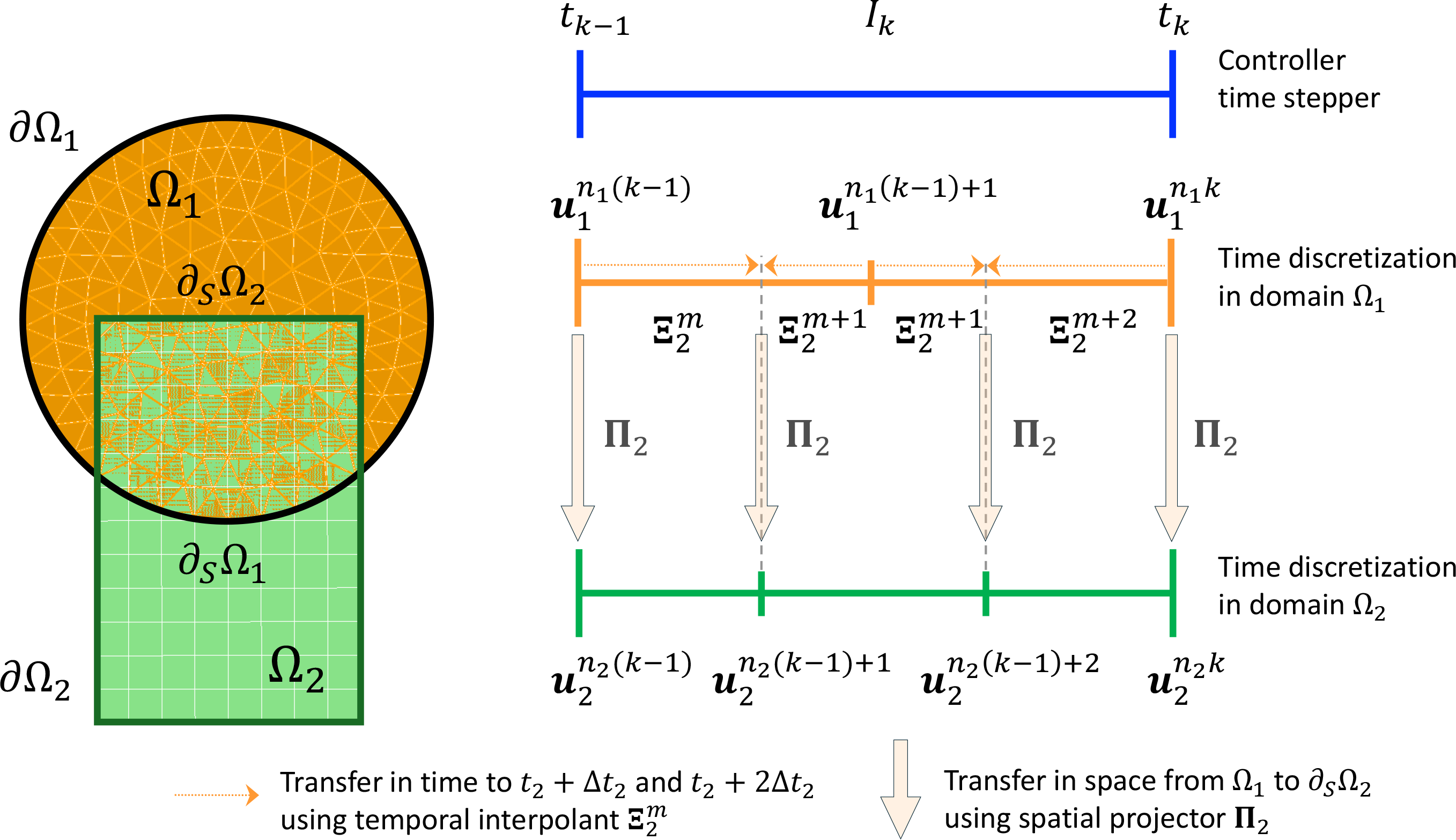}
        \subcaption{$\Omega_1$ as source, $\Omega_2$ as destination}
    \end{subfigure}
    \hspace{0.2cm}
       \begin{subfigure}{0.48\linewidth}
        \includegraphics[width=1.0\linewidth]{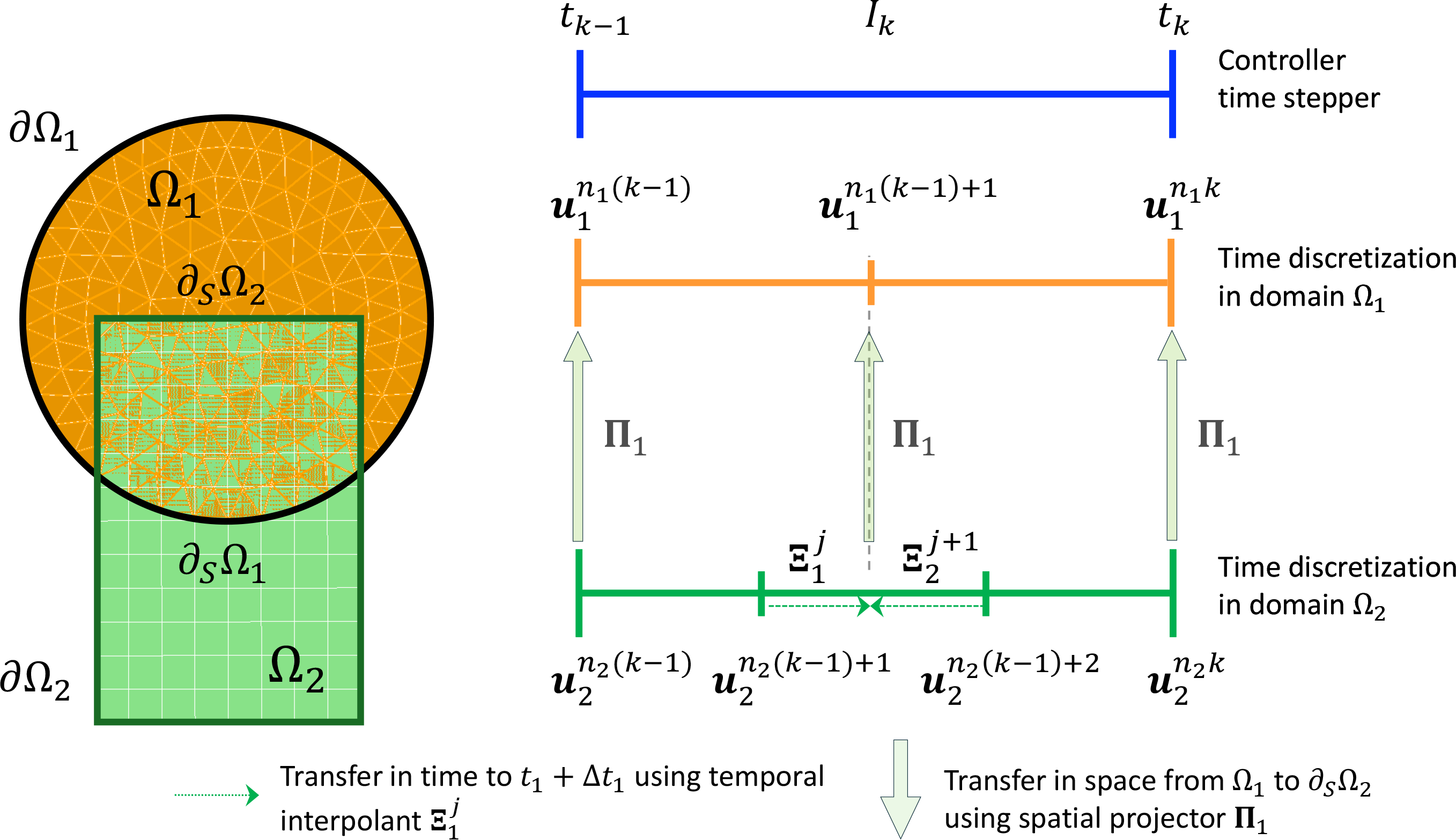}
        \subcaption{$\Omega_2$ as source, $\Omega_1$ as destination}
    \end{subfigure}
    \caption{Field transfer in the O-SAM algorithm within one controller time-step $I_k$ using the spatial projectors $\spatialSchwarzProjectionOperator_1$ and $\spatialSchwarzProjectionOperator_2$, and temporal interpolants $\temporalSchwarzProjectionOperator_1^j$ and $\temporalSchwarzProjectionOperator_2^m$. (a) $\Omega_1$ serves as the source subdomain, while $\Omega_2$ functions as the destination domain; (b) the roles of the domains are reversed.  The super-scripts $j+l$ and $m+l$ correspond to times $t_{k-1} + l\Delta t_1$ and $t_{k-1} + m\Delta t_2$, respectively, and the integers $n_i$ denote the number of time-steps  in $\Omega_i$, as described in~\eqref{eq:schwarz_discrete}.  Thus,  $~u_i^{n_i(k-1)} \approx ~u_i(t_{k-1}, ~X)$ and $~u_i^{n_ik} \approx ~u_i(t_{k}, ~X)$ for $i=1,2$.  %\ikt{I redid this figure.  Eric, please check that I got the indices right.  I don't know if I like going from $k$ to $m$ and $j$.  We can modify if there is a better idea.  I can also remove the projector/interpolants from the figures if you think it's too busy.}
    }
\label{fig:schwarz_time_discretization}
\end{figure}

In practice, we implement SAM at the fully discrete level, and potentially utilize disparate spatio-temporal discretizations on the different subdomains.
%\ikt{I started writing this, but am having second thoughts now...  Before writing our coupling scheme at the fully-discrete level for the specific case of the Newmark-$\beta$ time-integration scheme, we present the algorithm at the semi-discrete level.}
%\todo{We should say specifically we are NOT doing Schwarz in space-time, as this is too intrusive; we want 
%to get Schwarz to work with more traditional discretizations %first in space and in time.  We can say that this 
%is equivalent to space-time, as discussed in~\cite{Mota:2022}.}
Again, without loss of generality, we consider discretization in space using standard isoparametric finite elements and discretization in time using the Newmark-$\beta$ method. For time discretization, we partition each controller time-step $I_k$, $k=1,\ldots,\numControllerSteps$ into $\numTimeStepsOneArg{1}$ and $\numTimeStepsOneArg{2}$ time intervals with time-step $\timeStepOneArg{1}$ and $\timeStepOneArg{2}$ on $\Omega_1$ and $\Omega_2$, respectively; for notational simplicity, we assume a uniform number of time-steps in each time interval, as shown in Figure~\ref{fig:schwarz_time_discretization}. The total number of time-steps\footnote{Note that this assumes that the time-steps $\Delta t_i$ evenly divide the controller time-steps, $t_{k}-t_{k-1}$. This is depicted in Figure \ref{fig:schwarz_time_discretization}.} in $\domainOneArg{i}$ is thus $\numTotalTimeStepsOneArg{i} = \numTimeStepsOneArg{i}\numControllerSteps$. The fully discrete SAM iteration\footnote{To simplify the presentation, we have suppressed the dependence of stiffness matrices $\boldsymbol K_i$ on the Dirichlet and Schwarz boundary conditions.} for the unconstrained DoFs over the $k^{th}$ time interval is: find $\displacementThreeArg{1}{j+1}{(\schwarzIt)}$, $\displacementThreeArg{2}{m+1}{(\schwarzIt)}$ such that, for $j = \numTimeStepsOneArg{1}(k-1) + 1,\ldots,\numTimeStepsOneArg{1}k$, $m=\numTimeStepsOneArg{2} (k-1) ,\ldots,\numTimeStepsOneArg{2} k$, $\schwarzIt = 1,2,\ldots$, 
%\ikt{Doesn't $I_k$ go from $t_{k-1}$ to $t_k$?  I'm a bit confused mapping $j$ to $k$.  If $n_1=1$ (the time-step is equal to the controller time-step), you would have $j$ going from $k$ to $k$, when it should be $k-1$ to $k$, no?  I would expect in general for $j$ to start with $k-1$ and end with $k$, but then the intermediary $j$ would be different.  Am I missing something?  }
%\ikt{I will need to modify Section \ref{sec:solid_mechanics} to have f be a function of the $\chi$, but that's fine.} 
\begin{equation} \label{eq:schwarz_discrete}
\begin{split}
& \frac{\mass_1}{\timeStepOneArg{1}\beta} \displacementThreeArg{1}{j+1}{(\schwarzIt)} + \stiffness_1
\left( \displacementThreeArg{1}{j+1}{(\schwarzIt)}  \right) =
 \boldsymbol f_1 \left( 
\schwarzBoundaryBCThreeArg{1}{j+1}{(\schwarzIt-1)} ,\boldsymbol \chi_{h,1}^{j+1}
\right) +  \mass_1  \left[ \frac{1}{\timeStepOneArg{1}^2 \beta } \displacementTwoArg{1}{j} +  \frac{1}{\timeStepOneArg{1} \beta } \velocityTwoArg{1}{j} + \frac{1}{\beta } \accelerationTwoArg{1}{j} \right] ,  \\
& \frac{\mass_2}{\timeStepOneArg{2} \beta} \displacementThreeArg{2}{m+1}{(\schwarzIt)} + \stiffness_2 \left(
\displacementThreeArg{2}{m+1}{(\schwarzIt)} \right) = 
 \boldsymbol f_2 \left( 
\schwarzBoundaryBCThreeArg{2}{m+1}{(\schwarzIt)} ,\boldsymbol \chi_{h,2}^{m+1}
\right)
+  \mass_2  \left[ \frac{1}{\timeStepOneArg{2}^2 \beta } \displacementTwoArg{2}{m} +  \frac{1}{\timeStepOneArg{2} \beta } \velocityTwoArg{2}{m} + \frac{1}{\beta } \accelerationTwoArg{2}{m} \right]
,  \\
        &\left\{
        \begin{aligned}
 \accelerationThreeArg{i}{\ell+1}{(\schwarzIt)}
 &= \frac{1}{\timeStepOneArg{i} \beta} \left[ \displacementThreeArg{i}{\ell+1}{(\schwarzIt)} - \displacementTwoArg{i}{\ell}
\right]
- \frac{1}{\timeStepOneArg{i} \beta }\velocityTwoArg{i}{\ell} - \frac{1}{2 \beta} \left(1 - 2 \beta\right) \accelerationTwoArg{i}{\ell} \\
\velocityThreeArg{i}{\ell+1}{(\schwarzIt)} &= \velocityTwoArg{i}{\ell}  + (1 - \gamma) \timeStepOneArg{i} \accelerationTwoArg{i}{\ell} + \gamma \timeStepOneArg{i} \accelerationThreeArg{i}{\ell+1}{(\schwarzIt)} \\
\end{aligned} \qquad (i,\ell)=(1,j) ,(2,m) 
\right.
\end{split}
\end{equation}
until convergence is reached based on the absolute or relative error tolerance $\delta_{\text{abs}} , \delta_{\text{rel}} > 0$, 
\begin{equation}
\label{eq:conv_criterion}
\begin{split}
\displacementTwoArg{i}{\numTimeStepsOneArg{i} k} =   \displacementThreeArg{i}{\numTimeStepsOneArg{i} k}{(\schwarzIt)}   \\
\velocityTwoArg{i}{\numTimeStepsOneArg{i} k} =   \velocityThreeArg{i}{\numTimeStepsOneArg{i} k}{(\schwarzIt)} \\
\accelerationTwoArg{i}{\numTimeStepsOneArg{i} k} =   \accelerationThreeArg{i}{\numTimeStepsOneArg{i} k}{(\schwarzIt)} \\
\end{split} \;
\; \mathrm{s.t.} 
\begin{cases}
 & \sum\limits_{i=1}^2 \big\lVert \left( \displacementThreeArg{i}{\numTimeStepsOneArg{i} k}{(\schwarzIt)} -\displacementThreeArg{i}{\numTimeStepsOneArg{i} k}{(\schwarzIt -1)} \right) 
  + \timeStepOneArg{i}\, \left(\velocityThreeArg{i}{\numTimeStepsOneArg{i} k}{(\schwarzIt)}  - \velocityThreeArg{i}{\numTimeStepsOneArg{i} k}{(\schwarzIt -1 )} \right)
  \big\rVert_2^2
  < \delta_{\text{abs}}^2, \\
  & \mathrm{or} \\
  & \frac{  \sum\limits_{i=1}^2 \big\lVert \left( \displacementThreeArg{i}{\numTimeStepsOneArg{i} k}{(\schwarzIt)} -\displacementThreeArg{i}{\numTimeStepsOneArg{i} k}{(\schwarzIt -1)} \right) 
  + \timeStepOneArg{i}\, \left(\velocityThreeArg{i}{\numTimeStepsOneArg{i} k}{(\schwarzIt)}  - \velocityThreeArg{i}{\numTimeStepsOneArg{i} k}{(\schwarzIt -1 )} \right)
  \big\rVert_2^2}{\sum\limits_{i=1}^2
  \big\lVert  
  \displacementThreeArg{i}{\numTimeStepsOneArg{i} k}{(\schwarzIt)} 
  + \Delta t\, \velocityThreeArg{i}{\numTimeStepsOneArg{i} k}{(\schwarzIt)} 
 \big\rVert_2^2} 
  < \delta_{\text{rel}}^2. 
\end{cases}
\end{equation} 
%\ikt{Should $\ddot{~u}$ be in the previous equation?  The acceleration does not appear in the convergence criterion.}
%for selected tolerances $\delta_{\text{rel}}, \delta_{\text{abs}} > 0$.
In the above, $\gamma, \beta$ are parameters of the Newmark-$\beta$ method, $\displacementThreeArg{i}{j}{(\schwarzIt)}\in \mathbb{R}^{\numFreeDofsOneArg{i}}$ denotes the solution at the $n^{th}$ Schwarz iteration for the $\numFreeDofsOneArg{i}$ unconstrained DoFs for displacements on $\Omega_i$, $i=1,2$, at the $j^{th}$ time-step on the interval. Additionally, $\mass_i \in \mathbb{R}^{\numFreeDofsOneArg{i} \times \numFreeDofsOneArg{i}}$ and $\stiffness_i: \mathbb{R}^{\numFreeDofsOneArg{i}} \rightarrow \mathbb{R}^{\numFreeDofsOneArg{i}}$, $i=1,2$ are the mass matrix and stiffness, while $\boldsymbol f_i : \mathbb{R}^{\numSchwarzDirichletBoundaryDofsOneArg{i}} \times  \mathbb{R}^{\numDirichletBoundaryDofsOneArg{i} }  \rightarrow \mathbb{R}^{\numFreeDofsOneArg{i}}$, $i=1,2$ are the forcing vectors which account for the boundary conditions $\schwarzBoundaryBCThreeArg{i}{j}{(\schwarzIt)}$ on the $\numSchwarzDirichletBoundaryDofsOneArg{i}$ interior interface DoFs and $\boldsymbol \chi_{h,i}$ on the $\numDirichletBoundaryDofsOneArg{i}$ exterior Dirichlet interface DoFs. Critical to the Schwarz formulation, the interior Schwarz boundary terms $\schwarzBoundaryBCThreeArg{i}{j}{(\schwarzIt)}$ are given as
\begin{align*}
\schwarzBoundaryBCThreeArg{1}{j}{(\schwarzIt)} &= \left[ \ 
 \spatialSchwarzProjectionOperator_1  \displacementThreeArg{2}{1}{(\schwarzIt)},\ldots, 
  \spatialSchwarzProjectionOperator_1  \displacementThreeArg{2}{\numTimeStepsOneArg{2}}{(\schwarzIt)} \
\right] \temporalSchwarzProjectionOperator^j_1 \in \RR{\numSchwarzDirichletBoundaryDofsOneArg{1}}, \\
\schwarzBoundaryBCThreeArg{2}{m}{(\schwarzIt)} &=  \left[ \  
\spatialSchwarzProjectionOperator_2 \displacementThreeArg{1}{1}{(\schwarzIt)} ,\ldots, 
\spatialSchwarzProjectionOperator_2  \displacementThreeArg{1}{\numTimeStepsOneArg{1}}{(\schwarzIt)}  \
\right] \temporalSchwarzProjectionOperator^m_2 \in \RR{\numSchwarzDirichletBoundaryDofsOneArg{2}},
\end{align*}
%\ikt{I would consider adding a second index to $\spatialSchwarzProjectionOperator$ to show from where you are projecting, but I realize it will add more indices.  }
where $\spatialSchwarzProjectionOperator_1 \in \mathbb{R}^{ \numSchwarzDirichletBoundaryDofsOneArg{1} \times \numFreeDofsOneArg{2}}$ performs pointwise projection\footnote{%For non-uniform spatial discretizations, the projection operator interpolates pointwise in space. 
We note that it is possible to use variational projection, as discussed in \cite{Mota:2025}. %\ejp{Irina/Alejandro, is this OK?} \ikt{yes!}
} of the spatial DoFs in $\domainOneArg{2}$ onto the boundary $\domainBoundarySchwarzOneArg{1}$, and $\temporalSchwarzProjectionOperator_1^j \in \mathbb{R}^{\numTimeStepsOneArg{2} \times 1}$ is a temporal interpolant for the case where non-uniform time schemes are used on the different subdomains.  Similarly,  $\spatialSchwarzProjectionOperator_2 \in \mathbb{R}^{\numSchwarzDirichletBoundaryDofsOneArg{2} \times \numFreeDofsOneArg{1}}$ and $\temporalSchwarzProjectionOperator_2^m \in \mathbb{R}^{\numTimeStepsOneArg{1} \times 1}$. Note that interpolation in time is only required for non-equivalent temporal discretizations between the two domains; we again refer the reader to Figure~\ref{fig:schwarz_time_discretization}.  
%\todo{We should say somewhere we do pointwise projection (vs. variational projection).  Maybe this belongs in the numerical results section.}

%\input{schwarz_time_fig}

We emphasize that, while the discussion above details the implementation of O-SAM for the specific case of a Newmark-$\beta$ time integrator, the subdomain problems in $\Omega_1$ and $\Omega_2$ can be advanced in time using \textit{any} time-integration scheme.  This is made clear in Algorithm \ref{algo:sam-fom}, which outlines the general O-SAM scheme as applied to the subdomain-local PDEs in \eqref{eq:dynamic_elasticity_sam_bcs}.  As in the discussion above and in Figure \ref{fig:schwarz_time_discretization}, we assume in Algorithm \ref{algo:sam-fom} that the subdomain time-steps $\Delta t_1$ and $\Delta t_2$ evenly divide the controller time-steps.

\begin{remark} \label{remark:additive_sam}
    The present work restricts attention to the serial version of SAM, commonly known in the linear solver literature as ``multiplicative Schwarz" \cite{Gander:2008}.  In particular, the subdomain problems in \eqref{eq:dynamic_elasticity_sam_bcs} are solved sequentially. It is possible to improve the efficiency of SAM by implementing its ``additive" variant \cite{Gander:2008}, in which the subdomain problems in \eqref{eq:dynamic_elasticity_sam_bcs} are solved in parallel, with boundary data communicated between subdomains asynchronously as it becomes available.  Assessing the performance of a parallelized additive SAM on the test cases considered herein would be an interesting future research endeavor.  Our earlier work \cite{wentland2024Schwarz} has shown that, for some 2D nonlinear hyperbolic conservation law problems with a $2\times 2$ domain decomposition, speedups of $\sim 3\times$ are possible by employing parallel additive SAM over its sequential multiplicative analog. In general, the speedup expected will depend on the problem and the domain decomposition.
\end{remark}

\subsection{O-SAM for OpInf-FOM and OpInf-OpInf coupling} 
As outlined above, SAM is flexible with respect to different  discretizations employed within the different subdomains, and can be formulated for the case where an approximate ROM is used to model a subdomain. We now describe this process for OpInf in terms of the relevant offline and online stages. %For simplicity of presentation, we initially 
In the present work, we restrict attention to the case where we construct an OpInf ROM for $\domainOneArg{1}$ using snapshots of a corresponding FOM-FOM Schwarz simulation, often referred to in the literature as ``top-down" training~\cite{Chung:2024}. The decision to focus on ``top-down" training is based on the assumption that a FOM-FOM coupled solution on the full physical domain of interest is generally possible to obtain via SAM using our software frameworks.  
%
%We emphasize, however, that, as FOM-FOM Schwarz converges to the monolithic FOM-only solution (provided equivalent discretizations in space and time), an equivalent OpInf ROM could be constructed from snapshots of a monolithic FOM provide proper restriction operators. 
%We further note that the resulting 
We note that, since our formulation builds the OpInf ROM on $\domainOneArg{1}$ purely from information on $\domainOneArg{1}$ and $\domainBoundaryOneArg{1}$, it is, by construction, a \textit{subdomain-local} approach. 
%This type of approach 
As such, it is compatible with both ``top-down" training and ``bottom-up" training.
%; we only consider the former in the present work for simplicity. 
%\ejp{Ian, I think this is what you're referring to?}

As in Section \ref{sec:schwarz_fom}, we present our O-SAM-based OpInf-FOM and OpInf-OpInf formulation assuming all models have been discretized in time using the Newmark-$\beta$ time-integration scheme. We again emphasize that SAM is extremely flexible and allows for the use of \textit{any} time integration scheme in each subdomain.   For clarity, the online and offline stages of O-SAM for OpInf-OpInf coupling with a Newmark-$\beta$ time-stepping scheme is summarized in Algorithm \ref{algo:rom-offline-online}; extensions to OpInf-FOM coupling following the approach described in Section \ref{sec:opinf-fom} is straightforward, as is the extension to alternate time integrators.

\begin{remark} \label{remark:training}
    We intentionally do not perform training on the monolithic full domain $\Omega$ and generate subdomain-local POD bases by restricting or projecting the full domain snapshots onto the subdomains $\Omega_i$ (as done in e.g., \cite{Moore:2024}), since many of our problems of interest are posed on complex geometries that are extremely difficult to mesh without performing a domain decomposition into simpler subdomains.
\end{remark}

%\todo{Need to restructure the following. }

\label{sec:schwarz_rom}

\subsubsection{Offline stage: basis construction and operator inference}
\label{sec:schwarz_opinf_offline}
%Inspired by the form of the FOM, we postulate the OpInf model for the unconstrained DoFs on $\Omega_1$
As discussed in Section \ref{sec:opinf},  we postulate an OpInf model for the unconstrained DoFs on $\Omega_1$ of the form 
% \begin{equation} \label{eq:opinf_schwarz_discrete}
% \begin{split}
%  &\reducedAccelerationOneArg{1} + \reducedStiffness_1
% \reducedDisplacementOneArg{1} 
% + \reducedStiffnessQuadratic_1
% \left( \reducedDisplacementOneArg{1} \right)^{\otimes 2}
% + \reducedStiffnessCubic_1
% \left( \reducedDisplacementOneArg{1} \right)^{\otimes 3}
% =
%  \reducedBoundaryOperator_1 \mathbf{g}_1,%\left[ \reducedBoundaryBCOneArg{1} , \hat{\boldsymbol \chi}_{h} \right]  
%  \end{split}
% \end{equation}
\begin{equation} \label{eq:opinf_model_schwarz}
\begin{split}
 &\reducedAccelerationOneArg{1} + \reducedStiffness_1
\reducedDisplacementOneArg{1} 
+ \reducedStiffnessQuadratic_1
\reducedDisplacementOneArg{1}^{\otimes 2}
+ \reducedStiffnessCubic_1
\reducedDisplacementOneArg{1}^{\otimes 3}
=
 \reducedBoundaryOperator_1   \reducedBoundaryBCOneArg{1},%\left[ \reducedBoundaryBCOneArg{1} , \hat{\boldsymbol \chi}_{h} \right]  
 \end{split}
\end{equation}
where $\reducedDisplacementOneArg{1}(t) \in \mathbb{R}^{\romDimOneArg{1}}$ are the reduced coordinates for the solution on $\Omega_1$ with the associated basis $\trialBasis_1 \in \RRStar{\numFreeDofsOneArg{1}}{\romDimOneArg{1}}$, and $\reducedStiffness_1 \in \mathbb{R}^{\romDimOneArg{1} \times \romDimOneArg{1}}$,  $\reducedStiffnessQuadratic_1 \in \mathbb{R}^{\romDimOneArg{1} \times \romDimOneArg{1}^2}$, $\reducedStiffnessCubic_1 \in \mathbb{R}^{\romDimOneArg{1} \times \romDimOneArg{1}^3}$ are the linear, quadratic, and cubic reduced operators, respectively. Lastly, key to the Schwarz formulation, $\reducedBoundaryOperator_1$ and $  \reducedBoundaryBCOneArg{1}$ are used to model the Dirichlet boundary contributions; this will now be discussed.

Unlike the monolithic OpInf formulation presented in Section~\ref{sec:opinf}, in the Schwarz formulation the boundary term  must model both the essential Dirichlet boundary condition, $ {~\varphi}_1(t, ~X) = \boldsymbol{\chi} \; \text{ on } I \times \domainBoundaryDirichletOneArg{1}$, as well as the Schwarz overlap boundary condition, ${~\varphi}_1(t, ~X) = {~\varphi}_2(t, ~X) \; \text{ on } I \times \domainBoundarySchwarzOneArg{1}$. The approach employed in~\cite{Moore:2024} is to set
\begin{equation*}
    \reducedBoundaryBCOneArg{1} := \left( \
        \left[ \spatialSchwarzProjectionOperator_1 \displacement_2 \right]^{\intercal} \ , \
        \boldsymbol{\chi}_{h,1}^{\intercal}
    \ \right)^{\intercal} 
    \in \mathbb{R}^{\numSchwarzDirichletBoundaryDofsOneArg{1} + \numDirichletBoundaryDofsOneArg{1}}. 
\end{equation*}
This approach, however, results in an OpInf operator that scales with the FOM dimension as $\reducedBoundaryOperator_1 \in \RR{\romDimOneArg{1} \times \left( \numSchwarzDirichletBoundaryDofsOneArg{1} + \numDirichletBoundaryDofsOneArg{1}\right)}$. 

For 3D problems, the number of interface DoFs for $ \numSchwarzDirichletBoundaryDofsOneArg{1} + \numDirichletBoundaryDofsOneArg{1}$ can be large. 
To mitigate this issue, we employ a two-stage approach for the boundary operator wherein we employ POD to reduce the dimension of the boundary condition forcing vector and solve an OpInf problem in a reduced space. We set
\begin{equation*}
    \reducedBoundaryBCOneArg{1} := \left( \ 
    \left[ \trialBasisSchwarz^{\intercal} \spatialSchwarzProjectionOperator_1 \displacement_2 \right]^{\intercal} \ , \
     \left[ \trialBasisDirichlet^{\intercal} \boldsymbol{\chi}_{h,1} \right]^{\intercal} 
     \ \right)^{\intercal}
     \in \mathbb{R}^{\numReducedSchwarzDirichletBoundaryDofsOneArg{1} + \numReducedDirichletBoundaryDofsOneArg{1}},
\end{equation*}
where 
$\trialBasisDirichlet \in \RRStar{\numDirichletBoundaryDofsOneArg{1}}{ \numReducedDirichletBoundaryDofsOneArg{1}}$ with $\numReducedDirichletBoundaryDofsOneArg{1} \ll \numDirichletBoundaryDofsOneArg{1}$ 
and 
$\trialBasisSchwarz \in \RRStar{\numSchwarzDirichletBoundaryDofsOneArg{1}}{\numReducedSchwarzDirichletBoundaryDofsOneArg{1}}$ with $\numReducedSchwarzDirichletBoundaryDofsOneArg{1} \ll \numSchwarzDirichletBoundaryDofsOneArg{1}$ are orthonormal bases for the essential Dirichlet and Schwarz Dirichlet boundary conditions, respectively.
The above scales independently of the FOM.

%\subsubsection{Stage 1: Data collection and proper orthogonal decomposition}
%\label{sec:schwarz_pod}
To construct the OpInf operators and bases, we suppose that the Schwarz FOM problem~\eqref{eq:schwarz_discrete} has been solved in time for $\numTotalTimeStepsOneArg{1}+1 $ separate states\footnote{Note that the number of snapshots in $\Omega_1$ and $\Omega_2$, denoted by $\tau_1$ and $\tau_2$, respectively, can be different, as our algorithm allows one to use different time-steps in different subdomains when applying O-SAM.}, yielding the collection of unconstrained state snapshots on $\domainOneArg{1}$
\begin{equation} 
\displacementSnapshots_1 = [\displacementTwoArg{1}{0},\displacementTwoArg{1}{{1}}, \ldots, \displacementTwoArg{1}{\numTotalTimeStepsOneArg{1}}] \in \mathbb{R}^{\numFreeDofsOneArg{1} 
\times  (\numTotalTimeStepsOneArg{1}+1)},
 \label{eq:schwarz_snapshots}
\end{equation}
where $\displacementTwoArg{1}{j} \approx \displacementOneArg{1}(j \timeStepOneArg{1})$.
%\ikt{It's a little bit weird to imply $\Omega_1$ and $\Omega_2$ will have a different number of snapshots $\tau_1$ and $\tau_2$, if the snapshots are generated by doing a coupled FOM-FOM simulation, which will have a fixed number of times/parameter variations.  I suppose for some reason if you wanted to use a different number of these snapshots $\tau_i$ to build the POD bases in different subdomains, you could do that, so maybe it's OK to keep it like this.}
%\ejp{I was doing this for the case where we have different time-steps on the different domains - does this make sense or am I missing something?}
We additionally assume access to the set of snapshots for the essential Dirichlet boundary conditions and Schwarz Dirichlet boundary conditions,
\begin{equation*}
\begin{split}
&\dirichletBoundarySnapshots_1 = \left[ \boldsymbol \chi_{h,1}^0 ,\ldots,  \boldsymbol \chi_{h,1}^{\numTotalTimeStepsOneArg{1}} \right] \in \RR{\numDirichletBoundaryDofsOneArg{1} \times (\numTotalTimeStepsOneArg{1}+1)},\\
&{\schwarzBoundarySnapshots}_1 = \left[ \schwarzBoundaryBCTwoArg{1}{0}, \ldots , \schwarzBoundaryBCTwoArg{1}{\numTotalTimeStepsOneArg{1}}   \right] \in \RR{\numSchwarzDirichletBoundaryDofsOneArg{1} \times (\numTotalTimeStepsOneArg{1}+1)}.
\end{split}
\end{equation*}
Similarly to Section~\ref{sec:opinf}, we construct bases $\trialBasis_1$, $\trialBasisDirichlet$, 
$\trialBasisSchwarz $ for the
unconstrained DoFs on $\domainOneArg{1}$, essential Dirichlet boundary conditions, and Schwarz Dirichlet boundary conditions by employing POD on the snapshot matrices $\displacementSnapshots_1$, $\dirichletBoundarySnapshots_1$ and ${\schwarzBoundarySnapshots}_1$, respectively.We again emphasize that the training data could be constructed using alternate approaches, e.g., 
%be equivalently obtained from a monolithic FOM provided the appropriate restriction operators, or alternatively we could construct snapshot data in 
a ``bottom-up" approach by simulation the FOM on $\domainOneArg{1}$ only for, e.g., varying boundary conditions.
%\irm{Everything we say here is correct, but I think it is worth directly saying that this is a subdomain local approach, as that is a point in our favor as to the method's flexibility. We could comment that this is flexible with respect to the data source.}

\begin{remark}
    The reduced bases considered in this work are constructed from snapshots generated by the coupled FOM-FOM Schwarz simulation, i.e., they are constructed via a ``top-down" training approach. An alternative approach would be to generate training data independently on each subdomain and construct the reduced bases using only local simulations (the ``bottom-up" training approach \cite{Chung:2024}, thereby eliminating the need for coupled high-fidelity training runs. While such a strategy could reduce the offline training cost, the resulting bases may not fully capture the coupled interface dynamics present in the global system and may therefore need to be combined with online adaptation, model switching and/or basis enrichment techniques.  
        %require additional enrichment or adaptation strategies. 
        Since a coupled FOM-FOM Schwarz solver is available in our setting, we adopt the coupled-training approach, and leave the investigation of locally trained Schwarz ROMs to future work.
\end{remark}

Using the collected snapshot matrices $\displacementSnapshots_1$, $\dirichletBoundarySnapshots_1$, ${\schwarzBoundarySnapshots}_1$ and the reduced bases $\trialBasis_1$, $\trialBasisDirichlet$, and $\trialBasisSchwarz$, 
the OpInf operators are inferred from solving the least-squares problem
\begin{multline}\label{eq:opinf_schwarz_inference}
\argmin_{\reducedStiffness_1, \reducedStiffnessQuadratic_1,\reducedStiffnessCubic_1, \reducedBoundaryOperator_1} \left\|D^2_t(\bar{~U}_1) + \reducedStiffness_1 \bar{~U}_1 +\reducedStiffnessQuadratic_1 \bar{~U}_1^{\ast 2} + \reducedStiffnessCubic_1 \bar{~U}_1^{\ast 3} - \reducedBoundaryOperator_1 \left[ \left[ \trialBasisDirichlet^{\intercal} \dirichletBoundarySnapshots_1 \right]^{\intercal}  \; \left[ \trialBasisSchwarz^{\intercal} \schwarzBoundarySnapshots_1\right]^{\intercal} \right]^{\intercal} \right\|_F^2 + \\\regularizationParameter \big( \|\reducedStiffness_1\|_F^2 + \|\reducedStiffnessQuadratic_1\|_F^2 + \|\reducedStiffnessCubic_1\|_F^2 + \|\reducedBoundaryOperator_1\|_F^2 \big),
\end{multline}
where $\bar{~U}_1 = \trialBasis_1^{\intercal}~U_1$. 
Similar to standard OpInf presented in Section~\ref{sec:opinf}, the formulation~\eqref{eq:opinf_schwarz_inference} is non-intrusive and the OpInf ROM can be inferred directly from solving a convex minimization problem with a least-squares objective.

\subsubsection{Online stage: OpInf-FOM coupling} \label{sec:opinf-fom}

%\ikt{I just realized something we did not highlight at all in the paper.  In the case of OpInf-FOM coupling, the FOM is NOT the quadratic form \eqref{eq:quadratic_model_with_bcs_approx} -- it is the full nonlinear model.  The equations here make that clear, but we do not make it clear earlier.  We should.  I will add.}

Discretizing the OpInf ROM with Newmark-$\beta$, 
the fully discrete O-SAM iteration for an OpInf-FOM coupled system over the $k^{th}$ controller time interval are, for $j =
\numTimeStepsOneArg{1}(k-1) + 1,\ldots,\numTimeStepsOneArg{1}k$, $m=\numTimeStepsOneArg{2} (k-1) ,\ldots,\numTimeStepsOneArg{2} k$, $\schwarzIt = 1,2,\ldots$,
\begin{equation} \label{eq:schwarz_opinf_discrete}
\begin{split}
& \frac{\reducedDisplacementThreeArg{1}{j+1}{(\schwarzIt)}}{\timeStepOneArg{1} \beta}
+ \reducedStiffness_1
\reducedDisplacementThreeArg{1}{j+1}{(\schwarzIt)}
+ \reducedStiffnessQuadratic_1
{\reducedDisplacementThreeArg{1}{j+1}{(\schwarzIt)}}{^{\otimes 2}}
+ \reducedStiffnessCubic_1
{\reducedDisplacementThreeArg{1}{j+1}{(\schwarzIt)}}{^{\otimes 3}}
=
 \reducedBoundaryOperator_1    \reducedBoundaryBCThreeArg{1}{j+1}{(\schwarzIt-1)} + 
\frac{1}{\timeStepOneArg{1}^2 \beta } \reducedDisplacementTwoArg{1}{j} +  \frac{1}{\timeStepOneArg{1} \beta } \reducedVelocityTwoArg{1}{j} + \frac{1}{\beta } \reducedAccelerationTwoArg{1}{j}
 ,  \\
& \frac{\mass_2}{\timeStepOneArg{2} \beta} \displacementThreeArg{2}{m+1}{(\schwarzIt)} + \stiffness_2 \left(
\displacementThreeArg{2}{m+1}{(\schwarzIt)} \right) = 
 \boldsymbol f_2 \left( 
\approximateSchwarzBoundaryBCThreeArg{2}{m+1}{(\schwarzIt)} ,\boldsymbol \chi_{h}^{m+1} \right) 
+  \mass_2  \left[ \frac{1}{\timeStepOneArg{2}^2 \beta } \displacementTwoArg{2}{m} +  \frac{1}{\timeStepOneArg{2} \beta } \velocityTwoArg{2}{m} + \frac{1}{\beta } \accelerationTwoArg{2}{m} \right]
,  \\
& \reducedAccelerationThreeArg{1}{j+1}{(\schwarzIt)}
 = \frac{1}{\timeStepOneArg{1} \beta} \left[ \reducedDisplacementThreeArg{1}{j+1}{(\schwarzIt)} - \reducedDisplacementTwoArg{1}{j}
\right]
- \frac{1}{\timeStepOneArg{1} \beta }\reducedVelocityTwoArg{1}{j} - \frac{1}{2 \beta} \left(1 - 2 \beta\right) \reducedAccelerationTwoArg{1}{j}, \\
 &\reducedVelocityThreeArg{1}{j+1}{(\schwarzIt)}= \reducedVelocityTwoArg{1}{j}  + (1 - \gamma) \timeStepOneArg{1} \reducedAccelerationTwoArg{1}{j} + \gamma \timeStepOneArg{1} \reducedAccelerationThreeArg{1}{j+1}{(\schwarzIt)}, \\
& \accelerationThreeArg{2}{m+1}{(\schwarzIt)}
 = \frac{1}{\timeStepOneArg{2} \beta} \left[ \displacementThreeArg{2}{m+1}{(\schwarzIt)} - \displacementTwoArg{2}{m}
\right]
- \frac{1}{\timeStepOneArg{2} \beta }\velocityTwoArg{2}{m} - \frac{1}{2 \beta} \left(1 - 2 \beta\right) \accelerationTwoArg{2}{m}, \\
 &\velocityThreeArg{2}{m+1}{(\schwarzIt)}= \velocityTwoArg{2}{m}  + (1 - \gamma) \timeStepOneArg{2} \accelerationTwoArg{2}{m} + \gamma \timeStepOneArg{2} \accelerationThreeArg{2}{m+1}{(\schwarzIt)}, \\
\end{split}
\end{equation}
where the boundary couplings are defined as
\begin{equation*}
\begin{split}
\reducedBoundaryBCThreeArg{1}{j}{(\schwarzIt)} &= 
\left[  \left[ \trialBasisSchwarz^{\intercal} \schwarzBoundaryBCThreeArg{1}{j}{(\schwarzIt)} \right]^{\intercal}  ,\left[ \trialBasisDirichlet^{\intercal} \boldsymbol \chi_{h}^{j} \right]^{\intercal}\right]^{\intercal} \\
\approximateSchwarzBoundaryBCThreeArg{2}{m}{(\schwarzIt)} &=  \left[   
\spatialSchwarzProjectionOperator_2 \trialBasis_1 \reducedDisplacementThreeArg{1}{1}{(\schwarzIt)} ,\ldots, 
\spatialSchwarzProjectionOperator_2  \trialBasis_1 \reducedDisplacementThreeArg{1}{\numTimeStepsOneArg{1}}{(\schwarzIt)} 
\right] \temporalSchwarzProjectionOperator^m_2. \qquad
\end{split}
\end{equation*}
The iteration is continued until convergence, as described in ~\eqref{eq:conv_criterion}.   We emphasize the following aspects of the formulation~\eqref{eq:schwarz_opinf_discrete}:
\begin{enumerate}
\item As stated earlier in Remark \ref{remark:nonlinear_fom}, while the subdomain-local OpInf ROM in an OpInf-FOM coupling is based on a polynomial approximation of the form \eqref{eq:quadratic_ROM}, the subdomain-local FOM to which it is coupled is based on the full nonlinear PDEs \eqref{eq:dynamic_elasticity_sam_bcs}.  This is evident in \eqref{eq:schwarz_opinf_discrete}.  %\ikt{I added this to address my comment above.} 
\item The OpInf ROM and FOM are coupled through forcing terms that describe the Schwarz Dirichlet boundary condition.
\item When communicating the ROM information on $\Omega_1$ to the FOM on $\Omega_2$ through the term $\approximateSchwarzBoundaryBCTwoArg{2}{m}$, the boundary term needs to be computed in the full order space. This is achieved via the trial basis, $\trialBasis_1$, which in practice can be restricted to the relevant boundary DoFs to accelerate the computation.
\item When communicating the FOM information on $\Omega_2$ to the ROM on $\Omega_1$ through the term $\schwarzBoundaryBCTwoArg{1}{j}$, the FOM information is projected onto the reduced boundary space described by $\trialBasisSchwarz$.
%\item For clarity, the OpInf-OpInf coupling algorithm described above is summarized at a relatively high-level in Algorithm \ref{algo:rom-offline-online} in terms of the offline and online stages of the model reduction and coupling procedure. 
\end{enumerate}

\subsubsection{Online stage: OpInf-OpInf coupling} \label{sec:opinf-opinf}
We highlight that the above formulation can be easily translated to perform Schwarz-based couplings for various subdomain-local OpInf models. Specifically, one can repeat Section~\ref{sec:schwarz_opinf_offline} for $\Omega_2$, and a coupled OpInf-OpInf system can be formulated as
\begin{equation*} \label{eq:schwarz_opinf_opinf_discrete}
\begin{split}
& \frac{\reducedDisplacementThreeArg{1}{j+1}{(\schwarzIt)}}{\timeStepOneArg{1} \beta}
+ \reducedStiffness_1
\reducedDisplacementThreeArg{1}{j+1}{(\schwarzIt)}
+ \reducedStiffnessQuadratic_1
{\reducedDisplacementThreeArg{1}{j+1}{(\schwarzIt)}}{^{\otimes 2}}
+ \reducedStiffnessCubic_1
{\reducedDisplacementThreeArg{1}{j+1}{(\schwarzIt)}}{^{\otimes 3}}
=
 \reducedBoundaryOperator_1    \reducedBoundaryBCThreeArg{1}{j+1}{(\schwarzIt-1)} + 
\frac{1}{\timeStepOneArg{1}^2 \beta } \reducedDisplacementTwoArg{1}{j} +  \frac{1}{\timeStepOneArg{1} \beta } \reducedVelocityTwoArg{1}{j} + \frac{1}{\beta } \reducedAccelerationTwoArg{1}{j}
 ,  \\
& \frac{\reducedDisplacementThreeArg{2}{m+1}{(\schwarzIt)}}{\timeStepOneArg{2} \beta}
+ \reducedStiffness_2
\reducedDisplacementThreeArg{2}{m+1}{(\schwarzIt)}
+ \reducedStiffnessQuadratic_2
{\reducedDisplacementThreeArg{2}{m+1}{(\schwarzIt)}}{^{\otimes 2}}
+ \reducedStiffnessCubic_2
{\reducedDisplacementThreeArg{2}{m+1}{(\schwarzIt)}}{^{\otimes 3}}
=
 \reducedBoundaryOperator_2    \reducedBoundaryBCThreeArg{2}{m+1}{(\schwarzIt)} + 
\frac{1}{\timeStepOneArg{2}^2 \beta } \reducedDisplacementTwoArg{2}{m} +  \frac{1}{\timeStepOneArg{2} \beta } \reducedVelocityTwoArg{2}{m} + \frac{1}{\beta } \reducedAccelerationTwoArg{2}{m}
 ,  \\
\end{split}
\end{equation*}
along with the associated Newmark-$\beta$ updates. The boundary couplings are now defined as
\begin{equation*}
\begin{split}
\reducedBoundaryBCThreeArg{1}{j}{(\schwarzIt)} &= 
\left[  \left[ \trialBasisSchwarz^{\intercal} \approximateSchwarzBoundaryBCThreeArg{1}{j}{(\schwarzIt)} \right]^{\intercal} ,\left[ \trialBasisDirichlet^{\intercal} \boldsymbol \chi_{h}^{j} \right]^{\intercal}\right]^{\intercal}, \qquad 
\reducedBoundaryBCThreeArg{2}{m}{(\schwarzIt)} = 
\left[  \left[ \trialBasisSchwarzTwo^{\intercal} \approximateSchwarzBoundaryBCThreeArg{2}{m}{(\schwarzIt)} \right]^{\intercal}  ,\left[ \trialBasisDirichletTwo^{\intercal} \boldsymbol \chi_{h}^{m} \right]^{\intercal}\right]^{\intercal},
\\
\approximateSchwarzBoundaryBCThreeArg{1}{j}{(\schwarzIt)} &=  \left[   
\spatialSchwarzProjectionOperator_1 \trialBasis_2 \reducedDisplacementThreeArg{2}{1}{(\schwarzIt)} ,\ldots, 
\spatialSchwarzProjectionOperator_1  \trialBasis_2 \reducedDisplacementThreeArg{2}{\numTimeStepsOneArg{2}}{(\schwarzIt)} 
\right] \temporalSchwarzProjectionOperator_1^j, \qquad
\approximateSchwarzBoundaryBCThreeArg{2}{m}{(\schwarzIt)} =  \left[   
\spatialSchwarzProjectionOperator_2 \trialBasis_1 \reducedDisplacementThreeArg{1}{1}{(\schwarzIt)} ,\ldots, 
\spatialSchwarzProjectionOperator_2  \trialBasis_1 \reducedDisplacementThreeArg{1}{\numTimeStepsOneArg{1}}{(\schwarzIt)} 
\right] \temporalSchwarzProjectionOperator^m_2.
\end{split}
\end{equation*}
 For clarity, the OpInf-OpInf coupling algorithm is summarized at a relatively high-level in Algorithm \ref{algo:rom-offline-online} in terms of the offline and online stages of the model reduction and coupling procedure. 
 It is straightforward to modify Algorithm \ref{algo:rom-offline-online} to perform OpInf-FOM coupling.

\IncMargin{1.5em}
\begin{algorithm}[h]
	\SetKwInOut{Input}{Input}
	\SetKwInOut{Output}{Output}
    	\SetKwInOut{Offline}{\hspace{-0.5cm} Offline stage}
	\SetKwInOut{Online}{\hspace{-0.5cm} Online stage}
	%\Input{$\bm{\mathcal{D}} = \{(\mathbf{x}_{n}, \mathbf{y}_{n})\}_{n=1}^{N}$ - the set of observations of two r.v.s;\\ $\lambda$ - parameter of independence}
	%\Output{The causal direction}
	%\BlankLine
%\ikt{Feedback/corrections on/to this algorithm are welcome.} \;
 \Offline{}
Perform  O-SAM-based FOM-FOM coupled simulations on two overlapping subdomains $\Omega_1$ and $\Omega_2$\;
Collect snapshot matrices  $\displacementSnapshots_i$, $\dirichletBoundarySnapshots_i$, and $\schwarzBoundarySnapshots_i$ from the above coupled simulation in $\Omega_i$, $i=1,2$.\; 
Compute POD bases $~\Phi_1$ and $~\Phi_2$ in $\Omega_1$ and $\Omega_2$ using the snapshot sets $~U_1$ and $~U_2$. \;
Compute boundary POD bases $\trialBasisDirichlet$, $\trialBasisDirichletTwo$, $\trialBasisSchwarz$, and $\trialBasisSchwarzTwo$ from the snapshots $\dirichletBoundarySnapshots_1,\dirichletBoundarySnapshots_2$, $\schwarzBoundarySnapshots_1$, and $\schwarzBoundarySnapshots_2$. \;   
Assume a functional form for the ROM in each subdomain, e.g., the cubic form \eqref{eq:opinf_rom} \;
Compute  OpInf operators in each subdomain by solving the regularized least-squares minimization problems \eqref{eq:quadratic_opinf} in each subdomain with optimal regularization parameter selection described in Section \ref{sec:regularization}\; 
\BlankLine
\BlankLine 

  \setcounter{AlgoLine}{0} 
\Online{}
Apply the O-SAM coupling iteration procedure in Section \ref{sec:opinf-opinf}, with Schwarz BC transfer                           via the pre-learned boundary $\bar{~B}_1$ and $\bar{~B}_2$ in \eqref{eq:schwarz_opinf_opinf_discrete} \;

	\caption{O-SAM for OpInf-OpInf coupling in solid dynamics and the specific case of two overlapping subdomains advanced forward in time using a Newmark-$\beta$ time-integration scheme.}\label{algo:rom-offline-online}
\end{algorithm}
\DecMargin{1.5em}

\subsubsection{Offline stage: Selection of regularization parameter} \label{sec:regularization}
%\ikt{Can this section be generalized to have a discussion involving both FOM-ROM and ROM-ROM coupling?}
The stability and accuracy of the OpInf ROM is strongly dependent on the value of the regularization parameter, $\lambda$. Various strategies have been proposed to identify the optimal regularization parameter, and here we use an extension of the approach proposed in ~\cite{McQuarrie2021combustion}. In essence, the approach proposed in ~\cite{McQuarrie2021combustion} repeatedly solves the inference problem for various regularization parameters. For each regularization parameter, the resulting OpInf ROM is integrated in time using the same initial condition as the training data. The regularization parameter resulting in the OpInf ROM with the lowest trajectory error (as measured against the training data) is then selected.

The above approach is made more complicated in the present setting as the OpInf ROM is coupled to another subdomain model through the Schwarz Dirichlet boundary condition, and as such cannot be integrated independently. While one could solve the coupled system for each regularization parameter, this would be prohibitively expensive for the case where the coupled subdomain model is a FOM. To mitigate this issue, we employ a subdomain-local training process for the OpInf ROM regularization parameter by interpreting the OpInf ROM as a monolithic problem with boundary conditions taken from training data. That is, for each regularization parameter considered for the OpInf ROM, we advance the ROM according to the OpInf specific portion of~\eqref{eq:schwarz_opinf_discrete} with the Schwarz boundary condition taken from training data.

More specifically, let $\Lambda=\{\lambda_1,\ldots,\lambda_{n_{\lambda}}\}$ denote the candidate set of regularization parameters.  Assume without loss of generality that we are looking to regularize the OpInf model in $\Omega_1$.  For each $\lambda\in\Lambda$, we first solve the regularized OpInf least-squares problem \eqref{eq:opinf_schwarz_inference} using all reduced interior snapshots, acceleration snapshots, reduced essential Dirichlet boundary snapshots $\reducedDirichletBoundarySnapshots$, and reduced Schwarz boundary snapshots $\reducedSchwarzBoundarySnapshots$ available for the OpInf subdomain.  Writing the columns of these matrices as
\[
\bar{~U}=\left[\reducedDisplacementTwoArg{1}{0},\ldots,\reducedDisplacementTwoArg{1}{\numTotalTimeStepsOneArg{1}}\right],
\qquad
\reducedDirichletBoundarySnapshots=\left[\hat{\boldsymbol \chi}^{0},\ldots,\hat{\boldsymbol \chi}^{\numTotalTimeStepsOneArg{1}}\right],
\qquad
\reducedSchwarzBoundarySnapshots=\left[\hat{\boldsymbol s}^{0},\ldots,\hat{\boldsymbol s}^{\numTotalTimeStepsOneArg{1}}\right],
\]
the candidate ROM is then advanced fully discretely over the same time grid as the training data.  The initial reduced displacement, velocity, and acceleration are set to the training values at $j=0$.  For each time step $j=0,\ldots,\numTotalTimeStepsOneArg{1}-1$, the OpInf-specific Newmark-$\beta$ equation in \eqref{eq:schwarz_opinf_discrete} is solved with the boundary data at the new time level prescribed by the stored columns $\hat{\boldsymbol \chi}^{j+1}$ and $\hat{\boldsymbol s}^{j+1}$.  During this rollout, the Schwarz data are not updated by solving the coupled Schwarz problem.  Instead, the reduced Schwarz boundary vector $\hat{\boldsymbol s}^{j+1}$ is frozen to the value stored in the training data, while the reduced essential Dirichlet boundary vector $\hat{\boldsymbol \chi}^{j+1}$ is prescribed from the same training data.  Thus, the regularization selection isolates the subdomain-local OpInf model while still testing it under the same discrete interface forcing history that arises in the Schwarz training simulation.

For a stable rollout, the error assigned to a candidate parameter is the relative trajectory error
\begin{equation}
E(\lambda)=
\frac{\|\bar{~U}_{\lambda}^{\rm ROM}-\bar{~U}\|_F}
{\|\bar{~U}\|_F},
\label{eq:schwarz_regularization_error}
\end{equation}
where $\bar{~U}$ denotes the projected training trajectory and $\bar{~U}_{\lambda}^{\rm ROM}=[\reducedDisplacementTwoArg{1}{0,\rm ROM},\ldots,\reducedDisplacementTwoArg{1}{\numTotalTimeStepsOneArg{1},\rm ROM}]$ denotes the corresponding reduced ROM rollout obtained with regularization parameter $\lambda$.  If any discrete state in the rollout is non-finite, the candidate is assigned a large penalty.  We then select $\lambda^\star=\arg\min_{\lambda\in\Lambda}E(\lambda)$ and refit the OpInf operators with $\lambda^\star$ before using the ROM in the online Schwarz coupling.  If multiple independent training trajectories are available, the same procedure is applied to each trajectory and $E(\lambda)$ is replaced by the average of the relative trajectory errors.  Algorithm~\ref{algo:schwarz-reg-selection} summarizes the approach.

\begin{algorithm}
    \SetKwInOut{Input}{Input}
    \SetKwInOut{Output}{Output}
    \Input{Candidate set $\Lambda$; reduced training snapshots $\{\reducedDisplacementTwoArg{1}{j},\reducedVelocityTwoArg{1}{j},\reducedAccelerationTwoArg{1}{j}\}_{j=0}^{\numTotalTimeStepsOneArg{1}}$; reduced essential Dirichlet boundary columns $\{\hat{\boldsymbol \chi}^{j}\}_{j=0}^{\numTotalTimeStepsOneArg{1}}$; frozen reduced Schwarz boundary columns $\{\hat{\boldsymbol s}^{j}\}_{j=0}^{\numTotalTimeStepsOneArg{1}}$; time step $\timeStepOneArg{1}$.}
    \Output{Selected regularization parameter $\lambda^\star$ and fitted OpInf operators.}
    \ForEach{$\lambda \in \Lambda$}{
        \textcolor{black}{Fit the OpInf operators by solving \eqref{eq:opinf_schwarz_inference} with regularization parameter $\lambda$\;
        Set $\reducedDisplacementTwoArg{1}{0,\rm ROM}=\reducedDisplacementTwoArg{1}{0}$, $\reducedVelocityTwoArg{1}{0,\rm ROM}=\reducedVelocityTwoArg{1}{0}$, and $\reducedAccelerationTwoArg{1}{0,\rm ROM}=\reducedAccelerationTwoArg{1}{0}$\;}
        \For{$j=0,\ldots,\numTotalTimeStepsOneArg{1}-1$}{
            Form the reduced boundary vector at time level $j+1$ from the stored columns $\hat{\boldsymbol \chi}^{j+1}$ and $\hat{\boldsymbol s}^{j+1}$\;
            Solve the OpInf Newmark-$\beta$ equation in \eqref{eq:schwarz_opinf_discrete} for $\reducedDisplacementTwoArg{1}{j+1,\rm ROM}$ using this fixed boundary vector\;
            Update $\reducedVelocityTwoArg{1}{j+1,\rm ROM}$ and $\reducedAccelerationTwoArg{1}{j+1,\rm ROM}$ with the Newmark-$\beta$ formulas\;
        }
       \eIf{the rollout is unstable}{
    Set $E(\lambda)$ to a large penalty\;
}{
  Set $E(\lambda)=\|\bar{~U}_{\lambda}^{\rm ROM}-\bar{~U}\|_F/\|\bar{~U}\|_F$\;
}
    
    }
   \textcolor{black}{Set $\lambda^\star=\arg\min_{\lambda\in\Lambda}E(\lambda)$\;
    Refit \eqref{eq:opinf_schwarz_inference} using $\lambda^\star$\;}
    \caption{\textcolor{black}{Schwarz-aware regularization selection for a subdomain-local OpInf ROMs in $\Omega_1$ (without loss of generality).}}\label{algo:schwarz-reg-selection}
\end{algorithm}

\subsubsection{A note on O-SAM convergence} \label{sec:analysis}

Our previous works \cite{Mota:2017} and \cite{Mota:2022} performed a rigorous convergence analysis when applying O-SAM to solid mechanics problems in which each subdomain problem is solved exactly.  Geometric convergence was proven provided the underlying initial boundary value problem is well-posed, and the overlap region is non-empty \cite{Mota:2017}; for solid dynamics, it was also shown that convergence requires a small enough time-step within the Newmark-$\beta$ time-integration scheme used to advance the discretized problem forward in time \cite{Mota:2022}. 

%The convergence results discussed above apply to the classical Schwarz method in which each subdomain problem is solved exactly. 
Since, in the present work, high-fidelity subdomain solution operators are replaced by approximations obtained through Operator Inference, our past convergence analysis \cite{Mota:2017, Mota:2022} does not apply directly.  However, we remark that, when used to couple OpInf ROMs, the O-SAM iteration may be viewed as an inexact Schwarz method.
From this perspective, the convergence behavior is expected to depend on the accuracy of the ROM approximation, i.e., how well the ROM solution represents progress along a ``descent direction" for the functional being minimized at each time step. When the ROM solution operators accurately reproduce the corresponding high-fidelity subdomain responses, the O-SAM iteration is expected to converge similarly to the case of FOM-FOM Schwarz, with a fixed point that is not far from the true solution for in-distribution initial conditions.
% closely mimic the convergence behavior of the underlying FOM-FOM Schwarz method. 
As the ROM approximation error increases, however, the interface data exchanged between subdomains become increasingly perturbed, which may degrade the convergence rate and, in extreme cases, compromise convergence.

The numerical results presented in Section \ref{sec:results} suggest that the OpInf ROMs considered in this work are sufficiently accurate for the Schwarz-ROM iterations to remain stable and to produce solutions with small errors relative to the FOM-FOM reference solutions.  A rigorous analysis quantifying the relationship between ROM approximation errors and Schwarz convergence is beyond the scope of the present work, but remains an important topic for 
% a future publication.
future work.

%% file: 06-results.tex
% Numerical results
\section{Numerical results} \label{sec:results}

%\todo{Make sure we say what is IC for disp and velo for all problems!  I think it's missing for some of them now.}

%\ikt{I am thinking to focus on Norma results to differentiate from Ian's paper. What nonlinear results do we want to show? Would showing torsion problem undermine Eric's paper?}

%\ikt{Candidate test cases:
%\begin{itemize}
%\item Clamped problem - verification 
%\item Notched cylinder - can't get to work... worked before with April version of Norma. Need to look at a bit more.
%\item Laser weld - can't get to work... worked before with April version of Norma. Need to look at a bit more.
%\item Flapping problems (save for SP paper?).
%\item Torsion problem (overlap with Eric's paper?).
%\end{itemize}
%}

We now present some numerical results that illustrate the performance of our proposed O-SAM-based couplings involving disparate models, meshes, time integration schemes, and time-steps. Our Schwarz-based coupling approach has been implemented within {\tt Norma.jl}~\cite{Norma.jl}, an open-source\footnote{{\tt Norma.jl} is available on GitHub: \url{https://github.com/sandialabs/Norma.jl}.  For the numerical results presented herein, the following version of {\tt Norma.jl} was utilized: {\tt d11e20114a9aeb9d2427c03a946c89fee729b3d1}.} three-dimensional (3D) Julia-based finite element code designed for the rapid prototyping of algorithms and ideas for domain coupling and contact in solid mechanics. At the time of writing of this paper, the main branch of {\tt Norma.jl} is able to perform: (i) overlapping and non-overlapping coupling of high-fidelity quasistatic~\cite{Mota:2017} and dynamic~\cite{Mota:2022} solid mechanics models in solid mechanics, (ii) overlapping coupling of dynamic non-intrusive linear, quadratic and cubic OpInf models, %\todo{should we say something about structure preservation?} \adg{I added a few words below} 
(iii) overlapping coupling of dynamic non-intrusive NN-based models with prescribed symmetry and definiteness properties \cite{erics_inprep_paper}, and (iv) contact simulations using a recently-proposed non-overlapping SAM-based contact enforcement algorithm~\cite{Mota:2025}. The code additionally supports several elastic and hyper-elastic material models, including linear elastic (Appendix A.1), Saint Venant--Kirchhoff~\cite{holzapfel2000} (Appendix A.2) and Neohookean~\cite{Mota:2011} (Appendix A.3).  %, and Seth--Hill \ikt{question for Alejandro - is it worth putting a reference for this model?}.  
For time integration, {\tt Norma.jl} employs the usual Newmark-$\beta$ scheme  \cite{Mota2003}, which can be run either implicitly or explicitly.  In the explicit variant of the scheme, mass lumping is employed by default, so as to improve efficiency by avoiding linear solves during the time integration procedure.

The OpInf models coupled herein are trained using a ``top-down" approach, in which training data are generated by running an O-SAM-based FOM-FOM coupled simulation on the same physical geometries of interest. 
%Since our primary objective is to demonstrate our O-SAM-based coupling when mixing a variety of models and meshes, we do not consider bottom-up training~\cite{Chung:2024}, in which subdomain problems are simulated independently of each other, but remark that this would be a very interesting future research task. 
Once training data are generated by performing a coupled FOM-FOM O-SAM simulation using {\tt Norma.jl}, subdomain-local and boundary POD bases are constructed for the vector-valued displacement degrees of freedom and the appropriate OpInf operators are learned offline with the help of the open-source\footnote{{\tt norma-opinf} is available on GitHub: \url{https://github.com/sandialabs/norma-opinf}.  For the numerical results presented herein, the following version of {\tt norma-opinf} was utilized: {\tt 763e37198e197919f4f78dcee4b26efb4809c1cd}.} {\tt norma-opinf} Python package~\cite{norma_opinf} (Algorithm \ref{algo:rom-offline-online}, steps 3--4). Finally, an online coupled ROM-ROM or FOM-ROM simulation (Section \ref{sec:schwarz_rom}) is performed by running {\tt Norma.jl} once again, after feeding it the learned operators and bases obtained from {\tt norma-opinf}. In general, the POD basis dimension is selected based on an energy criterion, so that the modes retained capture some (large) percentage of the total variance in the snapshot set. 
That is, we define the POD energy as 
\begin{equation} \label{eq:pod_energy}
   E_r = \frac{\sum_{i=1}^{r}\sigma_i^2}{\sum_{i=1}^{N}\sigma_i^2} \in [0,1],
\end{equation}
where $\sigma_i$ are the ordered singular values of the snapshot matrix (see Section \ref{sec:pod}), and select the reduced dimension $r$ such that $E_r < \delta_E$, where $\delta_E \in (0,1)$.  The quantity $100E_r$ is typically referred to as the snapshot energy percentage captured by a given reduced basis $~\Phi_r$.  
In addition to representing the solution living in the subdomain interiors using a POD basis, we also construct separate POD bases for the DoFs at which either Schwarz or system Dirichlet boundary conditions are imposed to reduce further the online evaluation costs of the ROMs, as discussed earlier in Section \ref{sec:schwarz_opinf_offline}.  

All numerical experiments discussed in this paper were performed on a Linux RHEL9 cluster known as {\tt Rigel} and located at Sandia National Laboratories, which has two AMD EPYC 9634 (``Genoa") 84-core processors and 1.5 TB DDR5 ECC RAM. Since {\tt Norma.jl} does not, at the present time, possess MPI parallelism, all simulations were performed on one core. While {\tt Norma.jl} has hyper-threading capabilities, hyper-threading was not utilized in the experiments reported on herein. All test cases considered, including the 1D linear elastic wave propagation discussed in Section~\ref{sec:clamped}, were run as three-dimensional problems. Solution accuracy was assessed in terms of the mean square relative error within a given subdomain $\Omega$ for a given solution field $~u$ with respect to a reference solution $~u_{\text{ref}}$, defined as: 
\begin{equation} \label{eq:mse}
    \mathcal{E}_{\Omega}(~u):= \frac{\sqrt{\sum_{k=0}^{\numTime}|| ~u(t_k) - ~u_{\text{ref}}(t_k)||_2^2}}{\sqrt{\sum_{k=0}^{\numTime} || ~u_{\text{ref}}(t_k)||_2^2}}, 
\end{equation}
where $~u(t_k)$ is the solution at time $t_k$ (and similarly for $~u_{\text{ref}}(t_k)$) and $\numTime$ is the total number of time intervals considered. %\adg{we need to uniformize the notation for this. Right now it is $\tau, n_c$ and $S$}. 
We took as the reference solution $~u_{\text{ref}}(t)$ either the exact analytical solution (for the problem in Section~\ref{sec:clamped}) or the FOM solution in a given subdomain $\Omega_i$ obtained by performing a FOM-FOM coupling using O-SAM (for the problems in Sections~\ref{sec:bolted-joint}--\ref{sec:tension-specimen}). We did not perform comparisons with respect to a monolithic solution computed on the full domain $\Omega$, as it is often very difficult (and sometimes impossible) to robustly mesh a full geometry of practical interest (e.g., the bolted joint problem considered in Section \ref{sec:bolted-joint}).
For the 1D linear elastic wave propagation problem considered in Section~\ref{sec:clamped}, we also calculated projection errors associated with a given POD basis $~\Phi_r$, domain $\Omega$ and snapshot set $~u$, defined as: 
\begin{equation} \label{eq:projerr}
    e_{\Omega}(~u, ~\Phi_r) := \frac{||~u - ~\Phi_r ~\Phi_r^{\intercal} ~u||_2}{||~u||_2}.
\end{equation}
The projection error $e(~u, ~\Phi_r)$ is a measure of how well the basis $~\Phi_r$ is capable of representing the solution $~u$, and sets a lower bound on the error that it is achievable by a ROM defined in the latent space spanned by $~\Phi_r$. 

We build and assess subdomain-local OpInf ROMs of all three types described in Section~\ref{sec:opinf}: linear, quadratic and cubic, referred to as OpInf, QOpInf and COpInf, respectively. The first two problems, the 1D linear elastic wave propagation problem (Section~\ref{sec:clamped}) and the 3D bolted joint problem (Section~\ref{sec:bolted-joint}) specify material models that give rise to PDEs with linear and cubic nonlinearities, respectively; hence, we assess OpInf- and COpInf-based couplings for these problems, respectively. While the second two problems, the 3D torsion and 3D tension specimen problems, specify material models which give rise to PDEs with generic (non-polynomial) nonlinearities, we construct and assess couplings involving QOpInf-based models, which we demonstrate are reasonable local surrogates having moderate computational complexity. 
We utilize the algorithm detailed in Section \ref{sec:regularization} to find the ``optimal" values of the regularization parameters in the OpInf least-squares minimization problems.  In particular, we search over a logarithmically-spaced grid of regularization parameters ranging from $10^{-10}$ to 1, i.e., we search over $ \lambda = \{   10^{-10 + \Delta \left[j-1 \right] } \}_{j=1}^{n_{\lambda}} $,  with $\Delta = \frac{10}{n_{\lambda} - 1}$; in our experiments we set $n_{\lambda} = 11$.  For convergence of the Schwarz alternating method, we employ both an absolute and a relative convergence criterion based on the displacement as well as the velocity fields, see \eqref{eq:conv_criterion}.

%\eqref{eq:conv_criterion_abs} and \eqref{eq:conv_criterion_rel}.

When it comes to deciding which model (e.g., ROM vs. FOM) to assign to which subdomain, we generally take the strategy reducing the subdomains with the less complex dynamics.  
While, admittedly, smaller online CPU time improvements are anticipated with this approach, we expect that it will generate a coupled model that is more accurate and robust, a necessity for solid mechanics analyses.  
Our model assignment strategy is similar to the one considered in \cite{Gkimisis:2025} but contrary to the one in \cite{Goswami:2025} and \cite{parish2024embedded}, which opts to reduce subdomains discretized with finer meshes in order to achieve the best online speed-ups.

While this work builds on the preliminary studies presented in~\cite{Moore:2024, Rodriguez:2025}, we emphasize that this paper is the \textit{first} to demonstrate domain decomposition-based couplings involving non-intrusive OpInf ROMs using O-SAM for: (i) realistic 3D nonlinear problems in solid mechanics with geometries discretized using different mesh resolutions and element types, (ii) the coupling of disparate time integration schemes with possibly different time-steps, and (iii) the coupling of higher-order QOpInf and COpInf models with each other and with FOMs.  This also differentiates the present work from \cite{Farcas:2023} and \cite{Gkimisis:2025}.

Before showcasing our results, we note that SI units of measurement are employed throughout this paper, unless otherwise noted. All meshes employed in our simulations were generated using the {\tt CUBIT} meshing software~\cite{cubit} developed
at Sandia National Laboratories.  We remind the reader that, as stated in Section \ref{sec:intro}, the domain decompositions we consider in the present work and in general are physically-motivated, meaning that the number of subdomains being coupled will in generally be small (at most 5--10, but usually 2--3).  As discussed in Remark \ref{remark:many_sds}, our method is applicable to an arbitrary number of subdomains.  We additionally remind the reader that only the sequential multiplicative Schwarz variant is considered herein, as detailed in Remark \ref{remark:additive_sam}.  

More information about the code versions used and where to find input decks to reproduce the results presented herein can be found in the ``Code availability and reproducibility" section at the end of this paper.

\input{clamped}

\input{bolted-joint}

\input{torsion}

\input{tension-specimen}

%% file: clamped.tex
\subsection{1D linear elastic wave propagation problem} \label{sec:clamped}

%\ikt{There is a lot in this section. I am interested in folks' opinions on whether it is all needed. Maybe 2 reproductive variants of this problem is overkill?  If we wish to cut, my suggestion would be to show reproductive Rounded Square and predictive Symmetric Gaussian - thoughts?} \irm{It is definitely a lot. If we want 2, then I would say repro/predictive Gaussian. I'm not sure what the Rounded Square is giving besides proof that the result works for multiple ICs, but I think we have enough examples to demonstrate that the method is general. The 1D section also currently takes about 10 pages, and I think cutting the Rounded Square would be the most space saving, which would help the reader focus on the 3D results. Our other problems mostly follow a format of reproductive/predictive as well.} \adg{After reading carefully, I actually like having all 3 in the main body. I think the length can be mitigated with some Figure engineering\textemdash I've left comments below.}

The first test case we consider is the so-called 1D linear elastic wave propagation problem, variants of which can be found in~\cite{Mota:2022, Barnett:2022, Rodriguez:2025}. The purpose of this example is to verify our method's ability to couple disparate subdomain-local models (FOMs and OpInf ROMs) discretized using different time integration schemes (implicit and explicit Newmark-$\beta$) having potentially different time-steps.

Consider a simple beam geometry having a length of 1 m in the $z$--dimension and a cross-sectional area of 1 $\times$ 1 mm$^2$, so that 
% $\Omega = (-5.0\times 10^{-4}, 5.0\times 10^{-4}) \times (-5.0 \times 10^{-4}, 5.0 \times 10^{-4}) \times (-0.5, 0.5)$
$\Omega = (-5.0\times 10^{-4}, 5.0\times 10^{-4})^2 \times (-0.5, 0.5)$. The objective is to specify a 1D problem using the 3D {\tt Norma.jl} code. Toward this effect, we set a homogeneous Dirichlet boundary condition on the $x$-- and $y$--displacement at the 
% $x=-5.0\times 10^{-4}$, $x=5.0 \times 10^{-4}$, $y=-5.0 \times 10^{-4}$ and $y=5.0 \times 10^{-4}$ 
$x,y = \pm 5.0\times 10^{-4}$ boundaries. We  additionally assume that the beam is clamped at the ends in the $z$--dimension, which translates to a homogeneous Dirichlet boundary condition on the $z$--displacement at the $z = -0.5$ and $z = 0.5$ boundaries. We prescribe within $\Omega$ a simple linear elastic material model (Appendix A.1) with Young's modulus $E=1$ GPa, density $\rho = 1000$ kg/m$^3$ and Poisson's ratio $\nu = 0$. We initialize the problem by specifying an initial displacement of the form:
\begin{equation}
    ~u(~x, 0) = \left( \begin{array}{ccc} 0, & 0, & f(z) \end{array}\right)^{\intercal},
\end{equation}
for a specified function $f(z)$ with $z \in \Omega$, and a zero initial velocity $~v(~x, 0) = ~0$, where $~x:=\left( \begin{array}{ccc}x, & y, & z\end{array}\right)^{\intercal}$ is the coordinate vector. We will assume the problem is run from time $t=0$ to time $t = \timeFinal = 1.0 \times 10^{-3}$ s. It is straightforward to show using the method of characteristics that the exact analytical displacement solution to the problem described above is 
\begin{equation} \label{eq:disp_exact}
    ~u_{\text{ref}}(~x, t) = \left( \begin{array}{ccc}
    0, & 0, & \frac{1}{2}f(z-ct) + \frac{1}{2} f(z+ct) - \frac{1}{2} f(z-c(T-t)) - \frac{1}{2}\textcolor{blue}{f}(z+c(T-t)) 
    \end{array}\right)^{\intercal},
\end{equation}
where $c:=\sqrt{E/\rho}$ is the speed of sound ~\cite{Mota:2022}. The exact analytical solution for the velocity can be derived by differentiating~\eqref{eq:disp_exact} in time, and is not given here for the sake of brevity.

We will consider herein two initial conditions for our 1D linear elastic wave propagation problem, termed the ``Symmetric Gaussian" and the ``Rounded Square" initial condition, respectively, and summarized in Table~\ref{tab:clamped-ics}. As can be seen by comparing the solutions plotted in Figures~\ref{fig:clamped-rounded-square-solns-repro} and~\ref{fig:clamped-gaussian-solns-predi}, the Rounded Square variant of the problem gives rise to a solution with sharper gradients, which are more difficult to resolve, especially using data-driven models such as OpInf ROMs.

\begin{table}[ht!]
    \centering
    \caption{1D linear elastic wave propagation problem: initial conditions considered.}
    \label{tab:clamped-ics}
    \begin{tabular}{c|c|c|c|c}
    Initial Condition & $f(z)$ & $a$ & $b$ & $s$  \\
    \hline \hline 
    Symmetric Gaussian & $a \exp\left( -\frac{(z-b)^2}{2s^2}\right)$  & $1.0 \times 10^{-3}$ & $0$&  $2.0 \times 10^{-2}$  \\
    \hline 
    \multirow{2}{*}{Rounded Square} & $a \tanh(-b(z+0.5-s)) $ & \multirow{2}{*}{$5.0 \times 10^{-4}$} & \multirow{2}{*}{$100$} & \multirow{2}{*}{$0.6$}  \\
     & + $a \tanh(b(z-0.5 + s))$ & & &
    \end{tabular}
\end{table}

The first step in applying O-SAM to this problem is to generate a domain decomposition of the physical geometry $\Omega$ into overlapping subdomains $\Omega_i$, and to define spatial as well as temporal discretizations of the $\Omega_i$. Toward this effect, we will consider a simple decomposition of $\Omega$ into two subdomains, $\Omega_1$ and $\Omega_2$, where 
% $\Omega_1 = (-5.0\times 10^{-4}, 5.0\times 10^{-4}) \times (-5.0 \times 10^{-4}, 5.0 \times 10^{-4}) \times (-0.5, 0.25)$ and $\Omega_2 = (-5.0\times 10^{-4}, 5.0\times 10^{-4}) \times (-5.0 \times 10^{-4}, 5.0 \times 10^{-4}) \times (-0.25, 0.5)$
$\Omega_1 = (-5.0\times 10^{-4}, 5.0\times 10^{-4})^2 \times (-0.5, 0.25)$ and $\Omega_2 = (-5.0\times 10^{-4}, 5.0\times 10^{-4})^2 \times (-0.25, 0.5)$, so that the overlap region $\Omega_1 \cap \Omega_2$ has a length of $0.5$ m in the $z$--dimension. Since the problem considered in this section is effectively 1D, we will use uniform conformal hexahedral discretizations of $\Omega_1$ and $\Omega_2$ having spatial increments $\Delta x = \Delta y = \Delta z = 1.0 \times 10^{-3}$ m, but will focus on assessments involving various combinations of subdomain-local FOMs and OpInf ROMs, as well as time-steppers and time steps (see Table~\ref{tab:clamped-time-steppers}). All time-steppers utilized are of the Newmark-$\beta$ type, with $\beta = 0.25$, $\gamma=0.5$ for the implicit variant, and $\beta = 0$, $\gamma = 0.5$ for the explicit variant. It is straightforward to show that the stable time step according to the Courant--Friedrichs--Levy (CFL) condition for the explicit Newmark-$\beta$ stepper is $\Delta t = \Delta z / c = 1.0\times 10^{-6}$ s. The reader can observe from Table~\ref{tab:clamped-time-steppers} that all time steps employed are approximately one order of magnitude smaller than this value, so that stability and accuracy are assured. A controller time-step of $1.0 \times 10^{-7}$ s was employed for all runs.
%From this point forward, we will denote the time step in subdomain $\Omega_i$ by $\Delta t_i$.

%\begin{table}[ht!]
%    \centering
%    \caption{1D linear elastic wave propagation problem: models considered. \adg{I  recommend putting this alongside Table 3.}}
%    \label{tab:clamped-models}
%    \begin{tabular}{c|c}
    %\multicolumn{2}{c}{Time-steppe} & \multicolumn{2}{c}{Time step (s)}  \\
    % $\Omega_1$ & $\Omega_2 $ \\
%    \hline \hline 
%   FOM & FOM \\
%   \hline 
%   FOM & OpInf \\
%   \hline
%   OpInf & OpInf 
%    \end{tabular}
%\end{table}

\begin{table}[ht!]
    \centering
    \caption{1D linear elastic wave propagation problem: models, Newmark-$\beta$ time discretizations and time steps considered for the couplings evaluated herein.}
    \label{tab:clamped-time-steppers}
     \begin{minipage}{0.3\textwidth}
        \centering
          \begin{tabular}{c|c}
    %\multicolumn{2}{c}{Time-steppe} & \multicolumn{2}{c}{Time step (s)}  \\
    $\Omega_1$ & $\Omega_2 $ \\
    \hline \hline 
   FOM & FOM \\
   \hline 
   FOM & OpInf \\
   \hline
   OpInf & OpInf 
    \end{tabular}
        \end{minipage}
%\hspace{1cm}
 \begin{minipage}{0.45\textwidth}
        \centering
        \begin{tabular}{c|c||c|c}
    \multicolumn{2}{c}{Time-stepper} & \multicolumn{2}{c}{Time step (s)}  \\
    $\Omega_1$ & $\Omega_2 $ & $\Omega_1$ & $\Omega_2$ \\
    \hline \hline 
   Implicit & Implicit & \multicolumn{2}{c}{$1.0 \times 10^{-7}$}\\
   \hline 
   Explicit & Implicit & \multicolumn{2}{c}{$1.0 \times 10^{-7}$}\\
   \hline
   Explicit & Explicit & \multicolumn{2}{c}{$1.0 \times 10^{-7}$}\\
   \hline
   Explicit & Implicit & $1.0 \times 10^{-7}$ & $2.0 \times 10^{-7}$
    \end{tabular}
    \end{minipage}
\end{table}

For couplings involving OpInf ROMs, we consider three different training/testing scenarios: 
\begin{itemize}
    \item \textit{Scenario 1: reproductive Symmetric Gaussian.} Here, we train the OpInf ROM(s) by simulating the 1D linear elastic wave propagation problem with the Symmetric Gaussian initial condition, and assess the coupled models' ability to reproduce the problem solution with the same initial condition.
    \item \textit{Scenario 2: reproductive Rounded Square.} We next train the OpInf ROM(s) by simulating the 1D linear elastic wave propagation problem with the Rounded Square initial condition, and again evaluate the coupled models' skill in reproducing the problem solution with this initial condition. 
    \item \textit{Scenario 3: predictive Symmetric Gaussian.} Finally, we train the OpInf ROM(s) using solution data from the Rounded Square variant of the 1D linear elastic wave propagation problem, and study the coupled models' ability to predict the solution to the problem with the Symmetric Gaussian initial condition.
\end{itemize}
Although we performed FOM-FOM couplings involving all four time discretizations summarized in Table~\ref{tab:clamped-time-steppers}, we trained all the OpInf models using snapshot data obtained by running the Implicit-Implicit FOM-FOM. 

As mentioned earlier, errors for all the coupled models for the 1D linear elastic wave propagation problem were calculated with respect to the exact analytical solution, derived from \eqref{eq:disp_exact}. Since the problem is effectively one-dimensional and the meshes used to discretize the subdomains $\Omega_i$ are conformal, it is straightforward to obtain a single-domain solution from our O-SAM coupled solutions by simply averaging the solution values in the overlap region and applying the mean square error formula \eqref{eq:mse}.  This allows us to report a single error value, corresponding to the error in the full domain $\Omega$, in assessing the accuracy of our models.  
%We report such averaged errors in the plots and analysis below. 
For all O-SAM coupled runs, very tight Schwarz tolerances of $\delta_{\text{abs}} = 1.0 \times 10^{-8}$ and $\delta_{\text{rel}} = 1.0 \times 10^{-12}$ were utilized. For the implicit Newmark-$\beta$ runs, the relative and absolute linear solver tolerances were $1.0 \times 10^{-10}$ and $1.0 \times 10^{-6}$, respectively.

\begin{figure}[ht!]
    \centering
     \begin{subfigure}{0.49\linewidth}
\includegraphics[width=0.99\linewidth]
    {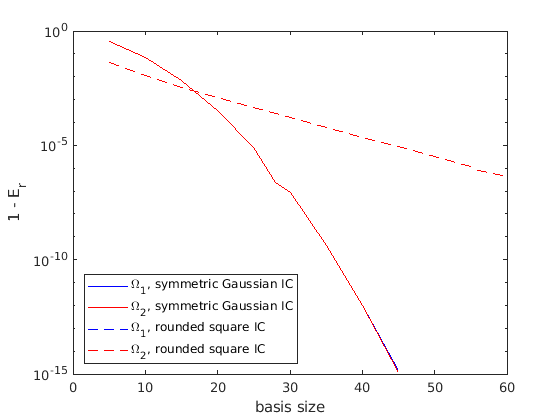}
     \subcaption{Singular value decay ($1-E_r)$}
    \end{subfigure}
    \begin{subfigure}{0.49\linewidth}
\includegraphics[width=0.99\linewidth]{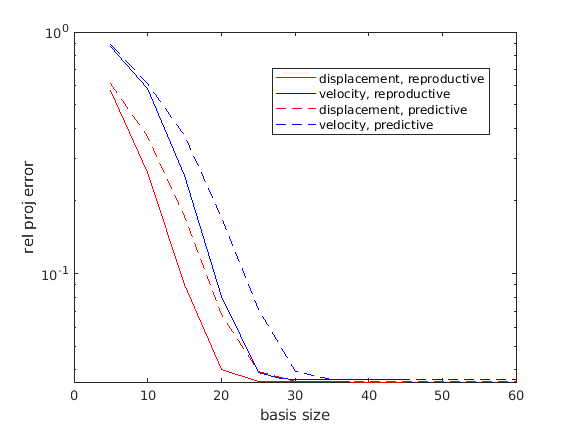}
        \subcaption{Projection errors, Symmetric Gaussian variant}
    \end{subfigure}
       \begin{subfigure}{0.49\linewidth}    \includegraphics[width=0.99\linewidth]{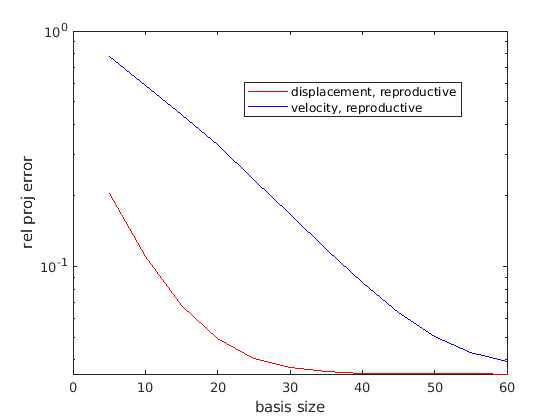}
        \subcaption{Projection errors, Rounded Square variant}
    \end{subfigure}
    
      \vspace{0.5cm}
    \caption{1D linear elastic wave propagation problem: singular value decay (a) and  projection errors calculated in $\Omega_2$ for the Symmetric Gaussian (b) and Rounded Square (c) problem variant, as functions of the POD basis size. } 
    \label{fig:clamped-proj-errs}
\end{figure}

Before presenting results for the three scenarios described above and the various couplings summarized in Table~\ref{tab:clamped-time-steppers}, we perform some diagnostic assessments to determine how accurate we expect our O-SAM coupled models involving OpInf ROMs to be. Figure~\ref{fig:clamped-proj-errs}(a) shows the singular value decay as a function of the POD basis size for the Symmetric Gaussian and Rounded Square variants of the 1D linear elastic wave propagation problem. The reader can observe that significantly more POD modes are needed to capture the same fraction of the snapshot energy for the Rounded Square version of this problem than for the Symmetric Gaussian version, as expected: to capture 99.9999\% of the snapshot energy, 56 POD modes are required for the former, compared to 28 modes for the latter. Moreover, the singular value decay is much slower for the Rounded Square problem variant. The fact that the singular value decay in $\Omega_2$ is identical to that in $\Omega_1$ is expected, given the symmetry of the problem solutions \eqref{eq:disp_exact} with respect to $z = 0$. For all test cases reported on here, we select boundary bases based on a 99.9999\% energy criterion, which gives rise to bases consisting of just one POD mode per boundary.

%\begin{figure}[ht!]
%    \centering
   % \begin{subfigure}{0.49\linewidth}
%\includegraphics[width=0.7\linewidth]{figures/clamped/clamped-energies.png}
     %   \subcaption{Displacement, FOM-OpInf}
   % \end{subfigure}
%      \vspace{0.5cm}
%    \caption{1D linear elastic wave propagation problem: POD singular value decay.
%    \todo{y-axis is not right I think. \adg{Could this figure be combined with Figure 2?}}} 
%    \label{fig:clamped-sv-decay}
%\end{figure}

Having examined the POD singular value decay, we now look at projection errors for the three problem scenarios described above, calculated using \eqref{eq:projerr} for $\Omega_2$ (without loss of generality) and reported in Figures~\ref{fig:clamped-proj-errs}(b)--(c). It can be seen from this figure that approximately 30 modes are sufficient to capture the displacement field for all three problem variants to an error of $\mathcal{O}(1\%)$, including the predictive one. A greater number of modes is needed achieve the same projection error for the velocity field, especially for the Rounded Square version of this problem which exhibits a steep gradient (Figure~\ref{fig:clamped-proj-errs}(c)).  This is expected for two reasons: (i) the velocity field is generally harder to capture, and (ii) the reduced basis $~\Phi_r$ is built from snapshots of only the displacement field, as discussed in Remark \ref{remark1}.  

%\crw{How is the projected velocity being computed? The basis is only for the displacement, right?} \ikt{Yes, exactly.  Figure \ref{fig:clamped-proj-errs}(c) is showing the projection errors for the velocity field as represented using the basis computed from displacement snapshots.  I made that clear above.  In terms of how velocity is obtained, it is as a part of the Newmark time-stepping scheme.  It comes for free there.  I didn't say that anywhere because I thought it was obvious, but I can if you think readers will be confused.} %\crw{I am personally confused on the latter part, usually projection error is independent of any time integrator, just snapshots and a basis.}  \ikt{Projection error looks at how well a basis can capture a field.  The field can be anything, whether or not it is even from the PDE.  Here, I am checking how well the basis represents the velocity, as the velocity is also represented with the same POD basis as the displacement.  This is independent of the time integration scheme.}

\subsubsection{Scenario 1: reproductive Symmetric Gaussian} \label{sec:clamped_scenario1}

\begin{figure}[ht!]
    \centering
     \begin{subfigure}{0.99\linewidth}
\hspace{6.5cm}\includegraphics[width=0.2\linewidth]{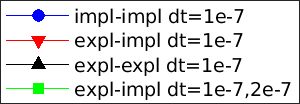}
    \end{subfigure}
    \begin{subfigure}{0.49\linewidth}
\includegraphics[width=0.99\linewidth]{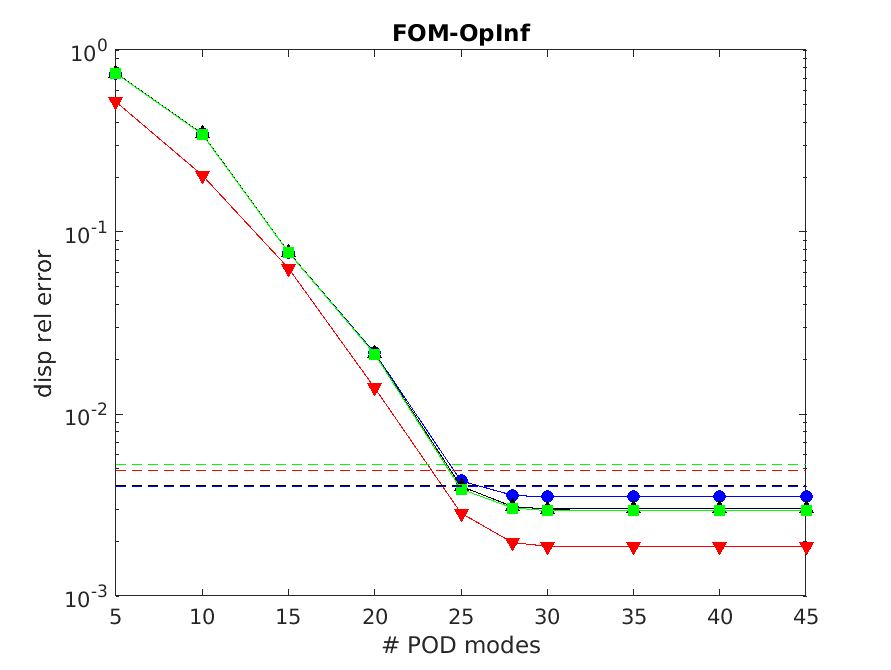}
        \subcaption{Displacement, FOM-OpInf}
    \end{subfigure}
       \begin{subfigure}{0.49\linewidth}    \includegraphics[width=0.99\linewidth]{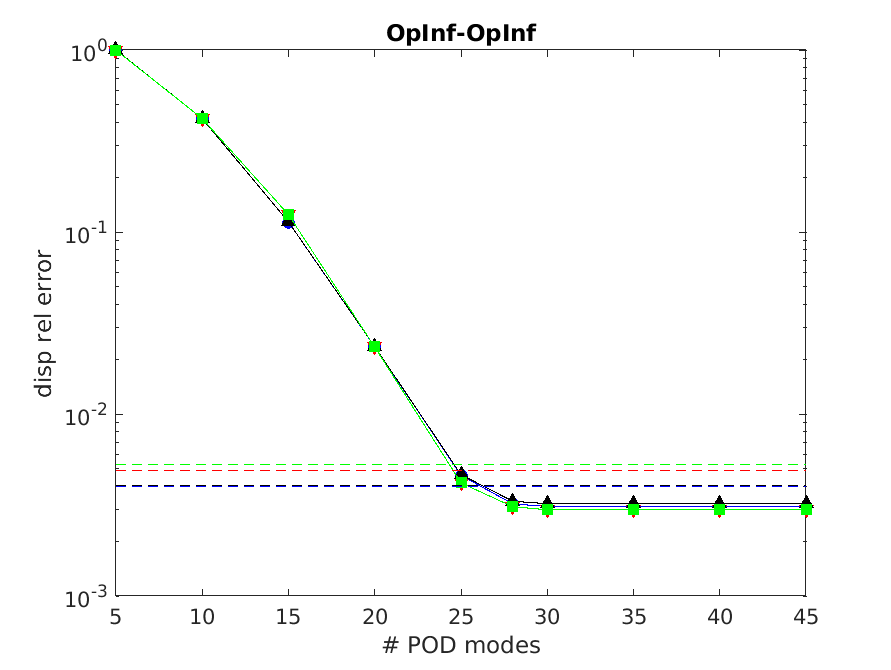}
        \subcaption{Displacement, OpInf-OpInf}
    \end{subfigure}
       \begin{subfigure}{0.49\linewidth}
\includegraphics[width=0.99\linewidth]{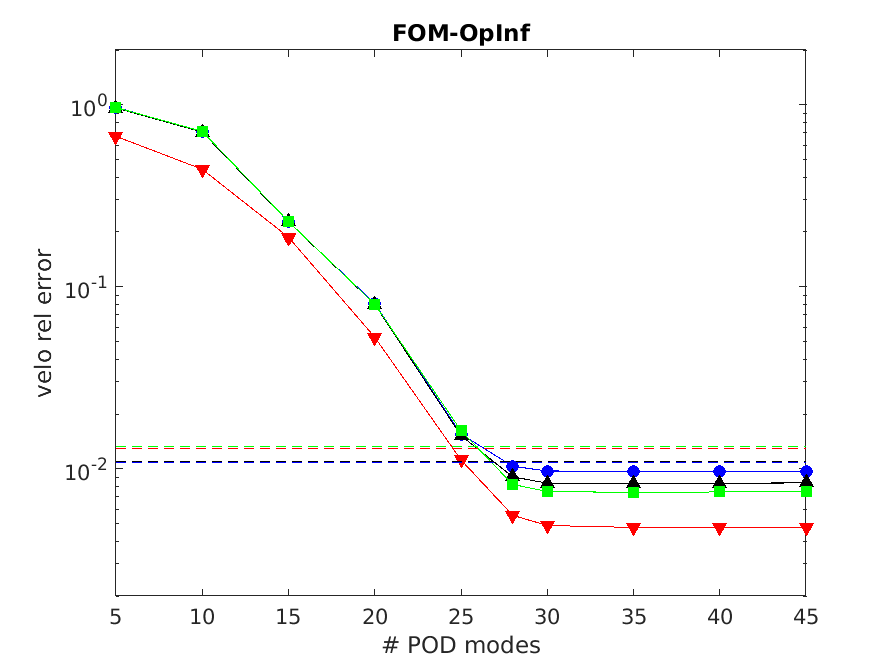}
        \subcaption{Velocity, FOM-OpInf}
    \end{subfigure}
       \begin{subfigure}{0.49\linewidth}    \includegraphics[width=0.99\linewidth]{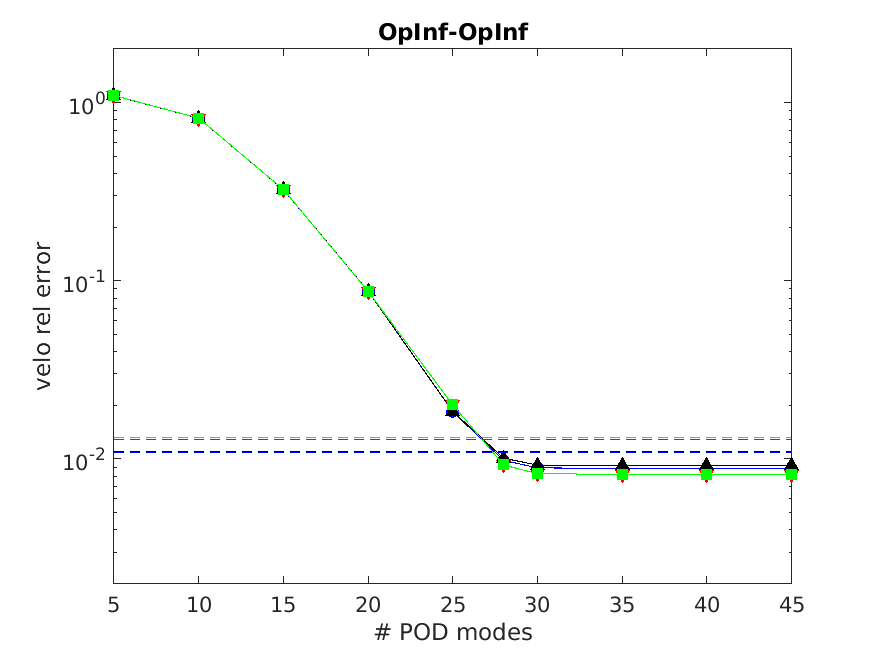}
        \subcaption{Velocity, OpInf-OpInf}
    \end{subfigure}
     %  \begin{subfigure}{0.49\linewidth}
%\includegraphics[width=0.99\linewidth]{figures/clamped/gaussian-fom-opinf-acce-repro.png}
 %       \subcaption{Acceleration, FOM-OpInf}
%       \begin{subfigure}{0.49\linewidth}    \includegraphics[width=0.99\linewidth]{figures/clamped/gaussian-opinf-opinf-acce-repro.png}
   %     \subcaption{Acceleration, OpInf-OpInf}
   % \end{subfigure}
    
      \vspace{0.5cm}
    \caption{1D linear elastic wave propagation problem, Symmetric Gaussian initial condition, reproductive regime: displacement (top row) and velocity (bottom row) relative errors with respect to the exact analytical solution for various FOM-OpInf (a) and OpInf-OpInf (b) O-SAM couplings. Dashed horizontal lines show relative errors for FOM-FOM O-SAM couplings with the colors designated in the legend. %\crw{I would remove figure titles if the subcaptions will already indicate the plot type}  \ikt{I see your point, but I am trying to anticipate a very lazy reader who is just skimming and doesn't want to read lengthy figure captions.  I am keeping it for now therefore.}
    } 
    \label{fig:clamped-gaussian-errors-repro}
\end{figure}

We first discuss results for Scenario 1, in which the Symmetric Gaussian variant of the 1D linear elastic wave propagation problem is run in the reproductive regime. Figure~\ref{fig:clamped-gaussian-errors-repro} is a convergence plot showing the displacement and velocity relative errors with respect to the exact analytic solution \eqref{eq:disp_exact} as a function of the number of POD modes for a variety of FOM-OpInf and OpInf-OpInf couplings. For comparison purposes, relative errors achieved by various analogous FOM-FOM couplings are indicated with dashed horizontal lines. It is interesting to observe that, for basis sizes larger than approximately 25 modes, all FOM-OpInf and OpInf-OpInf couplings deliver solutions that are actually \textit{more} accurate than analogous FOM-FOM couplings. We attribute this result to the fact that the couplings involving OpInf models rely on shape functions that are derived from solutions of the underlying PDE being solved, whereas the FOM-FOM couplings rely on generic and problem-agnostic finite element shape functions. The use of specialized problem-specific shape functions, namely POD modes, seems to enable our FOM-OpInf and OpInf-OpInf models to achieve a sort of super-convergence. The plateauing of the errors for large enough basis sizes can be attributed to the closure error inherent in the OpInf process, which arrests convergence in the pre-asymptotic regime.

Additionally, the reader can observe by examining Figures~\ref{fig:clamped-gaussian-errors-repro}(a) and (c) that the Explicit-Implicit FOM-OpInf coupling for which the same time-step is employed in each subdomain delivers the most accurate solution, whereas the Implicit-Implicit FOM-OpInf coupling is the least accurate. This trend is different from the trend observed for analogous FOM-FOM couplings, for which the Implicit-Implicit model is the clear ``winner" in terms of accuracy. Moreover, employing different time-steps in the two subdomains does not lead to the same decline in accuracy for the FOM-OpInf coupled models as it does for analogous FOM-FOM coupled models. From Figures~\ref{fig:clamped-gaussian-errors-repro}(b) and (d), it can be seen that all OpInf-OpInf coupled models, regardless of the time integration scheme or time step employed, deliver roughly the same accuracy for all basis sizes considered. This accuracy is comparable to that of all but one of the FOM-OpInf coupled models.

\begin{figure}[ht!]
    \centering
     \begin{subfigure}{0.99\linewidth}
\hspace{6.5cm}\includegraphics[width=0.2\linewidth]{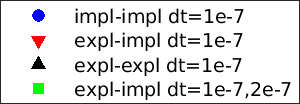}
    \end{subfigure}
    \begin{subfigure}{0.49\linewidth}
\includegraphics[width=0.99\linewidth]{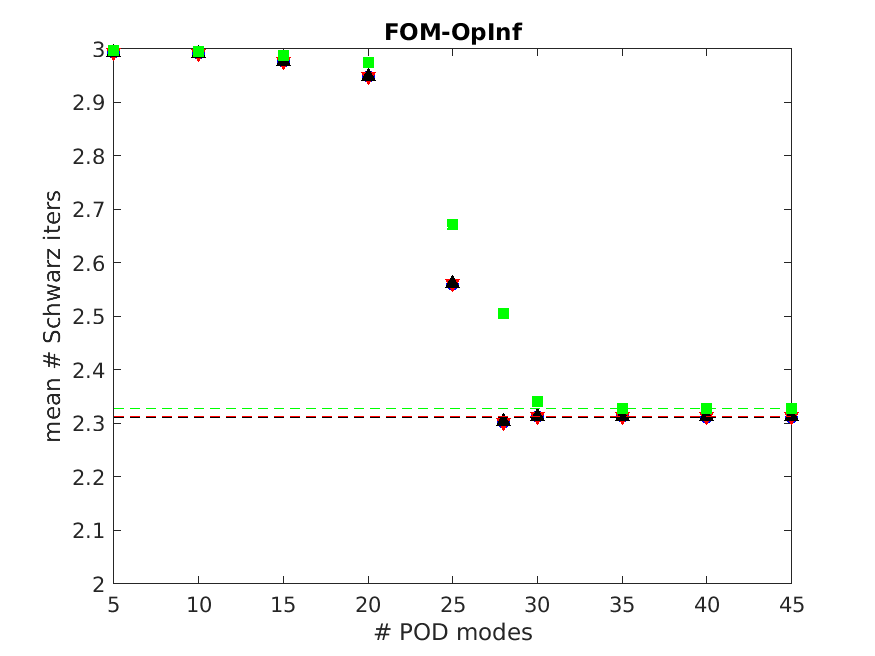}
        \subcaption{FOM-OpInf}
    \end{subfigure}
       \begin{subfigure}{0.49\linewidth}    \includegraphics[width=0.99\linewidth]{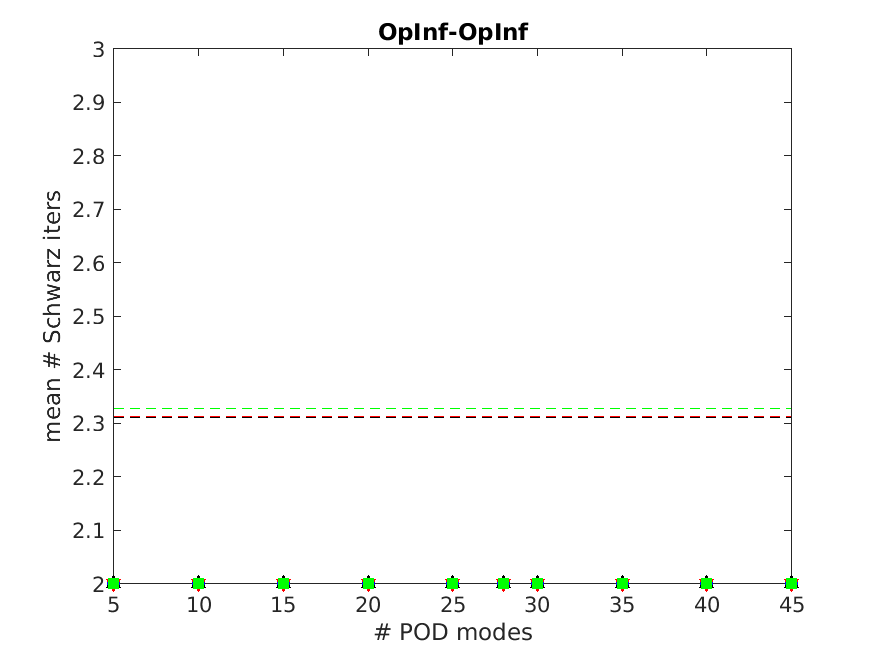}
        \subcaption{OpInf-OpInf}
    \end{subfigure}
        
      \vspace{0.5cm}
    \caption{1D linear elastic wave propagation problem, Symmetric Gaussian initial condition, reproductive regime: mean number of Schwarz iterations required to reach convergence for various FOM-OpInf (a) and OpInf-OpInf (b) O-SAM couplings. Dashed horizontal lines show the number of Schwarz iterations needed to reach convergence for FOM-FOM O-SAM couplings with the colors designated in the legend. 
    % \adg{Is the OpInf-OpInf number really constant at 2?} \ikt{Yes, believe it or not!  The text says this actually.}
    } 
    \label{fig:clamped-gaussian-schwarz-iters-repro}
\end{figure}

Having studied accuracy, we now turn our attention to convergence and performance. Figure~\ref{fig:clamped-gaussian-schwarz-iters-repro} shows the mean number of Schwarz iterations required to achieve convergence for a variety of FOM-OpInf and OpInf-OpInf coupled models, in comparison to analogous FOM-FOM couplings, indicated by the dashed horizontal lines. The reader can observe from Figure~\ref{fig:clamped-gaussian-schwarz-iters-repro}(a) that, for smaller basis sizes, more Schwarz iterations are required to converge the O-SAM method when performing FOM-OpInf couplings compared to analogous FOM-FOM couplings, but convergence of the former to the latter is observed with basis size refinement. The situation is different for the OpInf-OpInf couplings, however (Figure~\ref{fig:clamped-gaussian-schwarz-iters-repro}(b)): all OpInf-OpInf couplings converge in exactly two Schwarz iterations, compared to an average of $\sim 2.3$ Schwarz iterations for similar FOM-FOM couplings. Similar results were obtained for the other test cases considered in this paper. We believe there are two possible explanations for this behavior. First, since our OpInf models rely on POD modes to approximate the solution, these models give rise to solutions that are inherently smoother than their FEM analogues, which can aid convergence of the O-SAM coupling method. %It is likely that 
Convergence may also be accelerated due to the fact that the shape functions underlying the OpInf-OpInf couplings are problem-specific and data-driven, unlike the problem-agnostic and generic finite element shape functions underlying the FOM-FOM couplings.

\begin{figure}[ht!]
    \centering
   %  \begin{subfigure}{0.99\linewidth}
%\hspace{6.8cm}\includegraphics[width=0.3\linewidth]{figures/clamped/legend-pareto.png}
%    \end{subfigure}
    \begin{subfigure}{0.60\linewidth}
\includegraphics[width=0.99\linewidth]{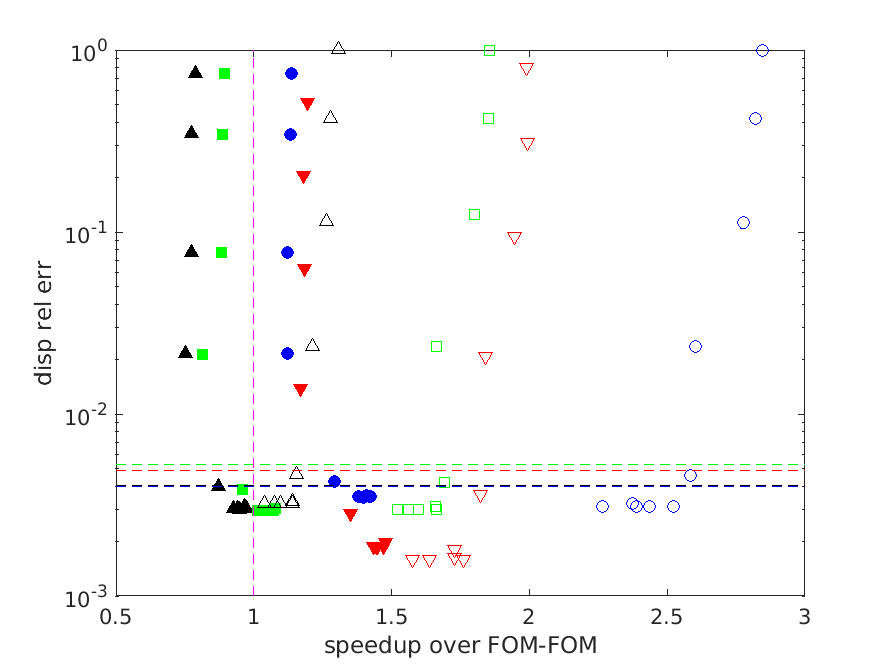}
        %\subcaption{Displacement}
        \end{subfigure}
         \raisebox{2cm}{
       \begin{subfigure}{0.3\linewidth} 
       \includegraphics[width=0.99\linewidth]{  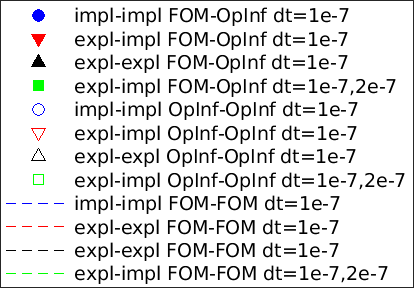}
      % \includegraphics[width=0.99\linewidth]{figures/clamped/gaussian-velo-pareto-repro.png}
       % \subcaption{Velocity}
          \end{subfigure}}
    %      \begin{subfigure}{0.49\linewidth}    \includegraphics[width=0.99\linewidth]{figures/clamped/gaussian-acce-pareto-repro.png}
   %     \subcaption{Acceleration}
   % \end{subfigure}
        
    %  \vspace{0.5cm}
    \caption{1D linear elastic wave propagation problem, Symmetric Gaussian initial condition, reproductive regime: Pareto plot showing the speed-up over an analogous FOM-FOM coupling vs. displacement relative error for various FOM-OpInf (filled symbols) and OpInf-OpInf (unfilled symbols) couplings. Dashed horizontal lines indicate relative errors for corresponding FOM-FOM O-SAM couplings. Dashed magenta vertical line indicates a speedup of 1.} %\adg{The two pareto plots are very similar. Does it make sense to present only one with the legend beside it?}} 
  
    \label{fig:clamped-gaussian-pareto-repro}
\end{figure}

Finally, we assess method performance. Figure~\ref{fig:clamped-gaussian-pareto-repro} shows a Pareto plot for the various couplings considered, in which the speedup over a corresponding FOM-FOM coupling is plotted vs. the displacement  relative errors. We do not show a Pareto plot for the velocity, as the conclusions in terms of relative accuracy and efficiency of the various coupled models evaluated would be the same.  A speed-up of unity is indicated by a vertical dashed magenta line, and the relative errors achieved by the four FOM-FOM couplings considered are indicated by horizontal dashed lines. It follows that all data points falling to the right of the vertical dashed line, and at/below the horizontal lines are competitive from an accuracy and efficiency standpoint. The reader can observe that the majority of the competitive models are OpInf-OpInf couplings. The largest speed-up achieved of $\sim 2.5 \times$ is achieved by the Implicit-Implicit OpInf-OpInf coupled model. As we show below, more substantial speedups are possible for larger, three-dimensional benchmarks, in which larger DoF reductions are possible through the use of subdomain-local OpInf ROMs.

\subsubsection{Scenario 2: reproductive Rounded Square} \label{sec:clamped_scenario2}

\begin{figure}[ht!]
    \centering
     \begin{subfigure}{0.99\linewidth}
\hspace{6.5cm}\includegraphics[width=0.2\linewidth]{figures/clamped/legend.png}
    \end{subfigure}
    \begin{subfigure}{0.49\linewidth}
\includegraphics[width=0.99\linewidth]{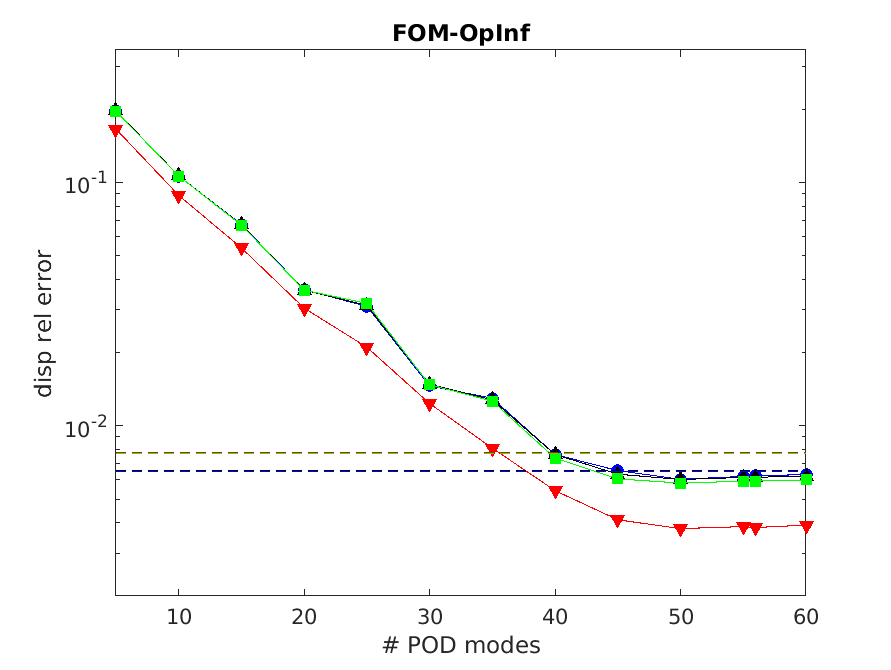}
        \subcaption{Displacement, FOM-OpInf}
    \end{subfigure}
       \begin{subfigure}{0.49\linewidth}    \includegraphics[width=0.99\linewidth]{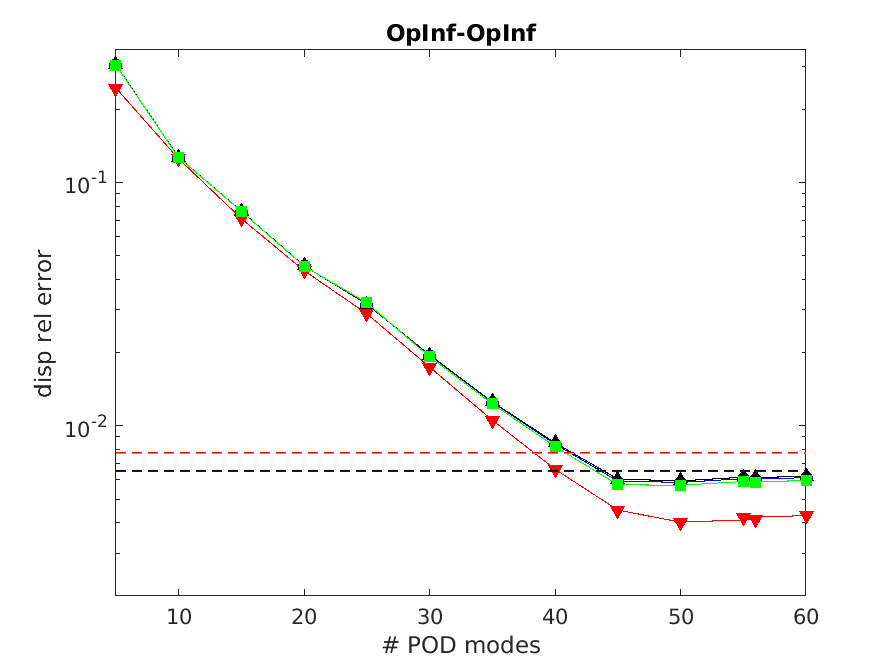}
        \subcaption{Displacement, OpInf-OpInf}
    \end{subfigure}
       \begin{subfigure}{0.49\linewidth}
\includegraphics[width=0.99\linewidth]{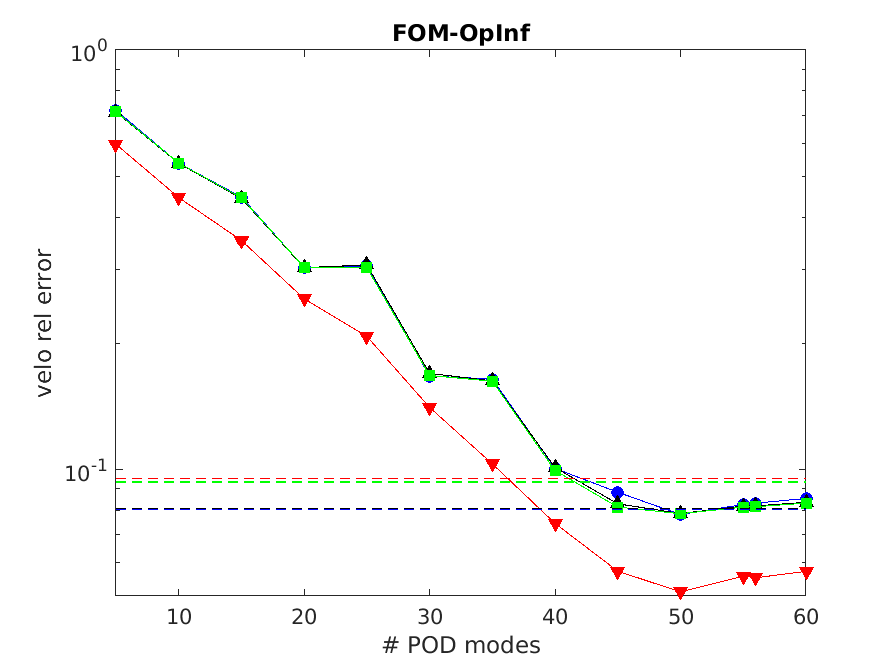}
        \subcaption{Velocity, FOM-OpInf}
    \end{subfigure}
       \begin{subfigure}{0.49\linewidth}    \includegraphics[width=0.99\linewidth]{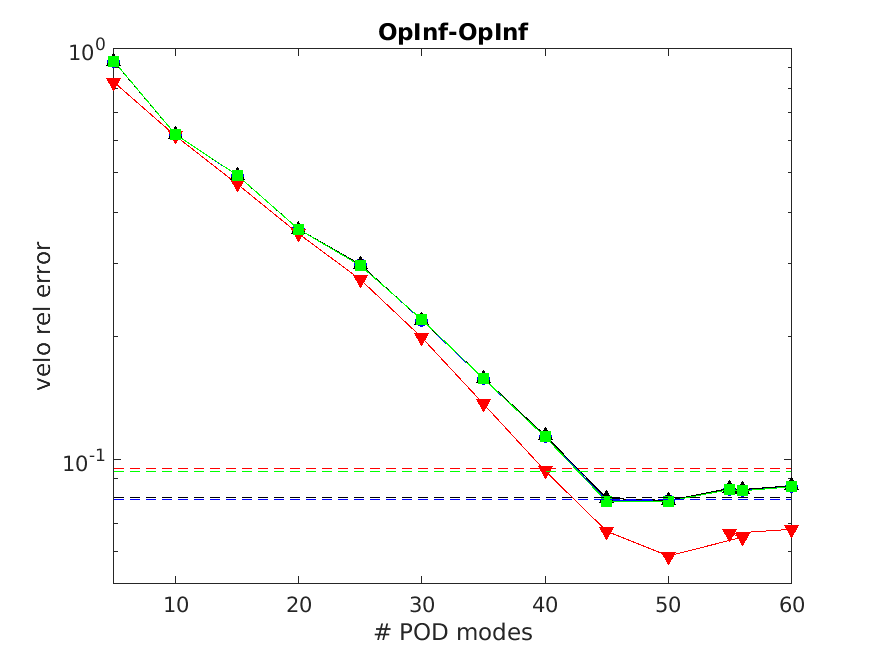}
        \subcaption{Velocity, OpInf-OpInf}
    \end{subfigure}
  %     \begin{subfigure}{0.49\linewidth}
%\includegraphics[width=0.99\linewidth]{figures/clamped/rounded-square-fom-opinf-acce-repro.png}
 %       \subcaption{Acceleration, FOM-OpInf}
 %   \end{subfigure}
  %     \begin{subfigure}{0.49\linewidth}    \includegraphics[width=0.99\linewidth]{figures/clamped/rounded-square-opinf-opinf-acce-repro.png}
 %       \subcaption{Acceleration, OpInf-OpInf}
 %   \end{subfigure}
    
      \vspace{0.5cm}
    \caption{1D linear elastic wave propagation problem, Rounded Square initial condition, reproductive regime: displacement (top row) and velocity (bottom row) relative errors with respect to the exact analytical solution for various FOM-OpInf (a) and OpInf-OpInf (b) O-SAM couplings. Dashed horizontal lines show relative errors for FOM-FOM O-SAM couplings with the colors designated in the legend.} 
    \label{fig:clamped-rounded-square-errors-repro}
\end{figure}

\begin{figure}[ht!]
    \centering
     \begin{subfigure}{0.99\linewidth}
\hspace{6.5cm}\includegraphics[width=0.2\linewidth]{figures/clamped/legend-schwarz.png}
    \end{subfigure}
    \begin{subfigure}{0.49\linewidth}
\includegraphics[width=0.99\linewidth]{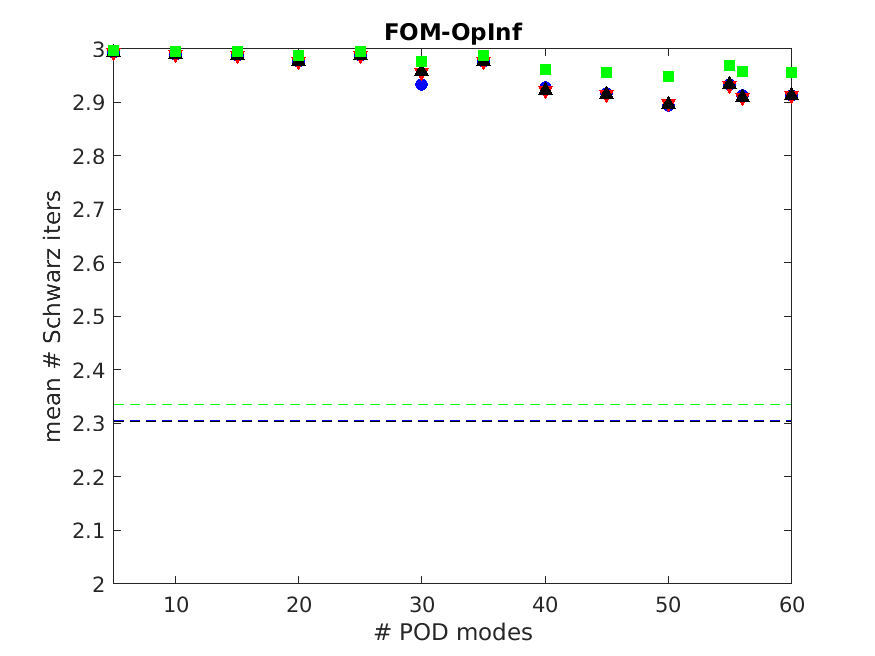}
        \subcaption{FOM-OpInf}
    \end{subfigure}
       \begin{subfigure}{0.49\linewidth}    \includegraphics[width=0.99\linewidth]{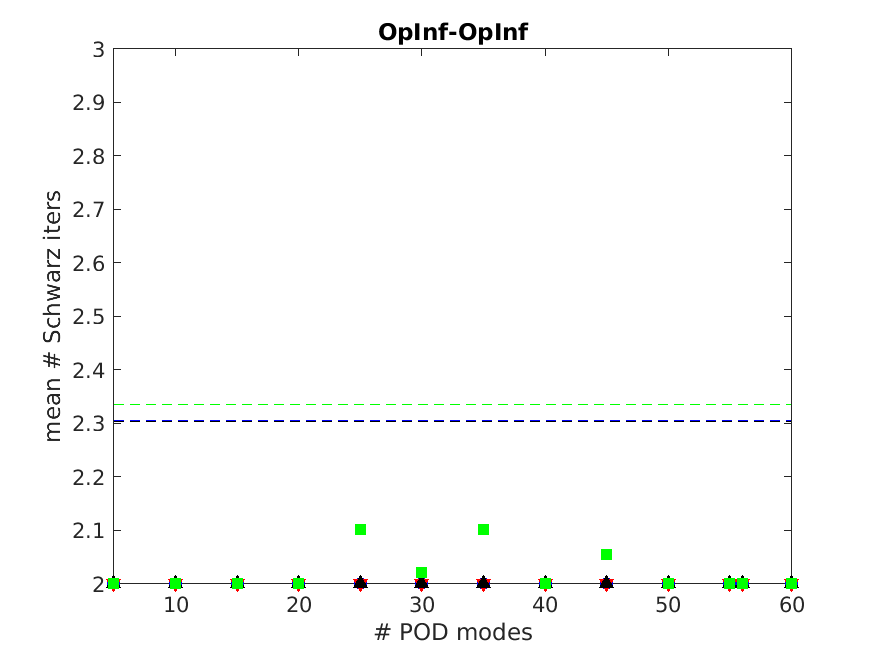}
        \subcaption{OpInf-OpInf}
    \end{subfigure}
        
      \vspace{0.5cm}
    \caption{1D linear elastic wave propagation problem, Rounded Square initial condition, reproductive regime: mean number of Schwarz iterations required to reach convergence for various FOM-OpInf (a) and OpInf-OpInf (b) O-SAM couplings. Dashed horizontal lines show the number of Schwarz iterations needed to reach convergence for FOM-FOM O-SAM couplings with the colors designated in the legend. } 
    \label{fig:clamped-rounded-square-schwarz-iters-repro}
\end{figure}

\begin{figure}[ht!]
    \centering
   %  \begin{subfigure}{0.99\linewidth}
%\hspace{6.8cm}\includegraphics[width=0.3\linewidth]{figures/clamped/legend-pareto.png}
%    \end{subfigure}
    \begin{subfigure}{0.6\linewidth}
\includegraphics[width=0.99\linewidth]{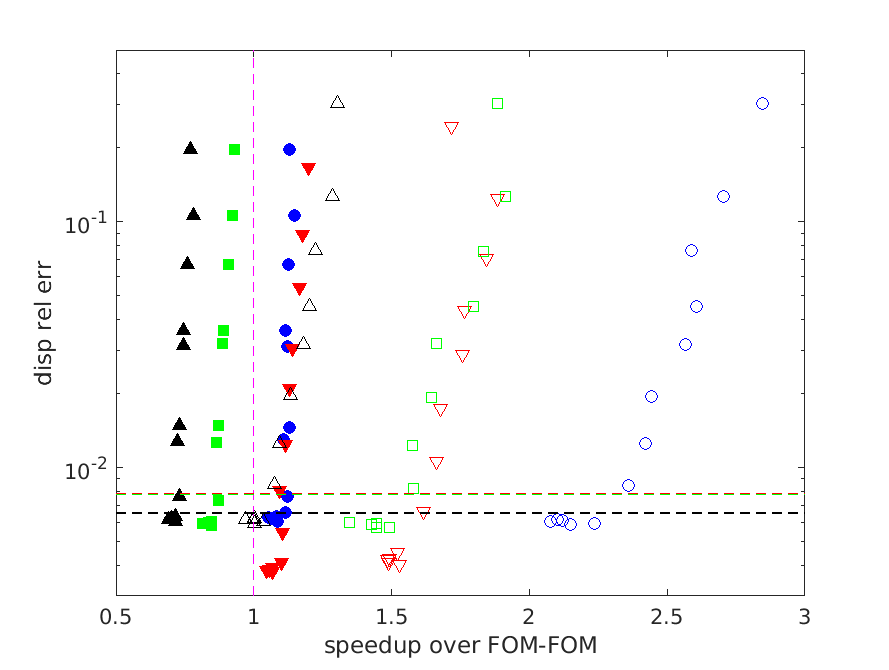}
       % \subcaption{Displacement}
        \end{subfigure}
        \raisebox{2cm}{
       \begin{subfigure}{0.3\linewidth}    \includegraphics[width=0.99\linewidth]{  figures/clamped/legend-pareto.png}
        %\subcaption{Velocity}
          \end{subfigure}}
        %  \begin{subfigure}{0.49\linewidth}    \includegraphics[width=0.99\linewidth]{figures/clamped/rounded-square-acce-pareto-repro.png}
        %\subcaption{Acceleration}
    %\end{subfigure}
        
    %  \vspace{0.5cm}
    \caption{1D linear elastic wave propagation problem, Rounded Square initial condition, reproductive regime: Pareto plot showing the speed-up over an analogous FOM-FOM coupling vs. the displacement relative errors for various FOM-OpInf (filled symbols) and OpInf-OpInf (unfilled symbols) couplings.  Dashed horizontal lines indicate relative errors for corresponding FOM-FOM O-SAM couplings. Dashed magenta vertical line indicates a speedup of 1.} %\adg{I have the same pareto plot comment here. Does it make sense to show just 1 with the legend beside it?}} 
  
    \label{fig:clamped-rounded-square-pareto-repro}
\end{figure}

\begin{figure}[ht!]
    \centering
     \begin{subfigure}{0.99\linewidth}
\hspace{6.5cm}\includegraphics[width=0.2\linewidth]{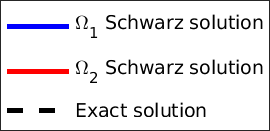}
    \end{subfigure}
    \begin{subfigure}{0.49\linewidth}
\includegraphics[width=0.99\linewidth]{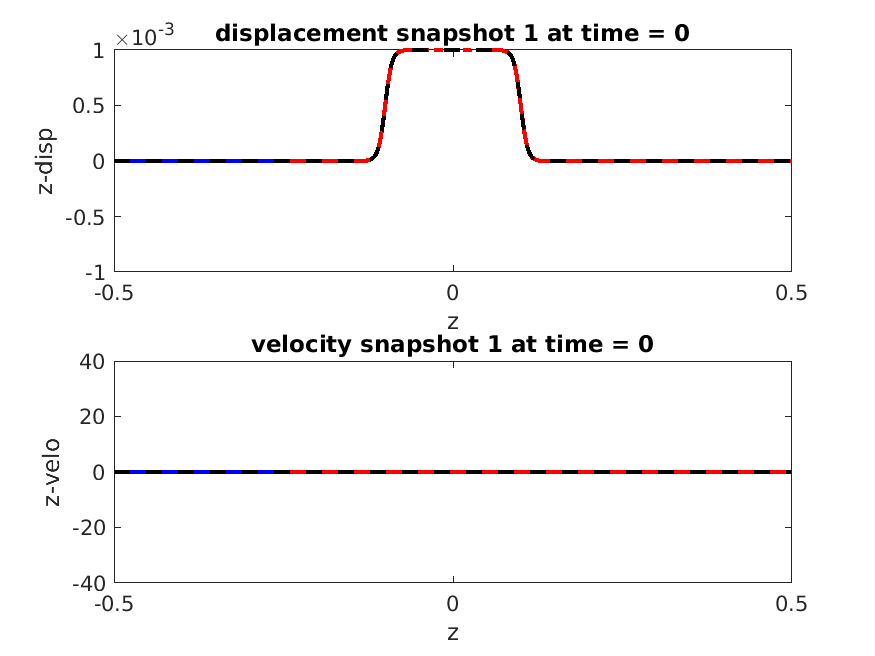}
        \subcaption{$t=0$ s}
    \end{subfigure}
       \begin{subfigure}{0.49\linewidth}    \includegraphics[width=0.99\linewidth]{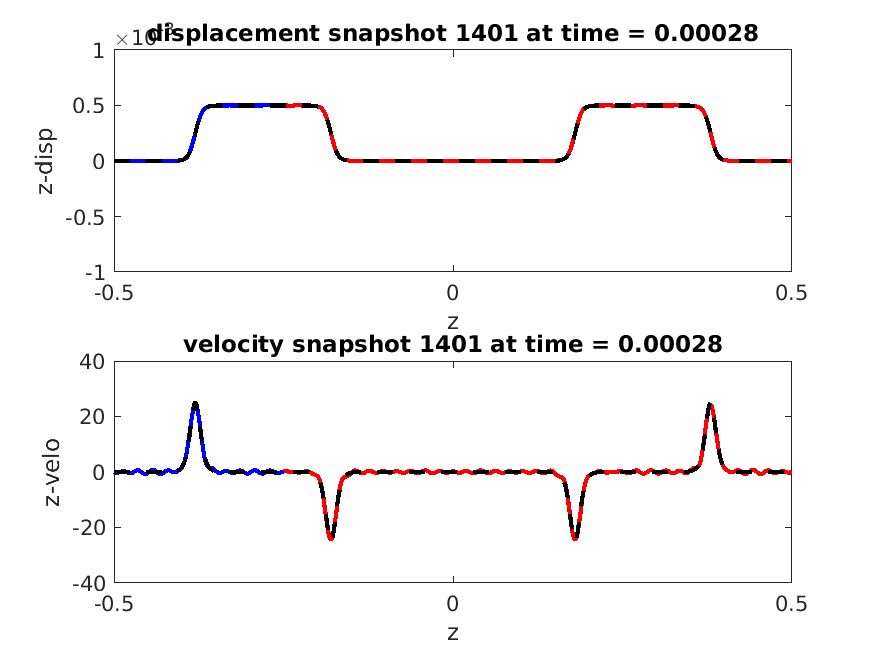}
        \subcaption{$t=2.80\times 10^{-4}$ s}
    \end{subfigure}
       \begin{subfigure}{0.49\linewidth}
\includegraphics[width=0.99\linewidth]{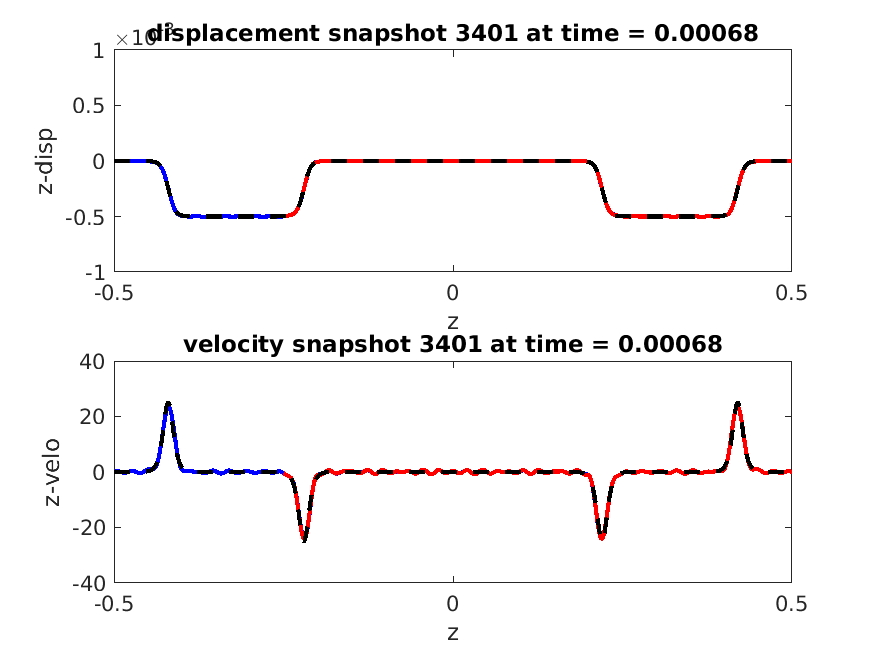}
        \subcaption{$t=6.80\times 10^{-4}$ s}
    \end{subfigure}
       \begin{subfigure}{0.49\linewidth}    \includegraphics[width=0.99\linewidth]{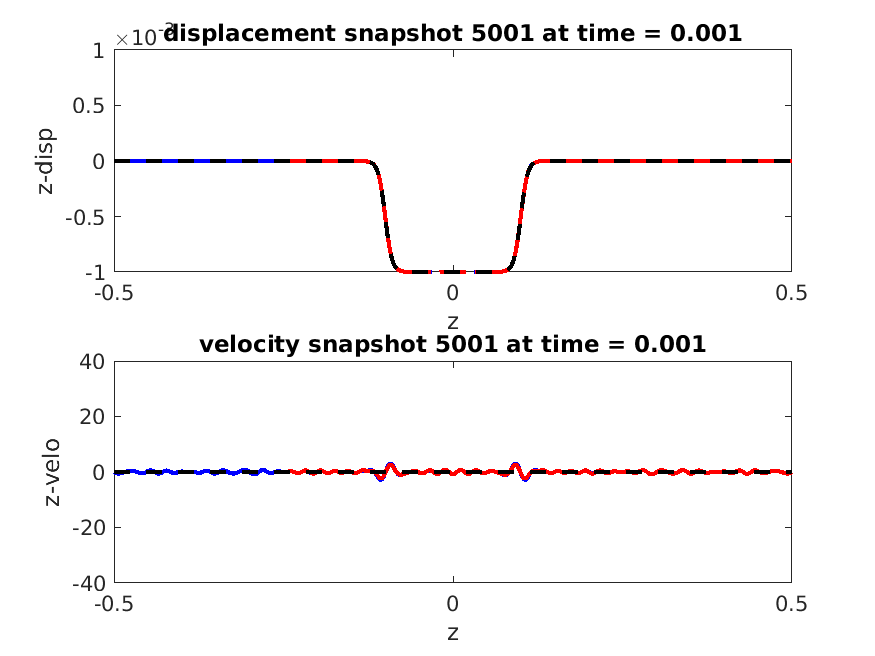}
        \subcaption{$t=1.0\times 10^{-3}$ s}
    \end{subfigure}
     %  \begin{subfigure}{0.49\linewidth}
%\includegraphics[width=0.99\linewidth]{figures/clamped/gaussian-fom-opinf-acce-repro.png}
 %       \subcaption{Acceleration, FOM-OpInf}
%       \begin{subfigure}{0.49\linewidth}    \includegraphics[width=0.99\linewidth]{figures/clamped/gaussian-opinf-opinf-acce-repro.png}
   %     \subcaption{Acceleration, OpInf-OpInf}
   % \end{subfigure}
    
      \vspace{0.5cm}
    \caption{1D linear elastic wave propagation problem, Rounded Square initial condition, reproductive regime: displacement and velocity computed using our coupled OpInf-OpInf models with $\romDimOneArg{1} = \romDimOneArg{2}=45$ and an explicit-implicit scheme with $\Delta t_1 = 1.0\times 10^{-7}$ s and $\Delta t_2 = 2.0 \times 10^{-7}$ s (blue, red), compared with the exact analytical solution (black). Solutions in $\Omega_1$ are shown in blue, whereas solutions in $\Omega_2$ are shown in red. %\crw{For these solution plots, to save space I'd suggest just showing the displacement initial conditions earlier in the section, and then here only showing two time step solutions.} \ikt{I disagree.  It's very useful to see the solution as it is crossing the Schwarz boundary, therefore I want to show it.}
    } 
    \label{fig:clamped-rounded-square-solns-repro}
\end{figure}

Although the takeaways are largely similar, we present, for completeness the reproductive Rounded Square variant of the 1D linear elastic wave propagation problem (Scenario 2), the solution to which is more difficult to represent using ROMs due to the presence of sharp gradients. We focus our discussion on the insights that are different from those discussed earlier in the context of Scenario 1. As can be seen from Figure~\ref{fig:clamped-rounded-square-errors-repro}, the Implicit-Implicit FOM-OpInf and OpInf-OpInf coupled models deliver solutions that are noticeably more accurate than their FOM-FOM analogs. This is true not only for the FOM-OpInf couplings but also the OpInf-OpInf couplings, unlike before. The reader can observe by examining Figure~\ref{fig:clamped-rounded-square-schwarz-iters-repro}(a) that the O-SAM convergence behavior is also slightly different for the FOM-OpInf models applied to Scenario 2. Whereas, for Scenario 1, the average number of Schwarz iterations required for O-SAM to converge when coupling FOM-OpInf models converged to the FOM-FOM value with basis refinement, all FOM-OpInf couplings require more Schwarz iterations to converge for Scenario 2. As expected, fewer FOM-OpInf and OpInf-OpInf coupled models are competitive from a combined accuracy and efficiency perspective for the Rounded Square variant of this problem (Figure~\ref{fig:clamped-rounded-square-pareto-repro}). This can be attributed to the problem being more difficult and requiring more modes to represent the solution, as inferred previously from Figure~\ref{fig:clamped-proj-errs}. Finally, in Figure~\ref{fig:clamped-rounded-square-solns-repro}, we plot the displacement and velocity solutions for the 45 mode mixed time-step Implicit-Explicit OpInf-OpInf coupled models at several times. While some oscillations can be seen in the velocity solutions, no coupling artifacts or propagation of these oscillations across subdomain boundaries are observed. 
% \adg{A general comment I have is that figures sometimes appear many pages away from where they are referenced. Can we fix this?}  \ikt{I looked at this and it is actually not that bad, relative to some other papers I've seen.  It is because there are a lot of figures and they are big.  I believe the `h!' syntax should put the figures as close to the text as possible, and I'm using that.  I can try to look up ways to improve this in a subsequent iteration of edits.}

\subsubsection{Scenario 3: predictive Symmetric Gaussian} \label{sec:clamped_scenario3}

\begin{figure}[ht!]
    \centering
     \begin{subfigure}{0.99\linewidth}
\hspace{6.5cm}\includegraphics[width=0.2\linewidth]{figures/clamped/legend.png}
    \end{subfigure}
    \begin{subfigure}{0.49\linewidth}
\includegraphics[width=0.99\linewidth]{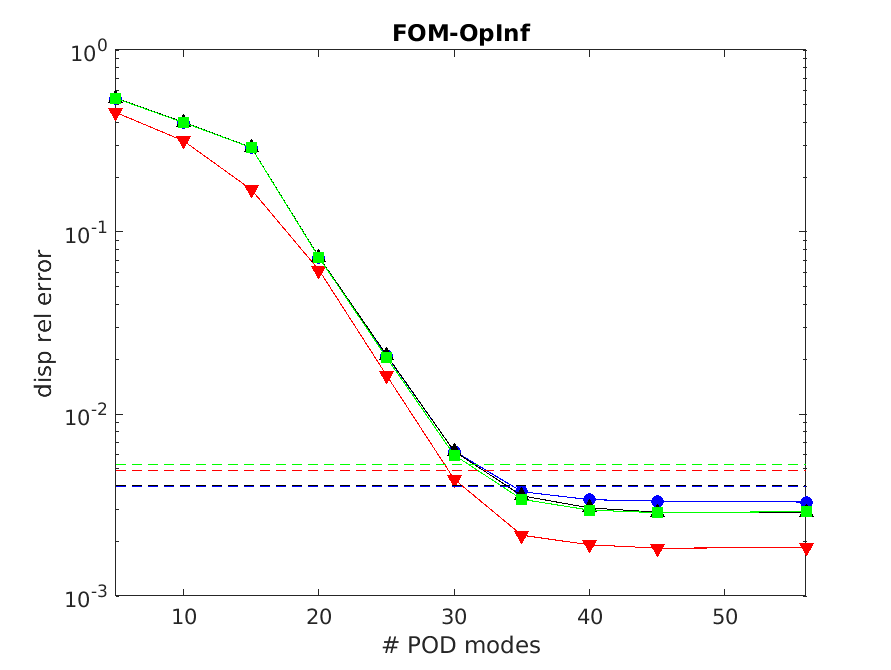}
        \subcaption{Displacement, FOM-OpInf}
    \end{subfigure}
       \begin{subfigure}{0.49\linewidth}    \includegraphics[width=0.99\linewidth]{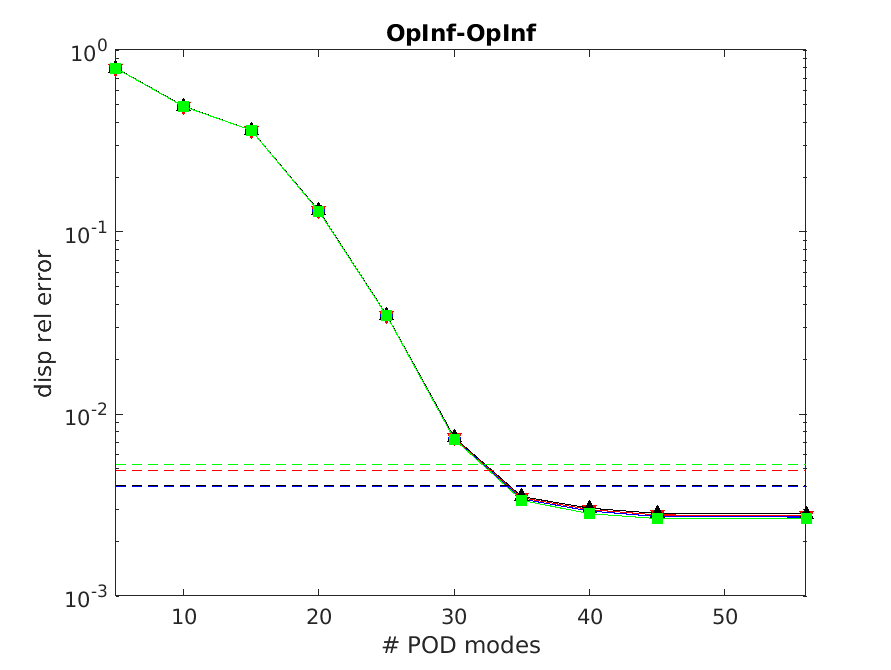}
        \subcaption{Displacement, OpInf-OpInf}
    \end{subfigure}
       \begin{subfigure}{0.49\linewidth}
\includegraphics[width=0.99\linewidth]{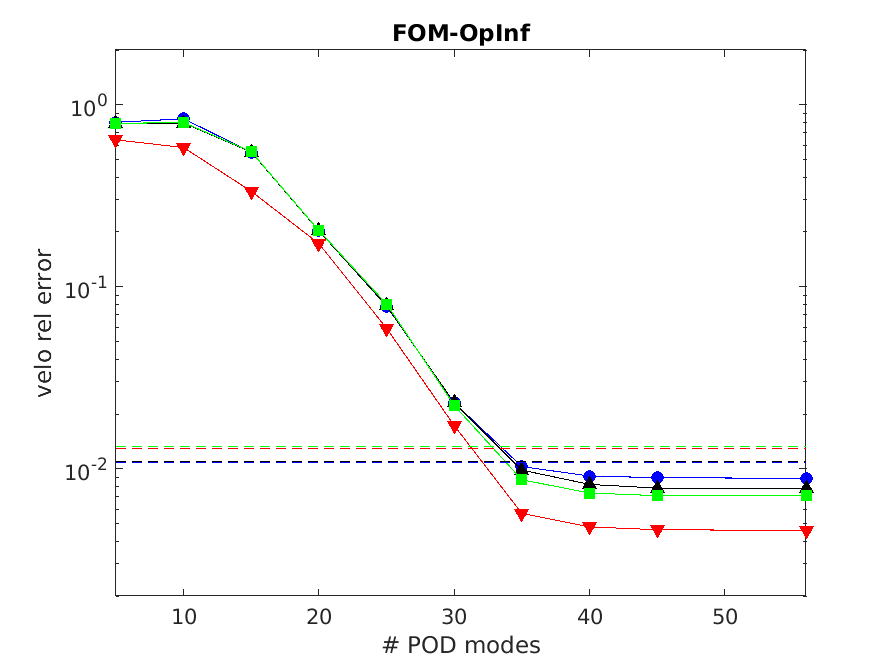}
        \subcaption{Velocity, FOM-OpInf}
    \end{subfigure}
       \begin{subfigure}{0.49\linewidth}    \includegraphics[width=0.99\linewidth]{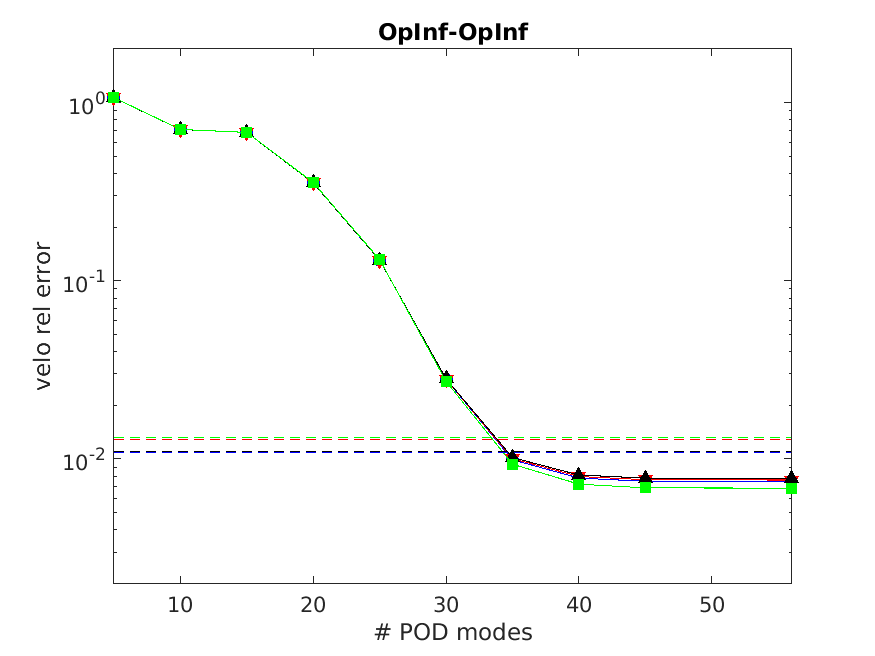}
        \subcaption{Velocity, OpInf-OpInf}
    \end{subfigure}
     %  \begin{subfigure}{0.49\linewidth}
%\includegraphics[width=0.99\linewidth]{figures/clamped/gaussian-fom-opinf-acce-predi.png}
 %       \subcaption{Acceleration, FOM-OpInf}
%    \end{subfigure}
 %      \begin{subfigure}{0.49\linewidth}    \includegraphics[width=0.99\linewidth]{figures/clamped/gaussian-opinf-opinf-acce-predi.png}
 %       \subcaption{Acceleration, OpInf-OpInf}
  %  \end{subfigure}
    
      \vspace{0.5cm}
    \caption{1D linear elastic wave propagation problem, Symmetric Gaussian initial condition, predictive regime: displacement (top row) and velocity (bottom row) relative errors with respect to the exact analytical solution for various FOM-OpInf (a) and OpInf-OpInf (b) O-SAM couplings. Dashed horizontal lines show relative errors for FOM-FOM O-SAM couplings with the colors designated in the legend.} 
    \label{fig:clamped-gaussian-disp-errors-predi}
\end{figure}

\begin{figure}[ht!]
    \centering
     \begin{subfigure}{0.99\linewidth}
\hspace{6.5cm}\includegraphics[width=0.2\linewidth]{figures/clamped/legend-schwarz.png}
    \end{subfigure}
    \begin{subfigure}{0.49\linewidth}
\includegraphics[width=0.99\linewidth]{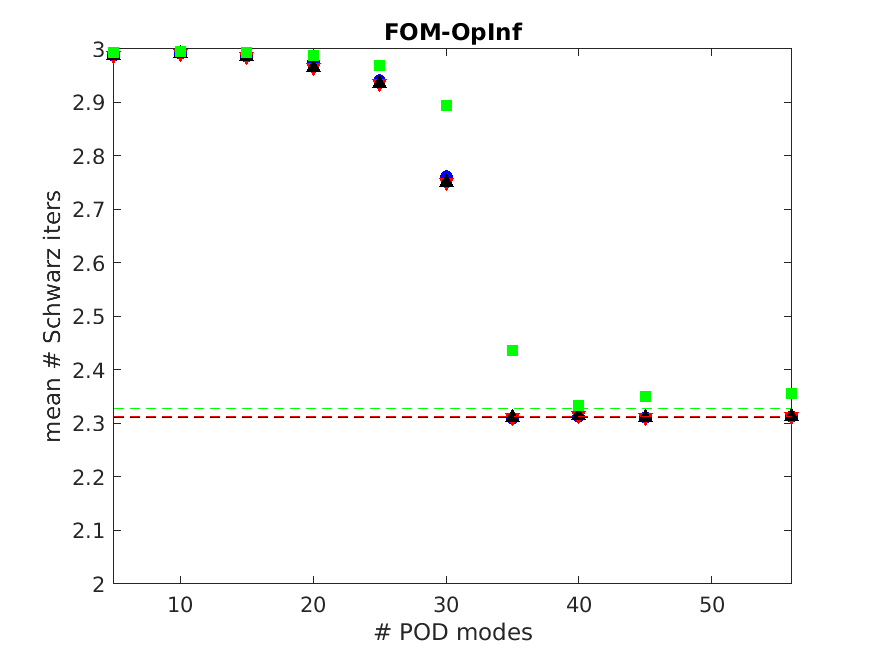}
        \subcaption{FOM-OpInf}
    \end{subfigure}
       \begin{subfigure}{0.49\linewidth}    \includegraphics[width=0.99\linewidth]{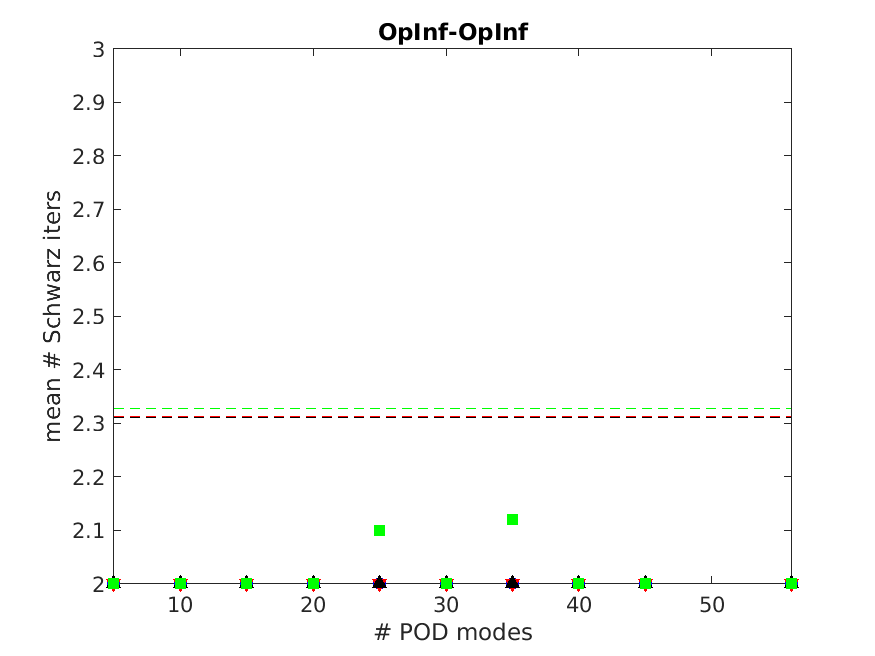}
        \subcaption{OpInf-OpInf}
    \end{subfigure}
        
      \vspace{0.5cm}
    \caption{1D linear elastic wave propagation problem, Symmetric Gaussian initial condition, predictive regime: mean number of Schwarz iterations required to reach convergence for various FOM-OpInf (a) and OpInf-OpInf (b) O-SAM couplings. Dashed horizontal lines show the number of Schwarz iterations needed to reach convergence for FOM-FOM O-SAM couplings with the colors designated in the legend.  } 
    \label{fig:clamped-gaussian-schwarz-iters-predi}
\end{figure}

\begin{figure}[ht!]
    \centering
 %   \begin{subfigure}{0.99\linewidth}
%\hspace{6.8cm}\includegraphics[width=0.3\linewidth]{figures/clamped/legend-pareto.png}
%    \end{subfigure}
    \begin{subfigure}{0.6\linewidth}
\includegraphics[width=0.99\linewidth]{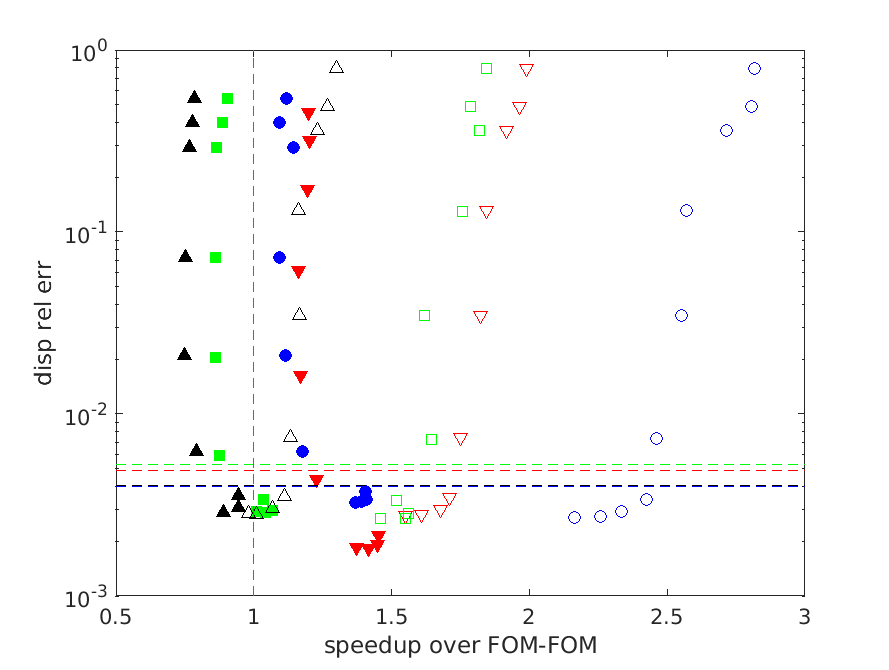}
       % \subcaption{Displacement}
        \end{subfigure}
         \raisebox{2cm}{
       \begin{subfigure}{0.3\linewidth}     \includegraphics[width=0.99\linewidth]{  figures/clamped/legend-pareto.png}
        %\subcaption{Velocity}
          \end{subfigure}}
       %   \begin{subfigure}{0.49\linewidth}    \includegraphics[width=0.99\linewidth]{figures/clamped/gaussian-acce-pareto-predi.png}
     %   \subcaption{Acceleration}
  %  \end{subfigure}
        
    %  \vspace{0.5cm}
    \caption{1D linear elastic wave propagation problem, Symmetric Gaussian initial condition, predictive regime: Pareto plot showing the speed-up over an analogous FOM-FOM coupling vs. displacement relative errors for various FOM-OpInf (filled symbols) and OpInf-OpInf (unfilled symbols) couplings. Dashed horizontal lines indicate relative errors for corresponding FOM-FOM O-SAM couplings. Dashed magenta vertical line indicates a speedup of 1.} %\adg{same comment about pareto plots here.} } 
  
    \label{fig:clamped-gaussian-pareto-predi}
\end{figure}

\begin{figure}[ht!]
    \centering
     \begin{subfigure}{0.99\linewidth}
\hspace{6.5cm}\includegraphics[width=0.2\linewidth]{figures/clamped/clamped-legend.png}
    \end{subfigure}
    \begin{subfigure}{0.49\linewidth}
\includegraphics[width=0.99\linewidth]{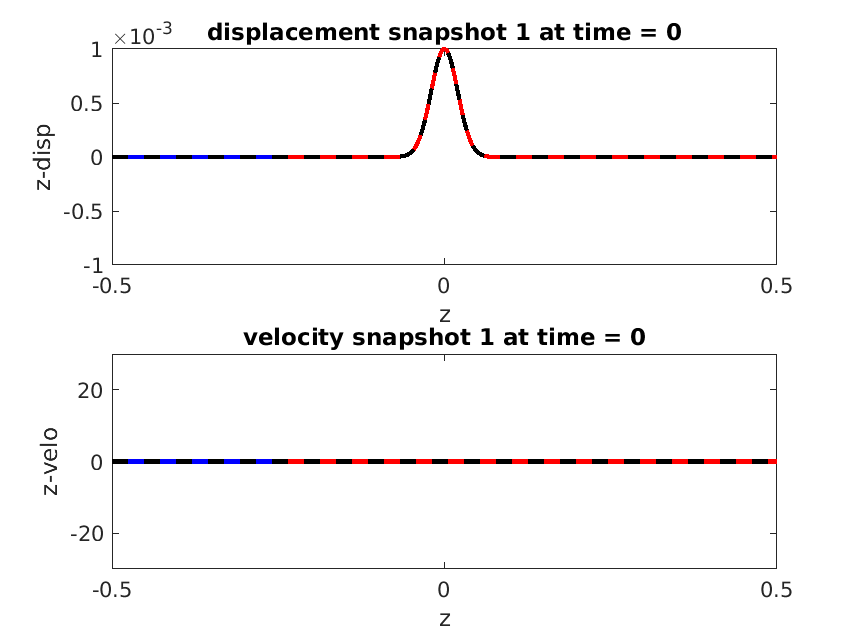}
        \subcaption{$t=0$ s}
    \end{subfigure}
       \begin{subfigure}{0.49\linewidth}    \includegraphics[width=0.99\linewidth]{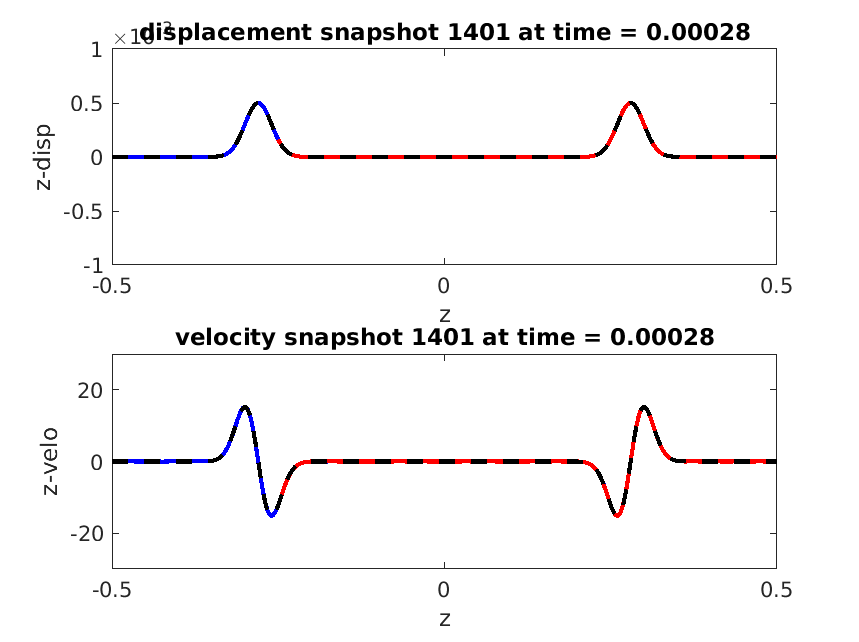}
        \subcaption{$t=2.80\times 10^{-4}$ s}
    \end{subfigure}
       \begin{subfigure}{0.49\linewidth}
\includegraphics[width=0.99\linewidth]{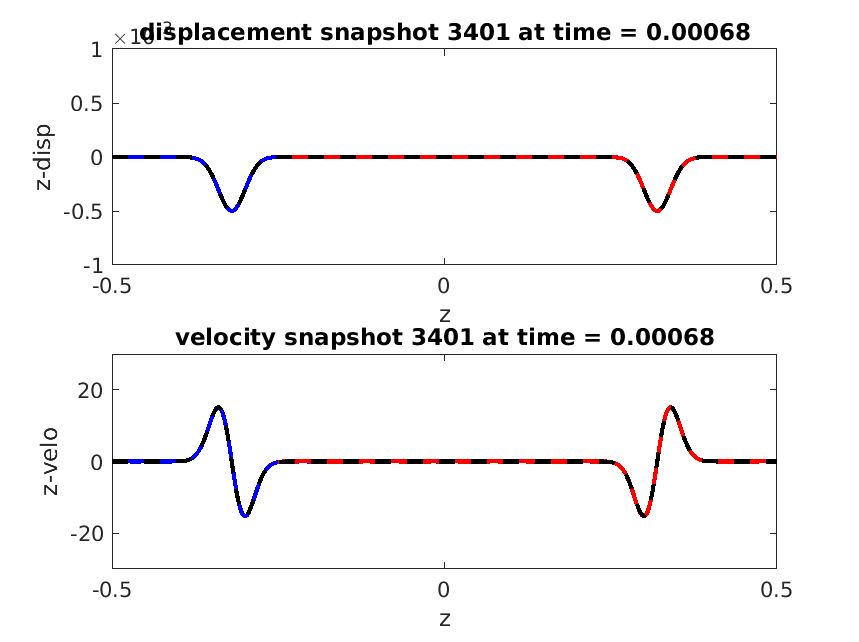}
        \subcaption{$t=6.80\times 10^{-4}$ s}
    \end{subfigure}
       \begin{subfigure}{0.49\linewidth}    \includegraphics[width=0.99\linewidth]{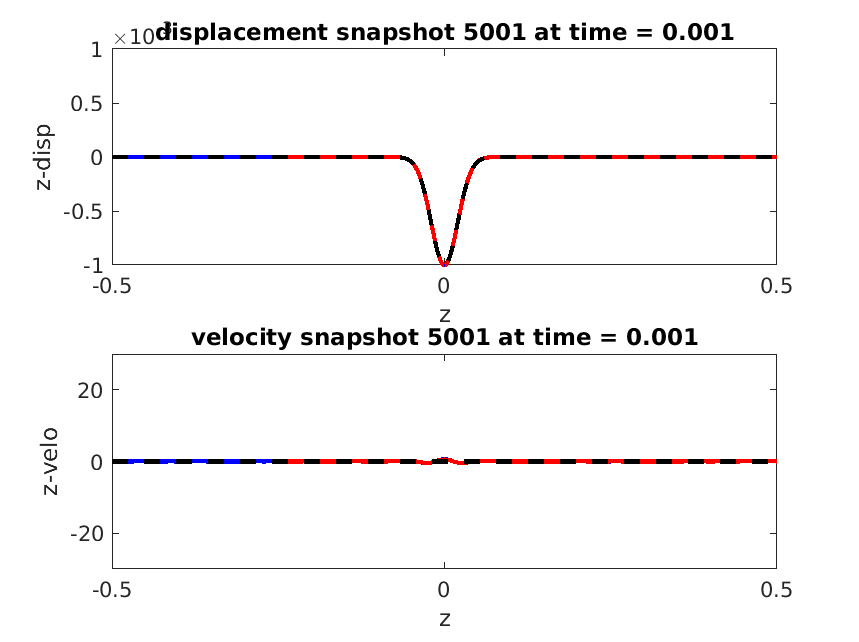}
        \subcaption{$t=1.0\times 10^{-3}$ s}
    \end{subfigure}
     %  \begin{subfigure}{0.49\linewidth}
%\includegraphics[width=0.99\linewidth]{figures/clamped/gaussian-fom-opinf-acce-repro.png}
 %       \subcaption{Acceleration, FOM-OpInf}
%       \begin{subfigure}{0.49\linewidth}    \includegraphics[width=0.99\linewidth]{figures/clamped/gaussian-opinf-opinf-acce-repro.png}
   %     \subcaption{Acceleration, OpInf-OpInf}
   % \end{subfigure}
    
      \vspace{0.5cm}
    \caption{1D linear elastic wave propagation problem, Symmetric Gaussian initial condition, predictive regime: displacement and velocity computed using our coupled OpInf-OpInf models with $\romDimOneArg{1} = \romDimOneArg{2}=35$ and an explicit-implicit scheme with $\Delta t_1 = 1.0\times 10^{-7}$ s and $\Delta t_2 = 2.0 \times 10^{-7}$ s (blue, red), compared with the exact analytical solution (black). Solutions in $\Omega_1$ are shown in blue, whereas solutions in $\Omega_2$ are shown in red.} 
    \label{fig:clamped-gaussian-solns-predi}
\end{figure}

Finally, we present results for Scenario 3, the predictive Symmetric Gaussian problem, in which the OpInf ROMs being coupled are trained on data from the Rounded Square variant of this problem and used to predict the solution with the Symmetric Gaussian initial condition. Remarkably, the key results/takeaways from Figures~\ref{fig:clamped-gaussian-disp-errors-predi}--\ref{fig:clamped-gaussian-pareto-predi} are largely the same as those for the Scenario 1 reproductive version of this problem. Figure~\ref{fig:clamped-gaussian-solns-predi} shows the displacement and velocity solutions produced by our mixed time-step Implicit-Explicit OpInf-OpInf couplings with $\romDimOneArg{1} = \romDimOneArg{2}=35$ POD modes at several times. The reader can observe that both solutions are smooth and artifact-free for all times despite this being a predictive problem.  The same result is generally \textit{not} observed for predictive monolithic ROMs. 
% \adg{We should emphasize that this is a benefit of the SAM and contrasts with the monolithic ROM case, where OoD generalization is very hard. }  \ikt{I agree, and I did this, however, I do worry a bit that it will open up a can of worms and the reader will ask why we are not comparing to a monolithic ROM.}

%% file: bolted-joint.tex
\subsection{3D nonlinear hyperelastic bolted joint problem} \label{sec:bolted-joint}

\begin{figure}[ht!]
    \centering
    \begin{subfigure}{0.40\linewidth}
        \includegraphics[width=0.85\linewidth]{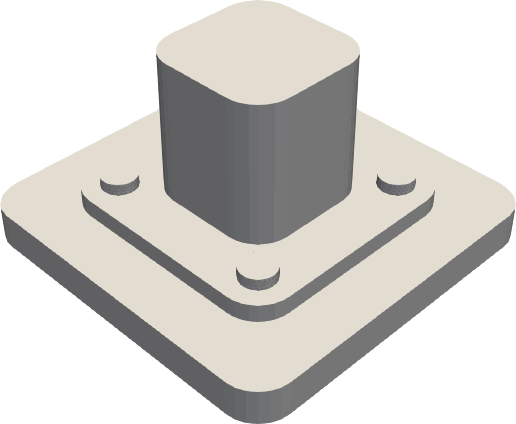}
        \vspace{0.5cm}
        \subcaption{Full geometry}
    \end{subfigure}
    \hspace{1cm}
       \begin{subfigure}{0.40\linewidth}
        \includegraphics[width=0.85\linewidth]{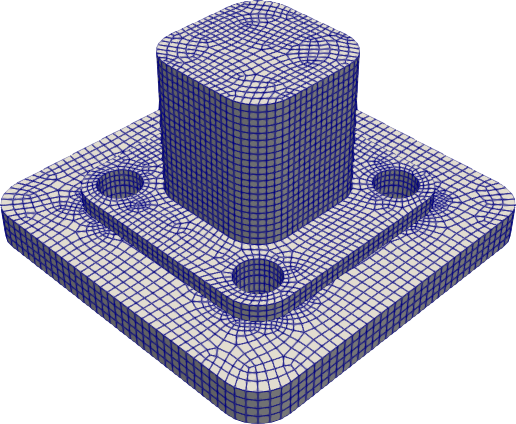}
        \vspace{0.5cm} \subcaption{$\Omega_p$ }
    \end{subfigure}
    \raisebox{-6cm}{
       \begin{subfigure}{0.40\linewidth}
        \includegraphics[width=0.85\linewidth]{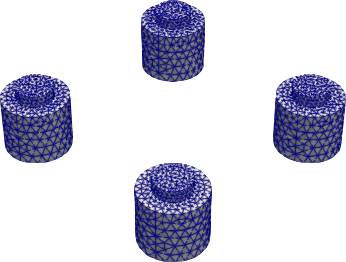}
        \vspace{0.5cm} \subcaption{$\Omega_b$ }
    \end{subfigure}
    \hspace{1cm}
       \begin{subfigure}{0.40\linewidth}
        \includegraphics[width=0.85\linewidth]{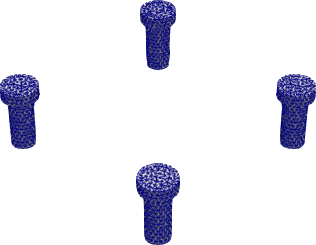}
        \vspace{0.5cm}
        \subcaption{Bolts $\subset \Omega_b$ }
    \end{subfigure}}
    \caption{3D nonlinear hyperelastic bolted joint problem: geometry and meshes.} 
    \label{fig:bolted-joint-dd}
\end{figure}

%The aim of the next three test cases, starting with the 3D nonlinear hyperelastic bolted joint problem considered in this subsection, is to study O-SAM's performance on several 3D realistic geometries.  
The aim of the second test case, the 3D bolted joint problem, is to illustrate the plug-and-play nature of the proposed coupling approach on a 3D realistic production-like geometry consisting of large-scale components joined together using small-scale fasteners (Figure \ref{fig:bolted-joint-dd}).  Such geometries are notoriously difficult to mesh due to their multiscale nature, and have the potential to benefit greatly from SAM: if one is interested in changing the shape of the geometry and/or including more detail in the bolts (e.g., by adding threading), one can create these alternate geometries/meshes offline and use SAM to seamlessly ``glue" them together without having to remesh the entire system conformally.

The specific geometry considered here consists of three pieces: a bottom plate, a top component, and four bolts which join together the plate and the component.  The square base of the joint is 127 mm $\times$ 127 mm, and the height of the joint is also 127 mm.  We prescribe a nonlinear hyperelastic Saint Venant--Kirchhoff material model \cite{holzapfel2000} with material properties corresponding to those of steel, namely a Young's modulus of $200\times 10^9$ Pa, a Poisson ratio of $\nu = 0.3$ and a density of $\rho = 7,800$ kg/m$^3$.  It can be shown that the Saint-Venant Kirchhoff material model gives rise to PDEs with cubic nonlinearities; the reader is referred to see Appendix A.2 for details.  Hence, in our numerical assessment of various FOM-ROM and ROM-ROM couplings, we will focus our attention on COpInf ROMs, which should be capable of representing the nonlinearities in the governing PDEs.

When applying SAM to the 3D bolted joint problem, we consider the natural domain decomposition of the geometry into two subdomains: one containing the four bolts, termed the ``bolts," and one containing the bottom plate and top component, termed the ``parts."  Since the focus of this paper is O-SAM, which requires a non-empty overlap region between the subdomains being coupled, the subdomain containing the bolts also includes a small circular region containing a piece of the bottom plate and top component, as shown in Figures \ref{fig:bolted-joint-dd}(c)--(d).  In typical analyses, it is critical to resolve the bolts part of the geometry, so as to accurately characterize  strain localization within the bolts and correctly predict possible fastener failure.  Toward this effect, we discretize the subdomain containing the bolts, known herein as $\Omega_b$, with a fine ten-node tetrahedral mesh containing 28,508 elements and 44,236 nodes, as shown in Figure \ref{fig:bolted-joint-dd}(c).  To demonstrate SAM's ability to couple different element types, mesh resolutions and models, we discretize the parts domain, known herein as $\Omega_p$, with a coarse eight-node tetrahedral mesh containing 22,027 elements and 29,327 nodes (Figure \ref{fig:bolted-joint-dd}(b)).  Since the accurate calculation of the solution within the bolts is of utmost importance, we are primarily interested in FOM-COpInf couplings, in which a high-fidelity model in $\Omega_b$ is coupled to an COpInf ROM in $\Omega_p$; however, for the sake of completeness, we also study herein COpInf-COpInf couplings in which an COpInf ROM is assigned to both subdomains.  The problem is initialized with a $~0$ initial condition for both the displacement and velocity fields.

\begin{remark} \label{remark:bj}
    As mentioned earlier, the bolted joint problem is an example of a test case where the monolithic governing geometry (Figure \ref{fig:bolted-joint-dd}(a)) is extremely difficult to mesh.  While we were able to generate a monolithic tetrahedral mesh of the full geometry, we were unsuccessful in producing a correct solution on this mesh using the classical Galerkin FEM.  Running this problem with four-node tetrahedral  finite elements yields an under-resolved model that produces non-physical stress patterns, whereas a ten-node tetrahedral mesh is too large and forces a very small time-step, making runs impractical.  
\end{remark}

%\begin{figure}[ht!]
%    \centering
%    \begin{subfigure}{0.49\linewidth}
        %\includegraphics[width=0.99\linewidth]{figures/bolted-joint/bolted-joint-cad.png}
        %\subcaption{Full geometry}
    %\end{subfigure}
     %  \begin{subfigure}{0.49\linewidth}
        %\includegraphics[width=0.99\linewidth]{figures/bolted-joint/bolted-joint-omega-p.png}
        %\subcaption{$\Omega_p$ }
    %\end{subfigure}
     %  \begin{subfigure}{0.49\linewidth}
        %\includegraphics[width=0.99\linewidth]{figures/bolted-joint/bolted-joint-omega-b.png}
        %\subcaption{$\Omega_b$ }
    %\end{subfigure}
     %  \begin{subfigure}{0.49\linewidth}
        %\includegraphics[width=0.99\linewidth]{figures/bolted-joint/bolted-joint-bolts.png}
        %\subcaption{Bolts $\subset \Omega_b$ }
    %\end{subfigure}
    %\caption{3D nonlinear hyperelastic bolted joint problem: geometry and meshes. \adg{We could probably put all 4 of these on one line by trimming whitespace in the figures.} \ikt{Is Figure \ref{fig:bolted-joint-dd-cropped} better?}} 
    %\label{fig:bolted-joint-dd}
%\end{figure}

\subsubsection{Reproductive variant of the bolted joint problem} \label{sec:bolted_joint_repro}

As an initial verification, we first study our mixed FOM-OpInf and OpInf-OpInf O-SAM-based couplings on a reproductive version of the bolted joint problem.  In this instance of the problem, a zero displacment boundary condition is prescribed on the bottom boundary of $\Omega_p$, while the following time-dependent boundary condition is applied on the top boundary of $\Omega_p$:
\begin{equation}
    ~u(~x,t) = \left( \begin{array}{ccc} 2-2\cos(500\pi t), & 0, & 2-2\cos(500\pi t) \end{array} \right)^{\intercal}, 
\end{equation}
for $ t \in [0, 5.8\times 10^{-5}]$ s.  %We set as an initial condition a zero initial displacement and velocity in both subdomains.  
The problem is advanced forward in time using an implicit Newmark-$\beta$ scheme with parameters $\beta = 0.25$ and $\gamma = 0.5$, and a constant time-step of $1.0 \times 10^{-6}$ s, giving rise to a total of 58 snapshots of the solution.  While our SAM-based methodology enables the coupling of disparate time-integrators with different time-steps, as shown earlier in Section \ref{sec:clamped}, we do not study such couplings in this section, as it is disadvantageous to use an explicit scheme for the bolted joint problem due to a highly-restrictive CFL condition, which gives rise to a stable time-step of $\mathcal{O}(1 \times 10^{-9})$.  Time-steps this small would not be used in production analyses involving this problem.  When applying O-SAM  in {\tt Norma.jl} \cite{Norma.jl}, we utilize a relative tolerance of $\delta_{\text{rel}} = 1.0\times 10^{-6}$ and an absolute tolerance of $\delta_{\text{abs}} = 1.0 \times 10^{-4}$ for the Schwarz convergence criteria \eqref{eq:conv_criterion}. For the Newton-based nonlinear solver and underlying iterative linear solver, we employ relative and absolute tolerances of $1.0 \times 10^{-7}$ and $1.0 \times 10^{-5}$, respectively.    To generate training data for our OpInf models, we run a coupled O-SAM simulation in which a subdomain-local FOM in $\Omega_b$ is coupled to a subdomain-local FOM in $\Omega_p$.   Running the FOM-FOM coupled model to completion takes 52m 49.5s on the {\tt Rigel} cluster, and requires an average number of 3.21 Schwarz iterations to converge, as reported in Table \ref{tab:bolted-joint-metrics-repro}.

\begin{figure}[ht!]
    \centering
    \begin{subfigure}{0.45\linewidth}
\includegraphics[width=0.99\linewidth]{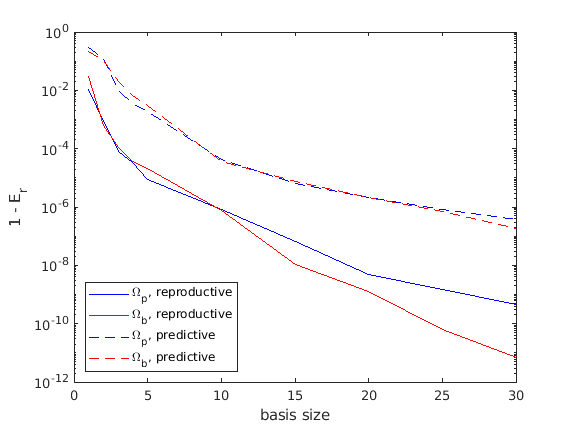}
      
    \end{subfigure}

    \caption{3D nonlinear hyperelastic bolted joint problem: POD singular value decay ($1-E_r$) in $\Omega_p$ and $\Omega_b$.  % \todo{the y-axis is wrong here I think.}
    } 
    \label{fig:bolted-joint-sv-decay}
\end{figure}

In order to build our subdomain-local COpInf models, we first apply the POD algorithm to our dataset consisting of 58 snapshots to create a reduced basis $\boldsymbol{\Phi}_r$.   The singular value decay for POD bases of varying sizes calculated for the subdomain interiors, $\Omega_p$ and $\Omega_b$, is plotted in Figure \ref{fig:bolted-joint-sv-decay}.  It is interesting to observe that the singular value decay is slightly slower in $\Omega_p$ than in $\Omega_b$ for the reproductive variant of this problem considered in the present subsection.  This is likely due to the fact that $\Omega_p$ is significantly larger than $\Omega_b$.
%\todo{Add projection errors to Figure \ref{fig:bolted-joint-sv-decay-proj-errs}, and add discussion of these figures.  It's interesting that the predictive singular value decay is pretty much the same for the bolts and the parts.  It's also interesting that the same number of modes capture less energy for the parts subdomain for larger basis sizes, as the parts subdomain is considered ``easier".  More modes is needed to capture the same energy percent for the predictive case, as expected.  Very few basis modes are needed to capture the bulk of the snapshot energy.}  
A POD basis of just 10 modes captures 99.9999\% of the snapshot energy in the interior of both subdomains as well as on the Schwarz boundary, while a basis having just 1 mode captures this same energy percentage on the remaining Dirichlet boundaries.

We begin by studying the convergence of the O-SAM algorithm when coupling a subdomain-local COpInf ROM in $\Omega_b$ with a subdomain-local FOM in $\Omega_p$ with respect to the basis size in $\Omega_b$, denoted by $\romDimOneArg{b}$.  We consider interior basis sizes that range from 5 to 30 modes for the $\Omega_b$ COpInf model, while fixing the boundary basis sizes to their 99.9999\% energy values. 
We fix the regularization parameters in the $\Omega_b$ and $\Omega_p$ subdomains, denoted by $\lambda_b$ and $\lambda_p$, respectively, to their ``optimal" values of $10^{-4}$ or $10^{-5}$, as calculated by our estimation algorithm described in Section \ref{sec:regularization}. These values are reported in Table \ref{tab:bolted-joint-metrics-repro}.  Curiously, for the FOM-COpInf couplings, increasing the basis size often leads to a decrease in the CPU time.  This can be attributed to the fact that fewer Schwarz iterations are required for convergence.

\begin{figure}[ht!]
    \centering
    \begin{subfigure}{0.45\linewidth}
        \includegraphics[width=0.99\linewidth]{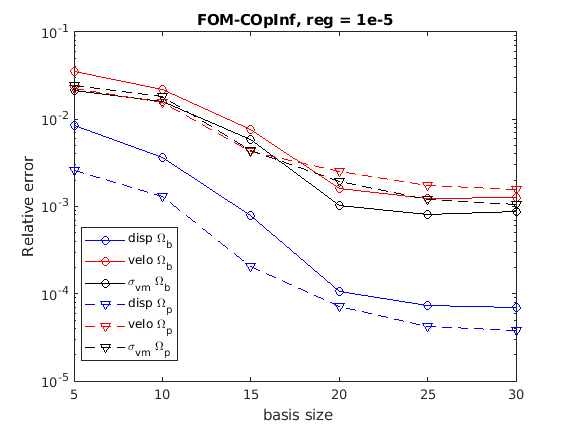}
        \subcaption{FOM-COpInf}
    \end{subfigure}
       \begin{subfigure}{0.45\linewidth}
        \includegraphics[width=0.99\linewidth]{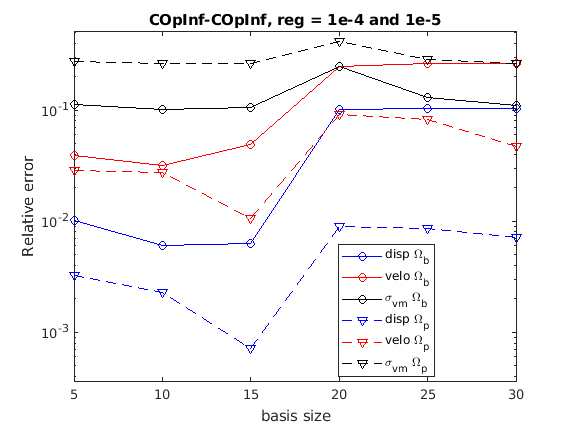}
        \subcaption{COpInf-COpInf}
    \end{subfigure}
    \caption{3D nonlinear hyperelastic bolted joint problem, reproductive regime: relative errors in displacement magnitude, velocity magnitude and average von Mises stress ($\sigma_{vm}$) for various FOM-COpInf and COpInf-COpInf couplings as a function of the reduced basis size.} 
    \label{fig:bolted-joint-repro-conv}
\end{figure}

Figure \ref{fig:bolted-joint-repro-conv}(a) plots the relative error in the displacement magnitude, the velocity magnitude and the average (over all integration points) von Mises stress ($\sigma_{vm}$) \cite{vonMises} in $\Omega_b$ and $\Omega_p$ as a function of $\romDimOneArg{b}$.  The reader can observe convergence with respect to the basis size, although this convergence plateaus for $\romDimOneArg{b} > 20$.
The accuracy and convergence of the von Mises stress follows closely that of the velocity, most likely due to the highly dynamic nature of the bolted joint problem.  For all three fields of interest, it is possible to achieve relative errors as low as $\mathcal{O}(0.1\%)$.  This is well below the typical acceptable relative error for most engineering analyses, which generally falls in the range of 1-5\%.  It is interesting to remark that the displacement solution in the ROM subdomain, $\Omega_p$, is more accurate for all values of $\romDimOneArg{p}$ than the solution in the FOM solution in $\Omega_b$.  This suggests that the solution in $\Omega_p$ is easier to represent than the solution in $\Omega_b$.

%\todo{Add explanation for why stresses have higher errors for bolted joint problem than the order of error of the displacements.  This is due to the highly dynamic nature of the problem.}

%It is interesting to remark that, while the solution in $\Omega_p$, the FOM subdomain, is in general more accurate than the solution in $\Omega_b$ when it comes to the displacement field, this is not generally the case when it comes to the velocity and acceleration fields.  Indeed, it is not possible to reduce the acceleration error below $\mathcal{O}(10\%)$.  \todo{Explain why, referring to projection error.} 

Next, we compare some performance metrics for our FOM-COpInf couplings, namely the CPU times and the mean/max numbers of Schwarz iterations required for convergence for the FOM-COpInf coupled models compared to the FOM-FOM coupled models.  These data are reported in Table \ref{tab:bolted-joint-metrics-repro}.  The reader can observe that the FOM-COpInf couplings require fewer Schwarz iterations to converge than the analogous FOM-FOM coupled model, as seen earlier in the context of the 1D linear elastic wave propagation problem (Section \ref{sec:clamped}).  
%, a result contrary to what was observed in some of our earlier work involving SAM-based coupling of intrusive projection-based ROMs \cite{Barnett:2022Schwarz}.  
While the mean number of Schwarz iterations does decrease as the basis size $\romDimOneArg{b}$ is increased, this decrease is not strictly monotonic.  For the FOM-COpInf couplings considered,   speedups of up to $1.80\times$ are observed with respect to the corresponding FOM-FOM coupling.  Curiously, employing a larger COpInf ROM in $\Omega_p$ does not necessarily lead to a larger CPU time.  This is because the total CPU time is dominated by the number of Schwarz iterations and increasing the ROM basis size can actually decrease the number of Schwarz iterations required for convergence (Table \ref{tab:bolted-joint-metrics-repro}).  Finally, it is interesting to observe that O-SAM in general converges in a very small number of iterations (between 2-5) despite the overlap region being quite small (see Figure \ref{fig:bolted-joint-dd}(c)).  This SAM convergence behavior is similar to what we observed earlier on a similar problem solved within the production {\tt SIERRA/SM} code base \cite{sierrasm}; the interested reader is referred to \cite{Mota:2022} for more details.

After assessing some FOM-COpInf couplings for the bolted joint problem, we now turn our attention to COpInf-COpInf couplings, in which a pre-trained COpInf model in $\Omega_p$ is coupled to a pre-trained COpInf model in $\Omega_b$.  We expect these couplings to give less accurate solutions, especially in the small-scale domain, $\Omega_b$, which contains fine scale features that a data-driven model may be incapable of representing.  Figure \ref{fig:bolted-joint-repro-conv}(b) shows relative errors in the displacement magnitude, velocity magnitude and the average von Mises stress in each of the two subdomains for COpInf-COpInf couplings having the same number of POD modes in each subdomain, denoted by $\romDimOneArg{} = \romDimOneArg{b} = \romDimOneArg{p}$, where $\romDimOneArg{p}$ is the number of modes retained in $\Omega_p$.  The value of $\romDimOneArg{}$ is varied between 5 and 30, as before.  For the COpInf-COpInf case, convergence with basis size is no longer observed, as the smaller ROMs for which $\romDimOneArg{} < 20$ are the most accurate (Figure \ref{fig:bolted-joint-repro-conv}(b)).  %Unfortunately, 
The von Mises stresses have a higher relative error of $\mathcal{O}(10\%)$.  
%This suggests that an analyst interested predicting very accurately the value of $\sigma_{vm}$ in the bolts should employ a FOM in $\Omega_b$, as originally conjectured.  
The reader can observe from Table \ref{tab:bolted-joint-metrics-repro} that greater CPU savings of up to $9.84\times$ are possible when utilizing COpInf-COpInf coupled models.  Like for the FOM-COpInf couplings, fewer Schwarz iterations are needed to reach convergence for the COpInf-COpInf coupled models than for the FOM-FOM coupled model.  Also, as before, employing a larger COpInf ROM does necessarily lead to a larger CPU time.

We conclude our discussion of the reproductive results by noting that it may be possible to obtain more accurate COpInf ROMs by further fine-tuning the regularization parameters used in the COpInf optimization problem.  It is well-known that, for quadratic and cubic OpInf ROMs  \cite{McQuarrie2021combustion, qian:2023nonlinearopinf}, it is best to employ different regularization parameters for different operators during the operator learning problem (see Remark \ref{remark3}).  Doing this is beyond the scope of the present work, which is focused on demonstrating SAM's ability to couple together a variety of different subdomain-local models including OpInf ROMs, rather than optimizing the individual OpInf models being coupled.

%\todo{Do we want to show monolithic results for ROMs to compare to?}

\begin{table}[ht!]
    \centering
    \caption{3D nonlinear hyperelastic bolted joint problem, reproductive regime: performance metrics for various O-SAM-based couplings.  The FOM-COpInf and COpInf-COpInf couplings requiring the lowest CPU time and the smallest number of Schwarz iterations are highlighted in green. %\adg{I like the shading!}  
    }
    \label{tab:bolted-joint-metrics-repro}
    \begin{tabular}{c|c|c|c|c|c|c}
     &  \multirow{2}{*}{$\romDimOneArg{b}$}& \multirow{2}{*}{$\romDimOneArg{p}$} & \multirow{2}{*}{$\lambda_b$} &\multirow{2}{*}{$\lambda_p$}&  \multirow{2}{*}{CPU time} & Mean/max  \\
     & & & & &  &  \# Schwarz iterations \\
     \hline  \hline 
      FOM-FOM  & $-$ & $-$  &$-$ &$-$ & 52m 49.5s& 3.21/5 \\ 
      \hline
      \multirow{6}{*}{FOM-COpInf} & $-$ & 5  & $-$& $1.0\times 10^{-5} $& 32m 58.1s& 2.62/3\\
                & $-$ & 10 &$-$& $1.0\times 10^{-5} $ &31m 24.6s & 2.48/3 \\
                & $-$ & 15 &$-$& $1.0\times 10^{-5} $ &31m 15.8s & 2.45/3 \\
                & $-$ & 20  & $-$& $1.0\times 10^{-5}$ &31m 11.3s &  2.47/3\\
                & $-$ &  25 & $-$& $1.0\times 10^{-5}$ &30m 29.3s & 2.40/3\\
                & \cellcolor{green!20}$-$ &  \cellcolor{green!20}30 & \cellcolor{green!20}$-$& \cellcolor{green!20}$1.0\times 10^{-5}$ &\cellcolor{green!20}29m 21.9s & \cellcolor{green!20}2.29/3 \\
                \hline 
        \multirow{6}{*}{COpInf-COpInf} & 5 & 5  &$1.0\times 10^{-5}$ & $1.0\times 10^{-5}$& 7m 19.0s& 2.24/3  \\
                    & \cellcolor{green!20}10 & \cellcolor{green!20}10  & \cellcolor{green!20}$1.0\times 10^{-5}$& \cellcolor{green!20}$1.0\times 10^{-5}$& \cellcolor{green!20}5m 21.9s & \cellcolor{green!20}2.22/3 \\
                    & 15 & 15  &$1.0\times 10^{-4}$& $1.0\times 10^{-4}$& 5m 40.1s&  2.05/3 \\
                    & 20 & 20  &$1.0\times 10^{-5}$& $1.0\times 10^{-5}$&6m 41.8s & 2.03/3\\
                    &25 & 25  &$1.0\times 10^{-5}$& $1.0\times 10^{-5}$& 7m 18.1s& 2.05/3\\
                    &30 & 30  &$1.0\times 10^{-5}$& $1.0\times 10^{-5}$&6m 14.7s& 2.05/3 
                    
    \end{tabular}
\end{table}

\subsubsection{Predictive variant of the bolted joint problem} \label{sec:bolted_joint_predi}

Having performed an initial assessment/verification on a reproductive variant of the 3D nonlinear hyperelastic bolted joint problem, we now evaluate our coupled hybrid models in the predictive regime.  To do so, consider the same problem formulation as above but with a more general time-dependent boundary condition applied on the top boundary of $\Omega_p$ having the form 
\begin{equation}
    ~u(~x,t) = \left( \begin{array}{ccc}a_1[1-\cos(500\pi t)], & a_2[1-\cos(500\pi t)], & a_3[1-\cos(500\pi t)]  \end{array} \right)^{\intercal},
\end{equation}
for $a_1, a_2, a_3 \in \mathbb{R}$.  In the predictive version of the bolted joint problem, we collect training data by solving the problem from time 0 to time $5.8\times 10^{-5}$ for three sets of parameters: $(a_1, a_2, a_3) = (2, 0, 0)$, $(a_1, a_2, a_3) = (0,0,2)$ and $(a_1, a_2, a_3) = (2, 2, 0)$.  This generates a total of 176 snapshots, from which we construct different size OpInf ROMs in $\Omega_b$ and $\Omega_p$.   After coupling these subdomain-local OpInf ROMs to each other and to subdomain-local FOMs, we predict the solution to the problem with $(a_1, a_2, a_3) = (2, 0, 2)$, as before. Effectively, we are generating training data by pushing the joint geometry in the positive $x$ direction, the positive $z$ direction, and at a 45 degree angle in the $x-y$ direction; we then predict the solution to the problem when the joint is pushed at a 45 degree angle in the $x-z$ direction.  We consider in our study reduced bases consisting of between 10 and 24 POD modes.  For the predictive bolted joint problem, $\romDimOneArg{b} = \romDimOneArg{p} = 15$ modes are needed to capture 99.999\% of the snapshot  energy in $\Omega_b$ and $\Omega_p$, whereas $\romDimOneArg{b} = \romDimOneArg{p} = 24$ modes are needed to capture 99.9999\% of the snapshot energy in $\Omega_b$ and $\Omega_p$.  On the Schwarz boundary of $\Omega_b$, 23 POD modes capture 99.9999\% of the snapshot energy, compared to 24 POD modes for the Schwarz boundary of $\Omega_p$.  As for the reproductive version of this problem, a single POD mode captures 99.9999\% of the snapshot energy for the Dirichlet boundaries.  

\begin{figure}[ht!]
    \centering
    \begin{subfigure}{0.45\linewidth}
\includegraphics[width=0.99\linewidth]{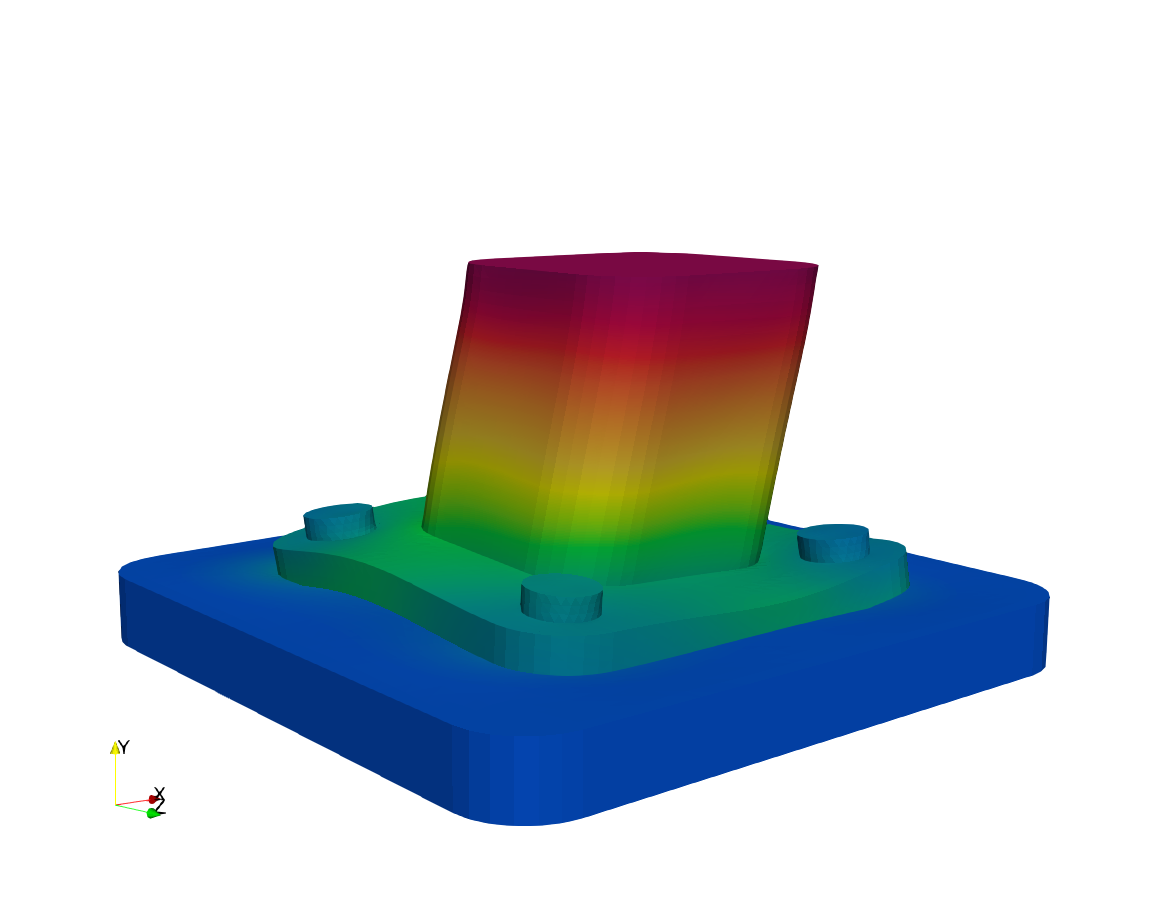}
        \subcaption{FOM-FOM}
    \end{subfigure}
       \begin{subfigure}{0.45\linewidth}    \hspace{0.5cm}\includegraphics[width=0.99\linewidth]{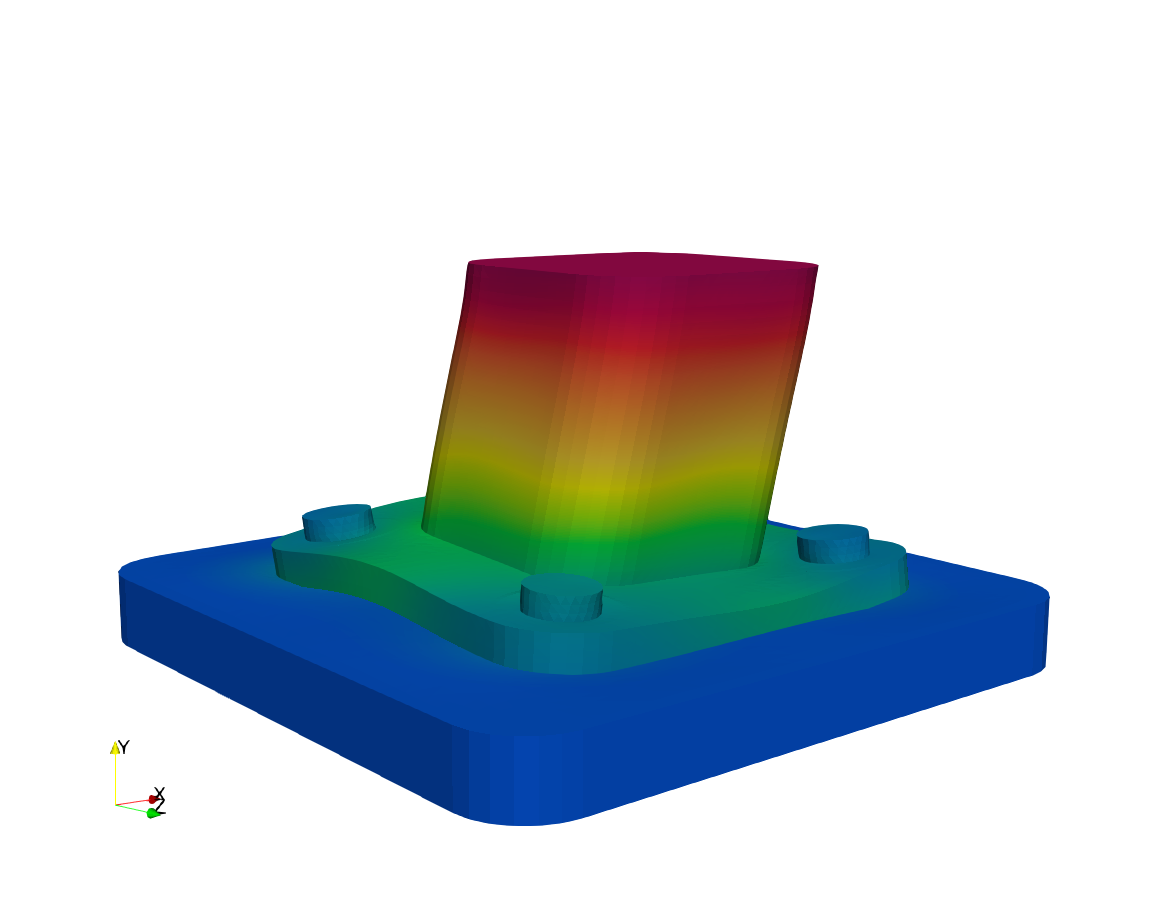}
        \subcaption{FOM-COpinf ($\romDimOneArg{b}=24$)}
    \end{subfigure}
       \begin{subfigure}{0.45\linewidth}
\includegraphics[width=0.99\linewidth]{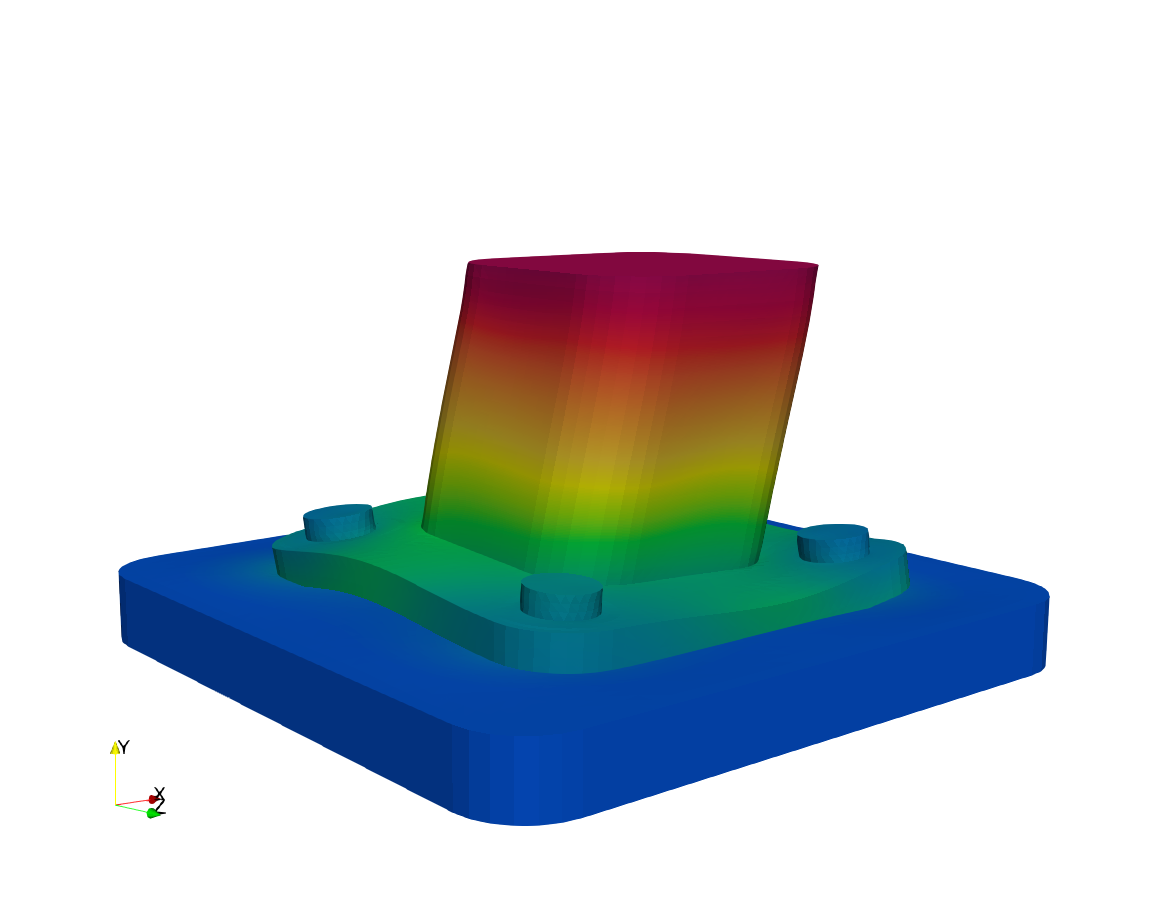}
        \subcaption{COpInf-COpInf ($\romDimOneArg{} = \romDimOneArg{b} = \romDimOneArg{p} = 10$)}
    \end{subfigure}
    \begin{subfigure}{0.45\linewidth}
\hspace{0.5cm}\includegraphics[width=0.99\linewidth]{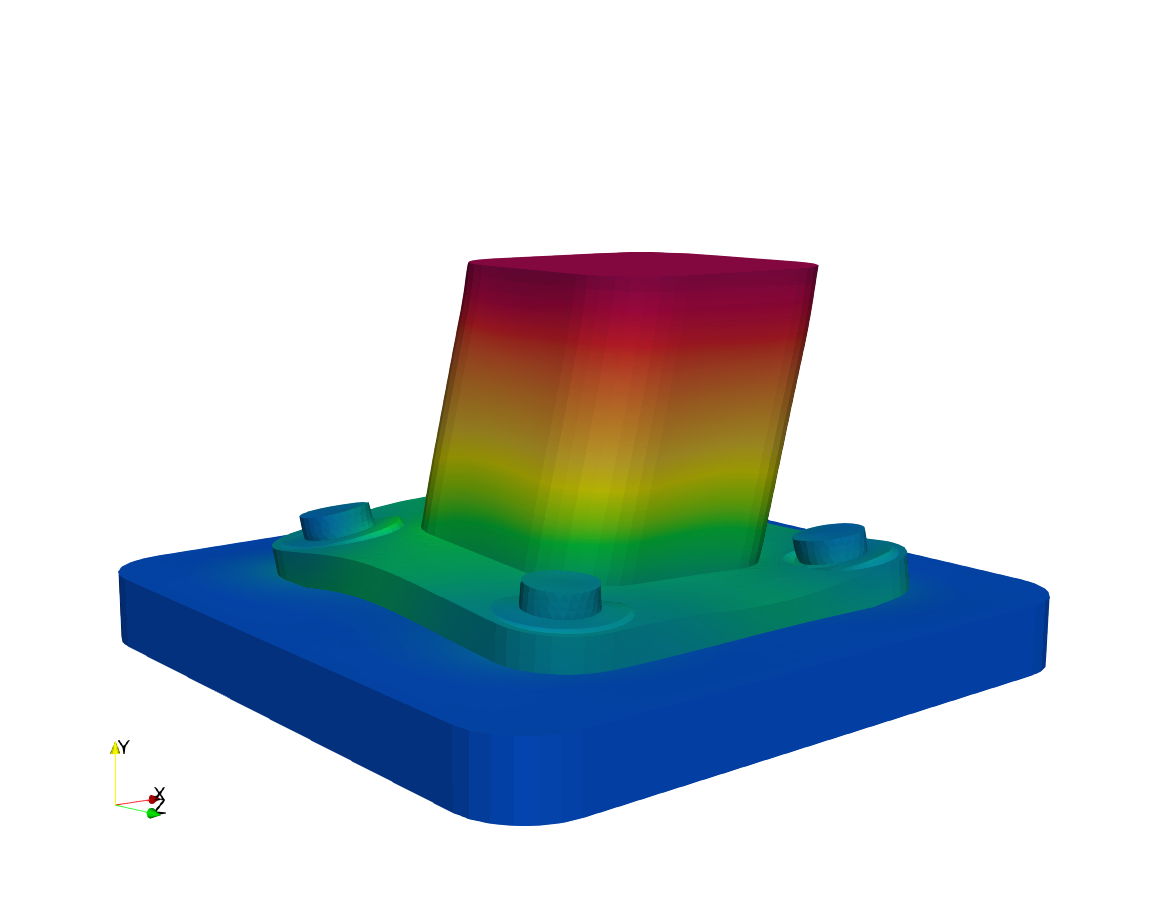}
        \subcaption{COpInf-COpInf ($\romDimOneArg{} = \romDimOneArg{b} = \romDimOneArg{p} = 20$)}
    \end{subfigure}
      \begin{subfigure}{0.99\linewidth}
      \vspace{0.5cm}
\hspace{4.5cm}\includegraphics[width=0.5\linewidth]{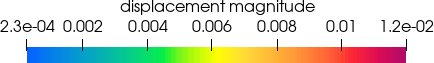}
    \end{subfigure}
     
      \vspace{0.5cm}
    \caption{3D nonlinear hyperelastic bolted joint problem, predictive regime: plots of the displacement magnitude for a FOM-FOM O-SAM-based coupling (a) compared to a COpInf-FOM Schwarz coupling with $\romDimOneArg{b} = 24$ (b), a COpInf-COpInf Schwarz coupling with $\romDimOneArg{}=\romDimOneArg{b}=\romDimOneArg{p}=10$ (c), and a COpInf-COpInf Schwarz coupling with $\romDimOneArg{}=\romDimOneArg{b}=\romDimOneArg{p}=20$ (d).}
    %\adg{I also think this figure could be on one line.  But, it's really cool and we should make sure the differences are still visible.} } 
    \label{fig:bolted-joint-disp-solns}
\end{figure}

\begin{figure}[ht!]
    \centering
    \begin{subfigure}{0.35\linewidth}
\includegraphics[width=0.99\linewidth]{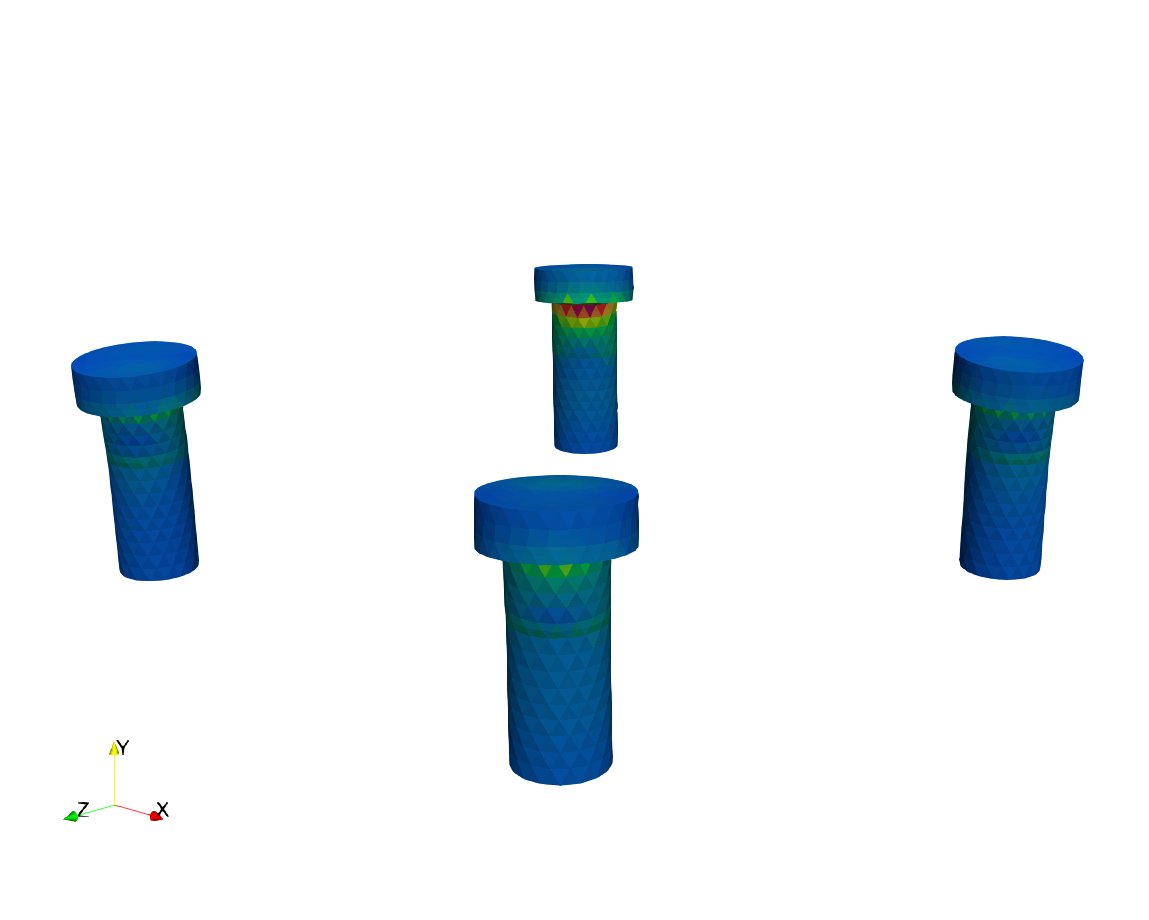}
        \subcaption{FOM-FOM}
    \end{subfigure}
      \hspace{1cm} \begin{subfigure}{0.35\linewidth}    \includegraphics[width=0.99\linewidth]{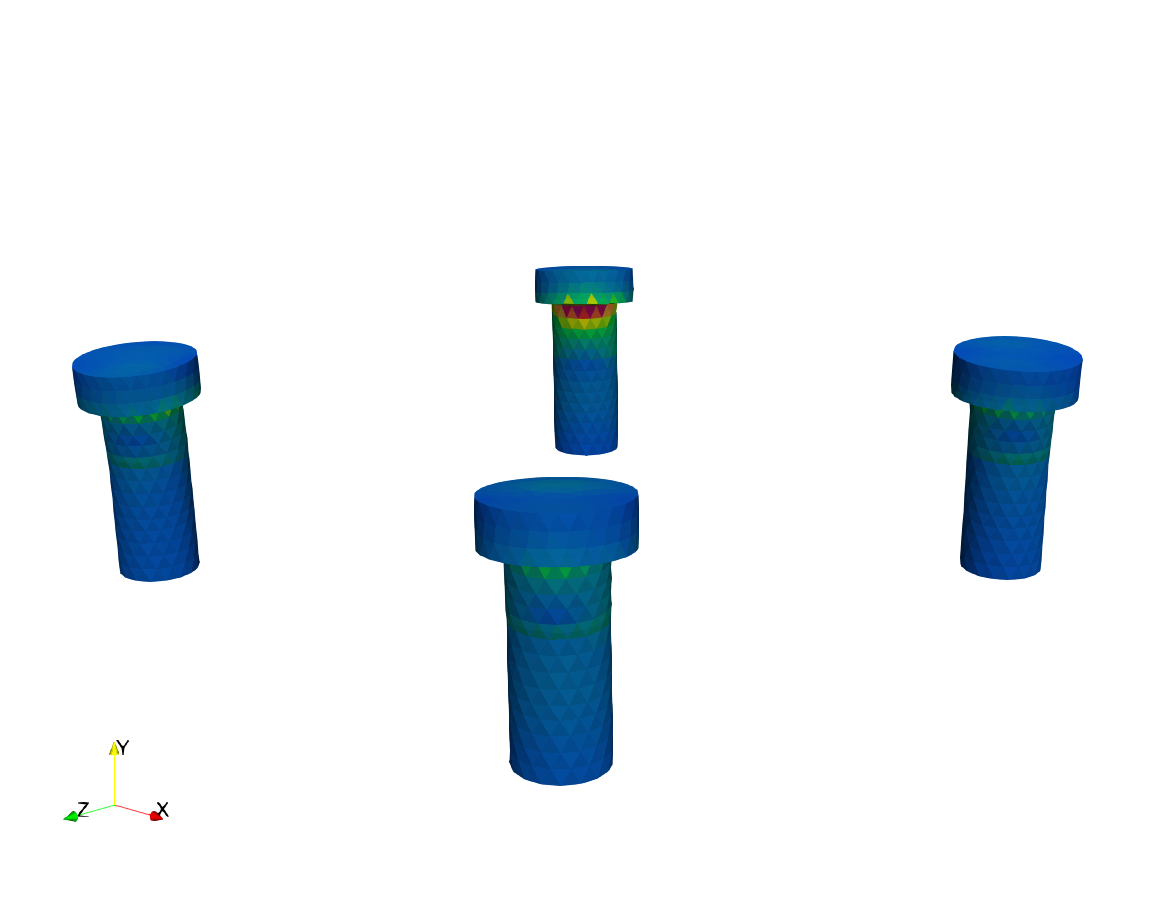}
        \subcaption{COpInf-FOM ($\romDimOneArg{b}=24$)}
    \end{subfigure}
       \begin{subfigure}{0.35\linewidth}
\includegraphics[width=0.99\linewidth]{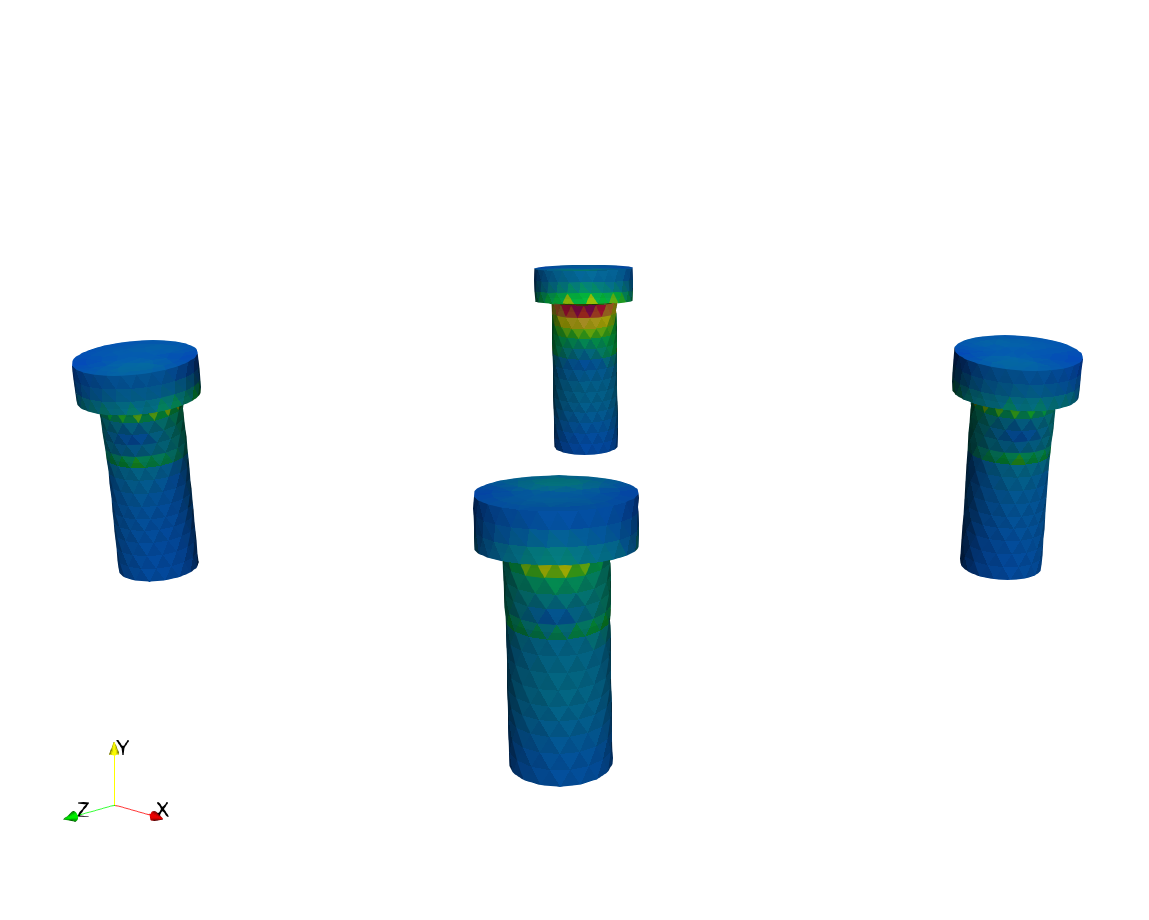}
        \subcaption{COpInf-COpInf ($\romDimOneArg{} = \romDimOneArg{b} = \romDimOneArg{p} = 10$)}
    \end{subfigure}
    \hspace{1cm} \begin{subfigure}{0.35\linewidth}
\includegraphics[width=0.99\linewidth]{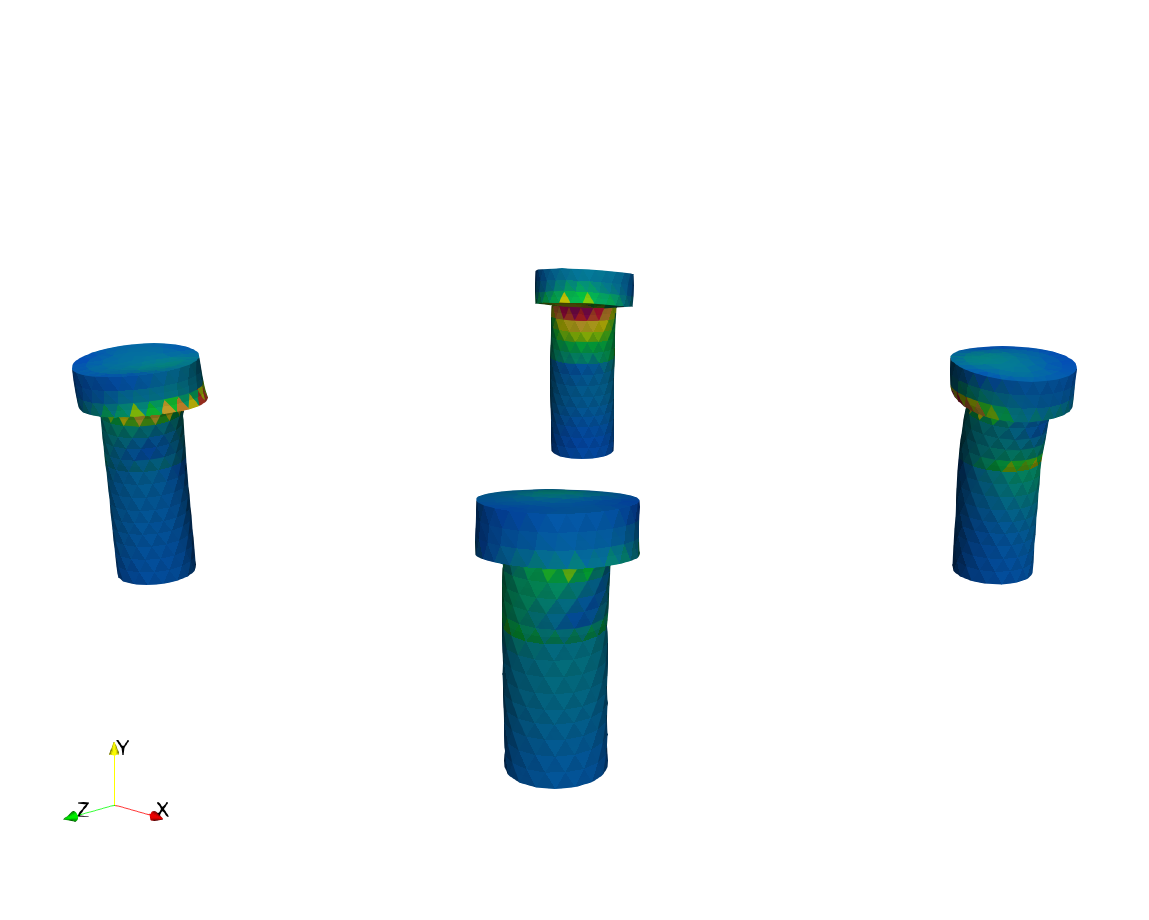}
        \subcaption{COpInf-COpInf ($\romDimOneArg{} = \romDimOneArg{b} = \romDimOneArg{p} = 20$)}
    \end{subfigure}
      \begin{subfigure}{0.99\linewidth}
      \vspace{0.5cm}
\hspace{4.5cm}\includegraphics[width=0.5\linewidth]{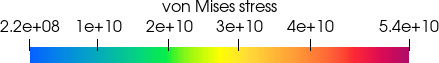}
    \end{subfigure}
     
      \vspace{0.5cm}
    \caption{3D nonlinear hyperelastic bolted joint problem, predictive regime: plots of the average von Mises stresses $\sigma_{vm}$ for a FOM-FOM O-SAM-based coupling (a) compared to a COpInf-FOM Schwarz coupling with $\romDimOneArg{b} = 24$ (b), a COpInf-COpInf Schwarz coupling with $\romDimOneArg{}=\romDimOneArg{b}=\romDimOneArg{p}=10$ (c), and a COpInf-COpInf Schwarz coupling with $\romDimOneArg{}=\romDimOneArg{b}=\romDimOneArg{p}=20$ (d). %\adg{This figure is awesome!} 
    } 
    \label{fig:bolted-joint-stress-vm-solns}
\end{figure}

%As for the reproductive variant of this problem, 
We report the relative errors in the displacement magnitude, velocity magnitude and von Mises stress solutions in each of the subdomains as a function of $\romDimOneArg{b}$ and $\romDimOneArg{p}$; see Tables \ref{tab:bolted-joint-metrics-predi-fom-rom} and \ref{tab:bolted-joint-metrics-predi-rom-rom} for FOM-COpInf and COpInf-COpInf results, respectively.  The ``optimal" values of the regularization parameters in $\Omega_b$ and $\Omega_p$, identified by our estimation algorithm (Section \ref{sec:regularization}) were again  $10^{-4}$ or $10^{-5}$.  For the FOM-COpInf couplings (Table \ref{tab:bolted-joint-metrics-predi-fom-rom}), relative errors of $\mathcal{O}(0.1)$--$\mathcal{O}(1\%)$ and $\mathcal{O}(1\%)$ can be achieved for the displacement and velocity fields, respectively, and relative errors of 7--8\% are possible for the von Mises stress $\sigma_{vm}$.  The accuracy of the COpInf-COpInf coupled models is worse and degrades with basis refinement, though errors of $4.41\%$, $10.7\%$ and $12.4\%$ are achievable for the displacement, velocity and $\sigma_{vm}$ fields in the bolts subdomain $\Omega_b$ %are achievable 
with COpInf ROMs based on 10 POD modes.  Figures \ref{fig:bolted-joint-disp-solns} and \ref{fig:bolted-joint-stress-vm-solns} show the displacement magnitude in the full domain, and $\sigma_{vm}$ in the bolts, respectively, for various O-SAM-based couplings.  For the COpInf-COpInf couplings, we show the best and worst cases in subfigures (c) and (d), respectively.  These correspond go $\romDimOneArg{}=\romDimOneArg{b}=\romDimOneArg{p} = 10$ and $\romDimOneArg{}=\romDimOneArg{b}=\romDimOneArg{p}=20$, respectively.  It can be seen that the large-scale deformation of the bolted joint is correct for all the models being assessed (Figure \ref{fig:bolted-joint-disp-solns}).  Figure \ref{fig:bolted-joint-stress-vm-solns} shows that the maximum stress occurs in the rear fastener in the negative $x-z$ coordinate plane, parallel to the direction in which the forcing is applied, as expected.  The  COpInf-FOM $\sigma_{vm}$ in the bolts (Figure \ref{fig:bolted-joint-stress-vm-solns}(b)) matches the FOM-FOM $\sigma_{vm}$ (Figure \ref{fig:bolted-joint-stress-vm-solns}(a)) remarkably well.  For the best COpInf-COpInf model with $\romDimOneArg{} = \romDimOneArg{b} = \romDimOneArg{p} = 10$ (Figure \ref{fig:bolted-joint-stress-vm-solns}(c)), the reader can observe that $\sigma_{vm}$ is slightly overestimated in all of the bolts.  The overestimation of $\sigma_{vm}$ is more profound for the worst COpInf-COpInf models, having $\romDimOneArg{}=\romDimOneArg{b}=\romDimOneArg{p}=20$ (Figure \ref{fig:bolted-joint-stress-vm-solns}(d)).  While the error in $\sigma_{vm}$ in the $\Omega_b$ subdomain for this coupling is close to 31\%, which suggests that the magnitude of $\sigma_{vm}$ is substantially off, it is important to remark  that the \textit{locations} of maximum stress are still correctly predicted by the model, as can be seen by examining Figure \ref{fig:bolted-joint-stress-vm-solns}(d). 
Hence, even the less accurate COpInf-COpInf coupled models can be useful, as they can provide a conservative estimate of the bolt failure.
%\adg{It is worth mentioning here or in the conclusion to the paper that the OpInf models can be useful even when they are wrong.  For example, here they give a conservative estimate of the bolt failure.} \ikt{Yes!  Done.}

%Additionally, some incorrect deformation around the bolts is observed in the displacement magnitude plot of the solution for this case (Figure \ref{fig:bolted-joint-disp-solns}(d)).  

%As before, we note that it may be possible to improve the accuracy of our hybrid models featuring COpInf ROMs by performing a careful tuning of the individual COpInf ROMs, e.g., by selecting alternate regularization parameters in the COpInf minimization problem; however, fine-tuning these COpInf models is beyond the scope of the present work, where the goal is to demonstrate that our O-SAM-based coupling of various fidelity models does not introduce nonphysical artifacts. \adg{I feel like this paragraph is not really needed-- we say this several times in various places.}

It is interesting to remark from Tables \ref{tab:bolted-joint-metrics-predi-fom-rom} and \ref{tab:bolted-joint-metrics-predi-rom-rom} that increasing the basis size does not necessarily reduce the error, which may seem counterintuitive.  We suspect that the main reason for this is ill-conditioning of the data matrices corresponding to the underlying least squares optimization problem \eqref{eq:quadratic_opinf}, together with the fact that this section considers a \textit{predictive} version of the bolted joint problem, meaning convergence with basis refinement is not guaranteed.  Since the OpInf ROM for the bolted joint problem is cubic, the relevant data matrix corresponding to subdomain $\Omega_i$ takes the form:
\begin{equation} \label{eq:data_matrix_cubic}
    ~D_i^c:= \left( \begin{array}{c} 
    \bar{~A}_i \\
    \bar{~U}_i \\
    \bar{~U}_i \otimes \bar{~U}_i \\
    \bar{~U}_i \otimes \bar{~U}_i \otimes \bar{~U}_i 
    \end{array}\right),
\end{equation}
where $\bar{~U}_i := ~\Phi_i ~U_i$ and $\bar{~A}_i := ~\Phi_i \ddot{~U}_i$, with $~U_i$ denoting the matrix of displacement snapshots for subdomain $\Omega_i$ with $i \in \{ b, p\}$.  Table \ref{tab:data_matrix_rank} reports the size, rank and condition number of the $~D_i^c$ matrix for each COpInf ROM considered.  The reader can observe that the data matrices are extremely ill-conditioned and rank deficient, which confirms the need for regularization in \eqref{eq:quadratic_opinf}.  While regularization is crucial for preventing overfitting
and ensuring long-term stability, it inevitably forces the learned system operators to deviate
from the theoretically exact, un-regularized values, which can have an impact on convergence
with respect to the basis size.

Also reported in Tables \ref{tab:bolted-joint-metrics-predi-fom-rom} and \ref{tab:bolted-joint-metrics-predi-rom-rom} are CPU times for the various couplings.  A modest speed-up of $1.65\times$ is achieved for the FOM-COpInf couplings, whereas a much greater speedup of $6.12\times$ is observed for the most accurate COpInf-COpInf coupling.  As for the reproductive version of this problem, we observe that O-SAM converges in fewer Schwarz iterations for the FOM-COpInf and COpInf-COpInf couplings than for the FOM-FOM coupling, which contributes to the observed CPU-time reduction.  It is likely that bigger cost reductions are possible when employing the cheaper QOpInf models, but a study involving QOpInf-FOM and QOpInf-QOpInf couplings for the bolted joint problem is beyond the scope of the present work.

%\todo{Say that, in Figure \ref{fig:tension-specimen-stress-vm-solns}, even the ``bad" ROM results would lead to correct prediction of failure.}

%\todo{Add the following to explain why we are not including monolithic results (from Alejandro): For the paper, I suggest we frame the monolithic case along these lines:
%\begin{enumerate}
%    \item It is very difficult to mesh robustly.
%\item A composite TET10 model is too large and forces a very small time step, making runs impractical.
%\item 	The TET4 model is under-resolved and produces non-physical stress patterns, so we do not use it for comparison.
%\end{enumerate}
%}

\begin{table}[ht!]
    \centering
    \caption{3D nonlinear hyperelastic bolted joint problem, predictive regime: relative errors and performance metrics for various FOM-COpInf O-SAM-based couplings.  The best coupling in terms of a combination of the overall accuracy and efficiency is highlighted in green.}
    \label{tab:bolted-joint-metrics-predi-fom-rom}
    \begin{tabular}{c|c|c|c|c|c|c}
    & Field & FOM-FOM & \multicolumn{4}{c}{FOM-COpInf} \\ 
    \hline
    $\romDimOneArg{b}$ & $-$& $-$ & 10 &  15 & 20 &\cellcolor{green!20} 24 \\
    \hline 
    $\lambda_b$ & $-$ & $-$ & $1.0\times 10^{-4}$ & $1.0\times 10^{-4}$ & $1.0\times 10^{-4}$ & \cellcolor{green!20}$1.0\times 10^{-4}$\\ 
    \hline 
    
    \multirow{3}{*}{$\Omega_b$ rel errors}& $~u$ & $-$ &$4.51\times 10^{-2}$&$3.94  \times 10^{-2}$ & $3.75 \times 10^{-2}$ & \cellcolor{green!20}$3.51\times 10^{-2}$\\ 
    & $\dot{~u}$ & $-$ &$9.40 \times 10^{-2}$   & $5.61 \times 10^{-2}$ & $4.79 \times 10^{-2}$  & \cellcolor{green!20}$4.25 \times 10^{-2}$\\ 
   & $\sigma_{vm}$& $-$ & $9.29 \times 10^{-2}$  & $8.40 \times 10^{-2}$& $8.14\times 10^{-2}$ & \cellcolor{green!20}$7.74 \times 10^{-2}$\\ 
    \hline
     \multirow{3}{*}{$\Omega_p$ rel errors}  & $~u$& $-$ & $ 9.71\times 10^{-3}$  &$7.26 \times 10^{-3}$ & $7.49\times 10^{-3}$ & \cellcolor{green!20}$6.92\times 10^{-3}$\\
      & $\dot{~u}$ & $-$ & $3.97 \times 10^{-2}$   & $3.89 \times 10^{-2}$ & $ 2.32\times 10^{-2}$ & \cellcolor{green!20}$1.94\times 10^{-2}$\\
      & $\sigma_{vm}$ & $-$ & $9.88\times 10^{-2}$ & $8.48\times 10^{-2}$ & $7.28\times 10^{-2}$ & \cellcolor{green!20}$7.17\times 10^{-2}$\\
     \hline
     CPU time & $-$& 52m 49.5s &33m 2.8s & 35m 41.5s & 31m 57.6s & \cellcolor{green!20}32m 46.7s\\
     \hline 
     
     Mean/max \# & & & & & & \cellcolor{green!20} \\
     Schwarz iterations & \multirow{-2}{*}{$-$} & \multirow{-2}{*}{3.21/5} & \multirow{-2}{*}{2.43/3}  & \multirow{-2}{*}{2.36/3} & \multirow{-2}{*}{2.34/3} & \cellcolor{green!20} \multirow{-2}{*}{ 2.28/3}
    \end{tabular}
\end{table}

\begin{table}[ht!]
    \centering
    \caption{3D nonlinear hyperelastic bolted joint problem, predictive regime: relative errors and performance metrics for various COpInf-COpInf O-SAM-based couplings.  The best coupling in terms of a combination of accuracy of the solution in $\Omega_b$ and efficiency is highlighted in green.}
    \label{tab:bolted-joint-metrics-predi-rom-rom}
    \begin{tabular}{c|c|c|c|c|c|c}
    & Field & FOM-FOM & \multicolumn{4}{c}{COpInf-COpInf} \\ 
    \hline
    $\romDimOneArg{}=\romDimOneArg{b}=\romDimOneArg{p}$ & $-$& $-$ & \cellcolor{green!20}10 &  15 & 20 & 24 \\
    \hline 
    $\lambda_b$ & $-$ & $-$ & \cellcolor{green!20}$1.0\times 10^{-4}$ & $1.0\times 10^{-4}$ &$1.0\times 10^{-4}$ & $1.0\times 10^{-5}$  \\
    \hline 
    $\lambda_p$ & $-$ & $-$ &\cellcolor{green!20} $1.0\times 10^{-4}$ & $1.0\times 10^{-4}$ & $1.0\times 10^{-4}$ & $1.0\times 10^{-4}$\\
    \hline 
    
    \multirow{3}{*}{$\Omega_b$ rel errors}& $~u$ & $-$ &\cellcolor{green!20}4.41$\times 10^{-2}$&$9.48 \times 10^{-2}$ & $ 1.29\times 10^{-1}$ & $4.51 \times 10^{-2}$\\ 
    & $\dot{~u}$ & $-$ &\cellcolor{green!20}$1.07 \times 10^{-1}$   & $3.60 \times 10^{-1}$ & $ 5.27\times 10^{-1}$  & $2.39 \times 10^{-1}$\\ 
   & $\sigma_{vm}$& $-$ &\cellcolor{green!20} $1.24 \times 10^{-1}$  & $2.38 \times 10^{-1}$& $5.29\times 10^{-1}$ & $1.72\times 10^{-1}$\\ 
    \hline
     \multirow{3}{*}{$\Omega_p$ rel errors}  & $~u$& $-$ & \cellcolor{green!20}$9.78 \times 10^{-3}$  &$ 1.62\times 10^{-2}$ & $1.25\times 10^{-2}$ & $7.26\times 10^{-3}$\\
      & $\dot{~u}$ & $-$ & \cellcolor{green!20}$4.36 \times 10^{-2}$   & $9.35 \times 10^{-2}$ & $8.57 \times 10^{-2}$ & $3.89\times 10^{-2}$\\
      & $\sigma_{vm}$ & $-$ & \cellcolor{green!20}$3.02\times 10^{-1}$ & $5.04 \times 10^{-1}$ & $5.96 \times 10^{-1}$& $3.09\times 10^{-1}$\\
     \hline
     CPU time & $-$& 52m 49.5s &\cellcolor{green!20} 8m 37.7s& 7m 41.2s & 8m 48.4s& 7m 35.1s \\
     \hline 

     Mean/max \# & & & \cellcolor{green!20} & & & \\
     Schwarz iterations & \multirow{-2}{*}{$-$} & \multirow{-2}{*}{3.21/5} & \cellcolor{green!20}\multirow{-2}{*}{2.43/3}  &\multirow{-2}{*}{2.47/3} &\multirow{-2}{*}{2.36/3} &  \multirow{-2}{*}{2.29/3}
    \end{tabular}
\end{table}

\begin{table}[ht!]
    \centering
    \caption{3D nonlinear hyperelastic bolted joint problem, predictive regime: sizes, ranks and condition numbers of the data matrices $~D_i^c$.}
    \label{tab:data_matrix_rank}
    \begin{tabular}{c||c||c|c||c|c}
     & & \multicolumn{2}{c ||}{$\Omega_b$} & \multicolumn{2}{c}{$\Omega_p$} \\ 
     \hline 
    $r_b=r_p = r$ & size of $~D_i^c$ & rank of $~D_b^c$ & condition number of $~D_b^c$  & rank of $~D_p^c$ & condition number of $~D_p^c$\\
    \hline 
    10 & $1120\times 58$ & 13 & $1.06\times 10^{29}$ & 13 & $2.91\times 10^{29}$ \\
   15 & $3630\times 58$ & 18 & $6.03\times 10^{27}$ & 16 & $5.28\times 10^{27}$ \\
   20 & $8440\times 58$ & 21 & $4.44\times 10^{26}$ & 21 & $7.85\times 10^{26}$ \\
   25 & $14448\times 58$ & 25 & $8.38\times 10^{25}$ & 24 & $8.88\times 10^{25}$
    \end{tabular}
\end{table}

%Remarkably, as seen in Figure \ref{fig:bolted-joint-predi-conv}(a), convergence with respect to the basis size is observed in the displacement and velocity fields for the FOM-COpInf models even in the presently-considered predictive regime.  %As before, the acceleration errors stagnate due to the POD bases' limited ability to represent the acceleration field. 
%Unfortunately, convergence with the basis size is not observed for the  OpInf-OpInf couplings considered (Figure \todo{add ref}(b)).  While an OpInf-OpInf coupled model with $M_b = M_p = 15$, corresponding to 99.999\% of the snapshot energy, gives a reasonably accurate solution, the accuracy degrades with increasing basis size.  

%\ikt{I am not including Pareto plots but I could.  Should I?  I didn't do it because all the FOM-COpInf and COpInf-COpInf models basically take close to the same CPU time.}

%% file: torsion.tex
\subsection{3D nonlinear hyperelastic torsion problem} \label{sec:torsion}

%\ikt{Question for Alejandro M.: I am not showing von Mises stresses for this problem b/c we never looked at it before for this problem.  Does that make sense?}

The next test case considers a nonlinear hyperelastic bar subjected to finite deformation by a high degree of torsion.  Our bar geometry has dimensions $0.05$ m $\times$ $0.05$ m $\times$ $1.0$ m, and is forced using the following initial conditions on the displacement and velocity, respectively:
\begin{equation} \label{eq:ic_torsion}
    ~u(~x,0) = ~0, \hspace{0.5cm} \dot{~u}(~x,0) = \left( \begin{array}{ccc} -b_1yz, & b_2xz, & 0\end{array}\right)^{\intercal}.
\end{equation}
In \eqref{eq:ic_torsion}, $b_1, b_2 \in \mathbb{R}$ are the rotation rates.  The boundary conditions are all homogeneous Neumann. We specify a nonlinear Neohookean-type material model with Young's modulus $E = 1.0 \times 10^9$ Pa, Poisson's ratio $\nu = 0.25 $ and density $\rho = 1000$ kg/m$^3$.  While these material properties are not realistic, they enable us to run this problem with an explicit time-stepper without having to resort to an extremely small time step to achieve CFL stability.  
For more details on the material model employed, the reader is referred to \cite{Mota:2011} and Appendix A.3.  
The problem is run from time 0 until time $2.0\times 10^{-3}$ s.  We intentionally choose not to run this problem with a linear elastic material model, as the solution to such a problem is known to be nonphysical.

\begin{figure}[ht!]
    \centering
    \begin{subfigure}{0.65\linewidth}
     \vspace{1cm}
        \includegraphics[width=0.99\linewidth]{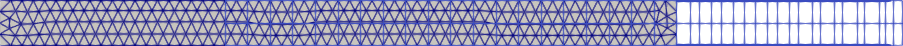}
        \subcaption{$\Omega_1$ }
         \vspace{1cm}
    \end{subfigure}
       \begin{subfigure}{0.65\linewidth}
        \includegraphics[width=0.99\linewidth]{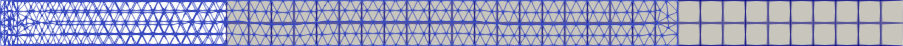}
        \subcaption{$\Omega_2$}
    \end{subfigure}
    \caption{3D nonlinear hyperelastic torsion problem: domain decomposition and meshes.  $\Omega_1$ and $\Omega_2$ and their meshes are shown in gray in subfigures (a) and (b), respectively.} 
    \label{fig:torsion-dd-meshes}
\end{figure}

The main objective of the torsion problem is to demonstrate that O-SAM is capable of coupling two regions of the 3D bar having different mesh resolutions, element types, time integration schemes, time-steps and models.  Toward this effect, we break up our geometry into two subdomains, 
% $\Omega_1 = (-0.025,0.025) \times (-0.025, 0.025) \times (-0.5, 0.25) $ and $\Omega_2 = (-0.025, 0.025) \times (-0.025, 0.025) \times (-0.25, 0.5)$ 
$\Omega_1 = (-0.025,0.025)^2 \times (-0.5, 0.25) $ and $\Omega_2 = (-0.025, 0.025)^2 \times (-0.25, 0.5)$ , and discretize them with a fine four-node tetrahedral and coarse eight-node hexahedral mesh, as shown in Figure \ref{fig:torsion-dd-meshes}.  The $\Omega_1$ mesh consists of 7333 elements and 1787 nodes (Figure \ref{fig:torsion-dd-meshes}(a)), whereas the $\Omega_2$ mesh has only 120 elements and 279 nodes (Figure \ref{fig:torsion-dd-meshes}(b)). Since the material model specified within the bar is nonlinear, a linear OpInf ROM will be incapable of capturing the problem dynamics.  Although the Neohookean-type material model utilized gives rise to a set of PDEs with generic nonlinearities, we choose to approximate these nonlinearities with a quadratic OpInf model, denoted by QOpInf.  %Since our goal herein is to study mixed couplings, 
We focus our attention on QOpInf-FOM couplings, in which a QOpInf ROM is prescribed in $\Omega_1$ and a FOM is prescribed in $\Omega_2$.  Since the solution in $\Omega_1$ is approximated by a relatively inexpensive ROM, we will advance this subdomain forward in time using an implicit Newmark-$\beta$ stepper with parameters $\beta = 0.25$ and $\gamma = 0.5$ and time-step $\Delta t = 2.0 \times 10^{-6}$ s.  For the FOM subdomain, $\Omega_2$, we utilize an explicit Newmark-$\beta$ scheme with $\gamma = 0.5$ and time-step $\Delta t = 1.0\times 10^{-6}$ s.  Within our SAM coupling algorithm, we specify a controller time-step of $2.0 \times 10^{-6}$ s.  We converge SAM to a relative tolerance of $\delta_{\text{rel}} = 1.0 \times 10^{-10}$ and an absolute tolerance of $\delta_{\text{abs}} = 1.0 \times 10^{-6}$.  For the implicit Newmark$-\beta$ runs, we employ the same relative and absolute tolerances for the Newton-based nonlinear solver and iterative linear solver.

As with the other test cases, we consider both a reproductive and a predictive variant of the torsion problem.  For both problem variants, we wish to predict the solution for $b_1 = b_2 = 5500$ in \eqref{eq:ic_torsion}.  For the reproductive version, an O-SAM-based explicit-implicit FOM-FOM coupling is performed using these parameter values with snapshots saved every $2.0 \times 10^{-5}$ s, yielding a total of 201 snapshots, from which a POD basis and QOpInf ROM is constructed.  For the predictive version of the problem, the parameters $b_1, b_2$ are each sampled from the parameter set $\{ 500, 1000, 5000, 8000\}$, and a total of 16 O-SAM-based implicit FOM-FOM coupling runs are performed.  As before, snapshots are saved every $2.0\times 10^{-5}$ s, this time yielding a total of 3216 snapshots, from which a QOpInf ROM is constructed in $\Omega_1$.  A singular value-based energy analysis reveals that $\romDimOneArg{1} = 30$ and $\romDimOneArg{1} = 27$ POD modes capture 99.999\% of the snapshot energy in $\Omega_1$ for the reproductive and predictive versions of the torsion problem, respectively; meanwhile, 3 and 6 POD modes capture 99.9999\% of the snapshot energy on the Schwarz boundary for the reproductive and predictive versions of this problem, respectively.  Our 
parameter sweep algorithm (see Section \ref{sec:regularization}) determined optimal regularization parameters of $1.0 \times 10^{-3}$ and $1.0 \times 10^{-11}$ respectively for these cases.

\begin{figure}[ht!]
    \centering
    \begin{subfigure}{0.45\linewidth}
        \includegraphics[width=0.99\linewidth]{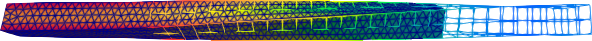}
        \subcaption{FOM-FOM}
    \end{subfigure}
       \begin{subfigure}{0.45\linewidth}
        \includegraphics[width=0.99\linewidth]{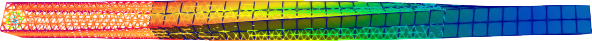}
        \subcaption{FOM-FOM}
    \end{subfigure}
       \begin{subfigure}{0.45\linewidth}
        \vspace{1cm}
        \includegraphics[width=0.99\linewidth]{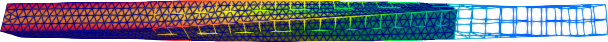}
        \subcaption{OpInf-FOM (reproductive)}
    \end{subfigure}
       \begin{subfigure}{0.45\linewidth}
       \vspace{1cm}
        \includegraphics[width=0.99\linewidth]{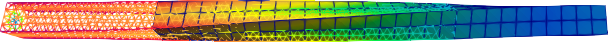}
        \subcaption{OpInf-FOM (reproductive)}
    \end{subfigure}
       \begin{subfigure}{0.45\linewidth}
        \vspace{1cm}
        \includegraphics[width=0.99\linewidth]{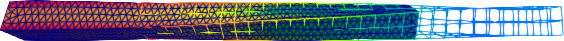}
        \subcaption{OpInf-FOM (predictive)}
    \end{subfigure}
       \begin{subfigure}{0.45\linewidth}
       \vspace{1cm}
        \includegraphics[width=0.99\linewidth]{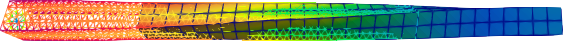}
        \subcaption{OpInf-FOM (predictive)}
    \end{subfigure}
     \begin{subfigure}{0.35\linewidth}
        \includegraphics[width=0.99\linewidth]{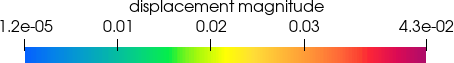}
       % \subcaption{blah}
    \end{subfigure}
      \vspace{0.5cm}
    \caption{3D nonlinear hyperelastic torsion problem: computed displacement solution magnitudes in $\Omega_1$ ((a), (c), (e)) and $\Omega_2$ ((b), (d), (f)) at the final time $2.0 \times 10^{-3}$ s for various O-SAM-based couplings.} 
    \label{fig:torsion-disp-solns}
\end{figure}

\begin{figure}[ht!]
    \centering
    \begin{subfigure}{0.45\linewidth}
        \includegraphics[width=0.99\linewidth]{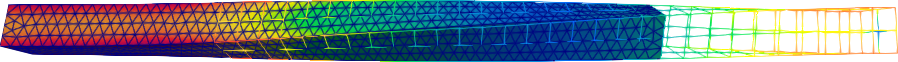}
        \subcaption{FOM-FOM}
    \end{subfigure}
       \begin{subfigure}{0.45\linewidth}
        \includegraphics[width=0.99\linewidth]{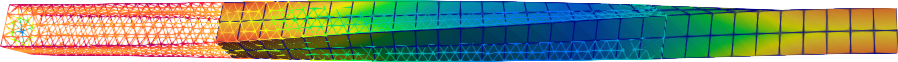}
        \subcaption{FOM-FOM}
    \end{subfigure}
       \begin{subfigure}{0.45\linewidth}
        \vspace{1cm}
        \includegraphics[width=0.99\linewidth]{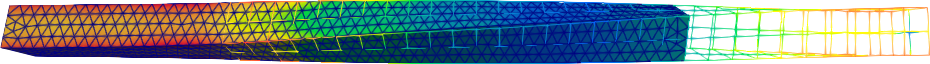}
        \subcaption{OpInf-FOM (reproductive)}
    \end{subfigure}
       \begin{subfigure}{0.45\linewidth}
       \vspace{1cm}
        \includegraphics[width=0.99\linewidth]{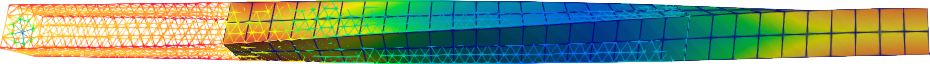}
        \subcaption{OpInf-FOM (reproductive)}
    \end{subfigure}
       \begin{subfigure}{0.45\linewidth}
        \vspace{1cm}
        \includegraphics[width=0.99\linewidth]{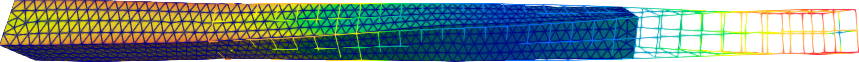}
        \subcaption{OpInf-FOM (predictive)}
    \end{subfigure}
       \begin{subfigure}{0.45\linewidth}
       \vspace{1cm}
        \includegraphics[width=0.99\linewidth]{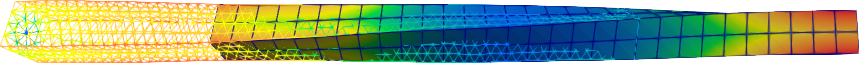}
        \subcaption{OpInf-FOM (predictive)}
    \end{subfigure}
     \begin{subfigure}{0.35\linewidth}
        \includegraphics[width=0.99\linewidth]{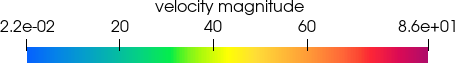}
       % \subcaption{blah}
    \end{subfigure}
      \vspace{0.5cm}
    \caption{3D nonlinear hyperelastic torsion problem: computed velocity solution magnitudes in $\Omega_1$ ((a), (c), (e)) and $\Omega_2$ ((b), (d), (f)) at the final time $2.0 \times 10^{-3}$ s  for various O-SAM-based couplings. 
    % \adg{We could consider showing only one of fig 20,21 if we need the space.}
    % \ikt{I think it's good to show both of them, to demonstrate there is nothing crazy going on in the velocity field.} \adg{OK, fine by me.}    
    } 
    \label{fig:torsion-velo-solns}
\end{figure}

Table \ref{tab:torsion-metrics} reports the displacement and velocity magnitude errors in both subdomains for QOpInf-FOM couplings in which the QOpInf ROM has either 30 or 27 POD modes, along with some performance metrics. The reader can observe that the QOpInf-FOM coupled models deliver accurate displacement solutions in both the reproductive and predictive regime, with maximum displacement magnitude errors of $0.267\%$ and $4.32\%$, respectively, and maximum velocity magnitude errors of $3.56\%$ and $14.9\%$, respectively  
%The velocity errors of $\mathcal{O}(1\%)$ and  $\mathcal{O}(10\%)$ for the reproductive and predictive models, respectively, are acceptable for some applications.  %The high errors in the acceleration field can be attributed to the fact that the FOM has general, rather than quadratic, nonlinearities, combined with the fact that the projection error for the acceleration snapshots is high, as shown in Figure \todo{add figure of projection errors}.  
As expected, the solutions are more accurate in the FOM subdomain, $\Omega_2$, than in the ROM subdomain, $\Omega_1$.  Figures \ref{fig:torsion-disp-solns}--\ref{fig:torsion-velo-solns} show plots of the reproductive and predictive QOpInf-ROM solutions compared to their FOM-FOM analogs for the displacement and velocity fields at the final simulation time, $2.0 \times 10^{-3}$ s.  The coupled QOpInf-FOM models are capable of capturing the displacement solution remarkably well.  %While there is some visible error in the velocity and acceleration fields, no coupling artifacts are observed.

The last three rows of Table \ref{tab:torsion-metrics} report the size, rank and condition number of the data matrix
\begin{equation} \label{eq:data_matrix_quadratic}
      ~D_1^q:= \left( \begin{array}{c} 
    \bar{~A}_1 \\
    \bar{~U}_1 \\
    \bar{~U}_1 \otimes \bar{~U}_1
    \end{array}\right),
\end{equation}
where $\bar{~A}_1$ and $\bar{~U}_1$ are as defined earlier in Section \ref{sec:bolted_joint_predi}.  As for the bolted joint problem, the reader can observe that the data matrices are extremely ill-conditioned and rank deficient.  This is a known issue for Operator Inference, and motivates the need for regularizing the OpInf minimization problem \eqref{eq:quadratic_opinf}.  

Finally, it is interesting to observe that the QOpInf-ROM couplings summarized in Table \ref{tab:torsion-metrics} require fewer Schwarz iterations to achieve convergence than their FOM-FOM analogs, as previously observed for our other test cases.  The reduced size of the QOpInf ROMs combined with this reduction in Schwarz iterations enables the coupled models to achieve speedups of up to $23.5\times$.

%\ikt{I can make plots of the error in each subdomain.  Is this worth adding?  Will lengthen the paper.}

%\ikt{Is it worth showing that linear OpInf-based couplings produce solutions that are completely wrong?}
%\irm{I think it is worth mentioning but probably enough just to state it doesn't work. Revisit after page length is close to final?}

\begin{table}[ht!]
    \centering
    \caption{3D nonlinear hyperelastic torsion problem: performance metrics for various O-SAM-based couplings.}
    \label{tab:torsion-metrics}
    \begin{tabular}{c|c|c|c|c}
    &  \multirow{2}{*}{Field}& \multirow{2}{*}{FOM-FOM} & QOpInf-FOM  & QOpInf-FOM \\ 
    & &  & reproductive & predictive \\
    \hline
    $\romDimOneArg{1}$  & $-$ & $-$ & 30 &  27 \\
    \hline 
    $\lambda_1$  & $-$ & $-$ & $1.0\times 10^{-3}$ & $1.0 \times 10^{-1}$\\
    \hline 
    \multirow{2}{*}{$\Omega_1$ rel errors}  & $~u$& $-$ &$2.67 \times 10^{-3}$ & $4.32\times 10^{-2}$ \\ 
    &$\dot{~u}$ &  $-$ &  $3.56\times 10^{-2}$ & $1.49 \times 10^{-1}$\\ 
    %& $-$ & $4.71\times 10^{-1}$ & $5.31 \times 10^{-1}$\\ 
    \hline
     \multirow{2}{*}{$\Omega_2$ rel errors} & $~u$&  $-$ & $1.13 \times 10^{-3}$ & $2.44 \times 10^{-2}$\\
     & $\dot{~u}$ &  $-$ &  $1.12 \times 10^{-2}$ & $9.52\times 10^{-2}$\\
    % & $-$ & $1.88\times 10^{-1}$ & $2.30 \times 10^{-1}$\\
     \hline
     CPU time & $-$& 39m 11.8s &1m 40.2s & 1m 39.5s \\
     \hline 
     
     Mean/max \#  & \multirow{2}{*}{$-$} &  \multirow{2}{*}{3.0/3} & \multirow{2}{*}{2.0/2} & \multirow{2}{*}{2.0/2} \\
     Schwarz iterations & & & & \\
     \hline 
    \textcolor{black}{ size of $~D_1^q$} & \textcolor{black}{$-$} & \textcolor{black}{$-$}& \multicolumn{2}{c}{ \textcolor{black}{$784\times 1000$}} \\
     \hline 
     \textcolor{black}{rank of $~D_1^q$} & \textcolor{black}{$-$} & \textcolor{black}{$-$}& \textcolor{black}{90} & \textcolor{black}{84}\\
     \hline 
     \textcolor{black}{condition number of $~D_i^q$} & \textcolor{black}{$-$} & \textcolor{black}{$-$}& \textcolor{black}{$3.53\times 10^{207}$} & \textcolor{black}{$2.97\times 10^{114} $}
    \end{tabular}
\end{table}

%% file: tension-specimen.tex
\subsection{3D nonlinear hyperelastic tension specimen problem} \label{sec:tension-specimen}

\begin{figure}[ht!]
    \centering
    \begin{subfigure}{0.3\linewidth}
\hspace{2cm}\includegraphics[width=0.2\linewidth]{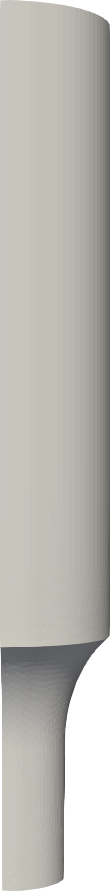}
        \subcaption{Geometry}
    \end{subfigure}
       \begin{subfigure}{0.3\linewidth}    \hspace{2cm}\includegraphics[width=0.2\linewidth]{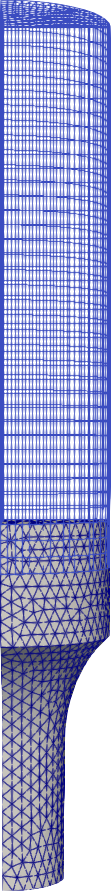}
        \subcaption{$\Omega_1$}
    \end{subfigure}
       \begin{subfigure}{0.3\linewidth}
\hspace{2cm}\includegraphics[width=0.2\linewidth]{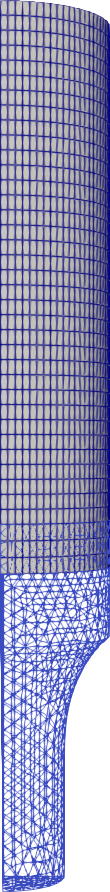}
        \subcaption{$\Omega_2$}
    \end{subfigure}
    
      \vspace{0.5cm}
    \caption{3D nonlinear hyperelastic tension-specimen problem: geometry (a) and domain decomposition into $\Omega_1$, discretized with a ten-node tetrahedral mesh (b), and $\Omega_1$, discretized using an eight-node hexahedral mesh (c). %\adg{this figure could be smaller, I think.} 
    } 
    \label{fig:tension-specimen-dd}
\end{figure}

The final test case on which we evaluate our coupling method is the tension specimen problem.  
Consider a uniaxial cylindrical tensile specimen made of aluminum.  In order to minimize computational 
cost, we perform a simulation of $\frac{1}{8}$ of the full tension specimen geometry (Figure \ref{fig:tension-specimen-dd}(a))
and apply symmetry boundary conditions to obtain a solution consistent with a simulation on the full geometry.
The height of the geometry on which the simulation is performed (Figure \ref{fig:tension-specimen-dd}(a)) is 4.445 cm.  
The smaller gauge radius is 3.81 m whereas the larger grip radius is 6.35 mm.  We decompose the geometry 
into two subdomains, one containing the gauge, denoted by $\Omega_1$, and one containing the grip, denoted by $\Omega_2$. 
The heights of $\Omega_1$ and $\Omega_2$ are 3.175 cm and 2.54 cm, respectively.  
In an effort to further assess O-SAM's ability to couple regions having different meshes and element types, we discretize $\Omega_1$ with 17,496 ten-node tetrahedral elements, while discretizing $\Omega_2$ with a coarser mesh having 10,800 eight-node hexahedral elements.  
The meshes contain 20,008 and 12,444 nodes, respectively.  As for the torsion problem, we specify a nonlinear Neohookean-type material model, 
this time with Young's modulus $E = 70 \times 10^9$ Pa, Poisson's ratio $\nu = 0.36$ and density $\rho = 2700$ kg/m$^3$.  Assuming our coordinate system is oriented such as ``up" is in the positive $y$--dimension and that the lower left corner of the geometry (Figure \ref{fig:tension-specimen-dd}(a)) is at the origin $(0,0,0)$, symmetry boundary conditions are specified for boundaries having $x=0$, $z=0$ and $y=0$, which amount to setting homogeneous Dirichlet boundary conditions for the $x$, $z$ and $y$ components of the displacement on these boundaries, respectively.  The problem is forced by setting the following dynamic boundary condition at the positive $y$ boundary: 
\begin{equation} \label{eq:tension-specimen-ic}
    ~u(~x,t) = \left(\begin{array}{ccc}0, &  
    \frac{1}{2}\alpha (1-\cos(\pi t)), & 0
    \end{array}\right)^{\intercal},
\end{equation}
for a specified scalar parameter $\alpha > 0$.  The problem is initialized with zero initial conditions for the displacement and the velocity, and advanced in time using an implicit Newmark-$\beta$ time-integrator having $\beta = 0.25$ and $\gamma = 0.5$ until time 1.0 s.  In applying O-SAM, we specify relative and absolute Schwarz tolerances of $\delta_{\text{rel}} = 1.0 \times 10^{-8}$ and 
$\delta_{\text{abs}} = 1.0 \times 10^{-6}$, respectively.  To prevent the Schwarz algorithm from taking too much CPU time, we cap the number of Schwarz iterations as 32.  For the Newton-based nonlinear and iterative linear solvers, we employ relative tolerances of $1.0 \times 10^{-7}$ and $1.0 \times 10^{-5}$, respectively.  %The problem is initialized with a $~0$ initial condition for both the displacement and velocity fields.

The objective herein is to construct FOM-QOpinf and QOpInf-QOpInf coupled models using O-SAM, and to evaluate 
them in the reproductive as well as the predictive regime.  Our goal is to capture the solution corresponding to the boundary condition \eqref{eq:tension-specimen-ic} with $\alpha = 0.005$.  For our FOM-QOpInf coupling, we assign the QOpInf ROM to the 
subdomain containing the holder, $\Omega_2$, where less complex dynamics are expected.  For the reproductive version of the problem, training data are generated by simulating the problem using a FOM-FOM O-SAM coupling with $\alpha = 0.005$ and a time-step of $1.0 \times 10^{-2}$ s in both subdomains.  For the predictive variant, training data are generated by performing FOM-FOM O-SAM couplings with $\alpha \in \{ 0.003, 0.006\}$.  We used the resulting snapshot sets to build POD bases in $\Omega_1$ and $\Omega_2$, each capturing 99.9999\% of the snapshot energy, which required a mere 2 modes in the subdomain interiors, 2 modes on the Schwarz boundary and 1 mode on the symmetry boundaries for both the reproductive and predictive cases.  We utilized our regularization parameter optimization algorithms to determine the optimal values of $\lambda_1$ and $\lambda_2$ in the two subdomains; these are reported in Tables \ref{tab:tension-specimen-metrics-repro} and \ref{tab:tension-specimen-metrics-predi}.  As for the torsion problem, we calculated the sizes, ranks and condition numbers for the quadratic data matrices $~D_i^q$ \eqref{eq:data_matrix_quadratic}.  These are $8 \times 100$, $7$ and $\mathcal{O}(10^{24})$, respectively.  While the matrices are not as rank deficient as for the previous examples, they are still poorly conditioned, which indicates regularization of the OpInf minimization problem \eqref{eq:quadratic_opinf} is a must.

The main results of our analysis are reported in Tables \ref{tab:tension-specimen-metrics-repro} and \ref{tab:tension-specimen-metrics-predi} for the reproductive and predictive cases, respectively.  The tables give relative errors in the displacement, velocity and von Mises stress fields, as well as the CPU times required to run the cases on {\tt Rigel}, and and mean and maximum number of Schwarz iterations required to reach convergence.  In addition to the classical iterative Schwarz algorithm advocated in this paper, termed ``full SAM" in Tables \ref{tab:tension-specimen-metrics-repro} and \ref{tab:tension-specimen-metrics-predi}, in which the Schwarz iteration is converged up to the specified $\delta_{\text{rel}}$ and $\delta_{\text{abs}}$ tolerances,
we consider also a simplified version of the Schwarz algorithm, in which a single Schwarz iteration is performed in each time step.  This second variant of O-SAM is roughly equivalent to the approach proposed in \cite{Farcas:2023} by Farcas \textit{et al.}, and is labeled ``1-iter SAM" in the tables and figures that follow. 

\begin{figure}[ht!]
    \centering
    \begin{subfigure}{0.3\linewidth}
\hspace{2cm}\includegraphics[width=0.2\linewidth]{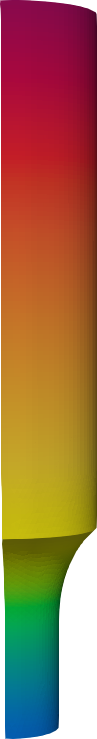}
        \subcaption{FOM-FOM, full SAM}
    \end{subfigure}
       \begin{subfigure}{0.3\linewidth}    \hspace{2cm}\includegraphics[width=0.2\linewidth]{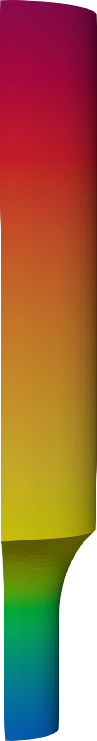}
        \subcaption{QOpInf-QOpInf, full SAM}
    \end{subfigure}
       \begin{subfigure}{0.3\linewidth}
\hspace{2cm}\includegraphics[width=0.2\linewidth]{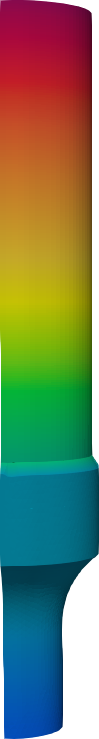}
        \subcaption{QOpInf-QOpInf, 1-iter SAM}
    \end{subfigure}
      \begin{subfigure}{0.99\linewidth}
      \vspace{0.5cm}
\hspace{4.5cm}\includegraphics[width=0.5\linewidth]{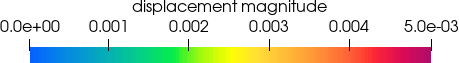}
    \end{subfigure}
     
      \vspace{0.5cm}
    \caption{3D nonlinear hyperelastic tension specimen problem, reproductive regime: plots of the displacement magnitude at the final time for a FOM-FOM O-SAM-based coupling (a) compared to a QOpInf-QOpInf full Schwarz coupling (b) and a QOpInf-QOpInf Schwarz coupling with only 1 Schwarz iteration (c).  } 
    \label{fig:tension-specimen-disp-solns}
\end{figure}

\begin{figure}[ht!]
    \centering
    \begin{subfigure}{0.3\linewidth}
\includegraphics[width=0.99\linewidth]{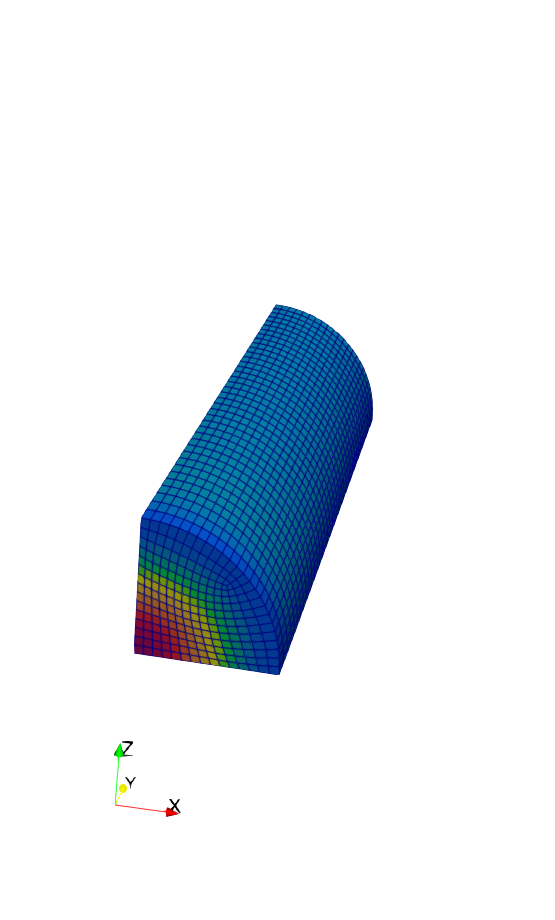}
        \subcaption{FOM-FOM, full SAM}
    \end{subfigure}
       \begin{subfigure}{0.3\linewidth}    \includegraphics[width=0.99\linewidth]{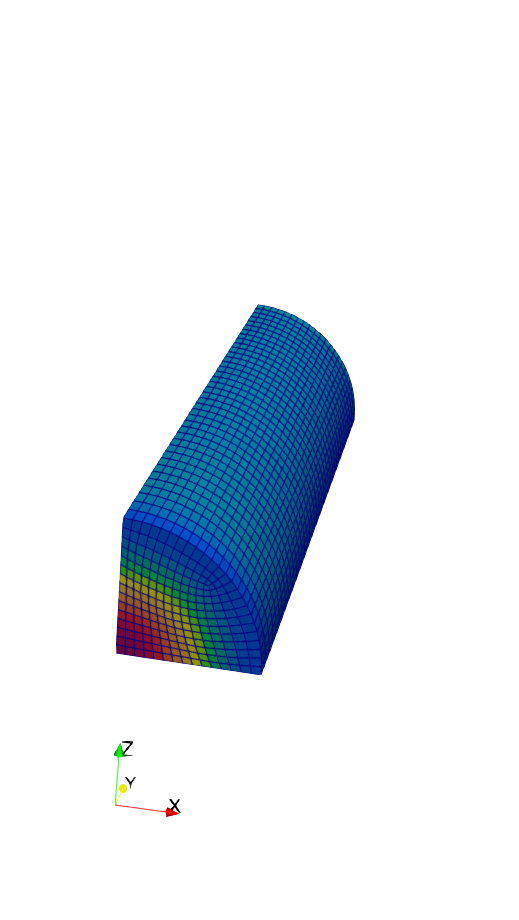}
        \subcaption{QOpInf-QOpInf, full SAM}
    \end{subfigure}
       \begin{subfigure}{0.3\linewidth}
\includegraphics[width=0.99\linewidth]{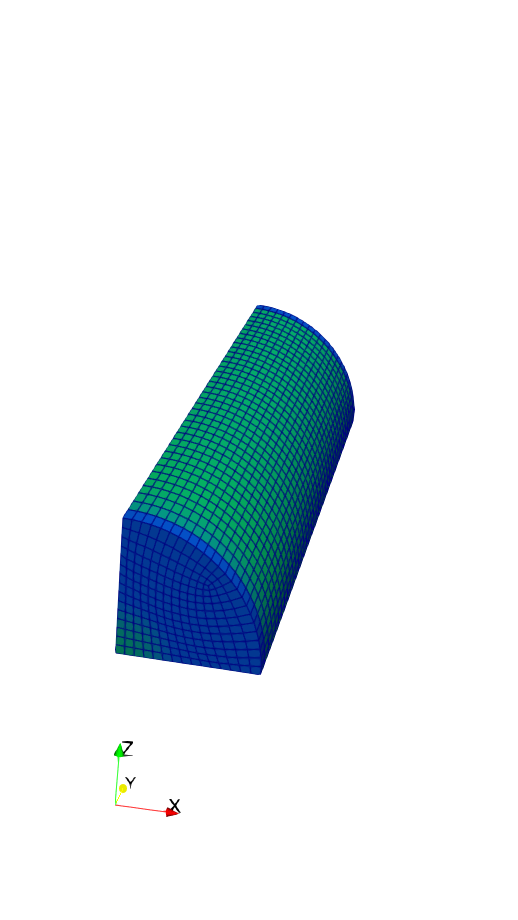}
        \subcaption{QOpInf-QOpInf, 1-iter SAM}
    \end{subfigure}
      \begin{subfigure}{0.99\linewidth}
      \vspace{0.5cm}
\hspace{4.5cm}\includegraphics[width=0.5\linewidth]{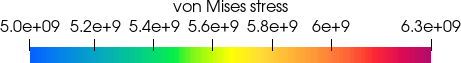}
    \end{subfigure}
     
      \vspace{0.5cm}
    \caption{3D nonlinear hyperelastic tension specimen problem, predictive regime: plots of the average von Mises stress $\sigma_{vm}$ in $\Omega_2$ at the final time for a FOM-FOM O-SAM-based coupling (a) compared to a QOpInf-QOpInf full Schwarz coupling (b) and a QOpInf-QOpInf Schwarz coupling with only 1 Schwarz iteration (c).  } 
    \label{fig:tension-specimen-stress-vm-solns}
\end{figure}

The reader can observe by examining Tables \ref{tab:tension-specimen-metrics-repro} and \ref{tab:tension-specimen-metrics-predi}  that all FOM-QOpInf and QOpInf-QOpInf models evaluated deliver relative errors of $\mathcal{O}(10^{-4})-\mathcal{O}(10^{-3})$ for the displacement and von Mises stress, and errors of $\mathcal{O}(10^{-3})-\mathcal{O}(10^{-2})$ %\adg{Do we need the $1.0\times$ in these big $\mathcal{O}'s$?} 
for the velocity when applying the full Schwarz algorithm.  When applying the 1-iteration version of SAM, the errors are several orders of magnitude higher, indicating convergence to the wrong solution, a result corroborated by Figures \ref{fig:tension-specimen-disp-solns}(c) and \ref{fig:tension-specimen-stress-vm-solns}(c), which show the displacement and von Mises stress solutions at the final time $1.0$ s for several such couplings.  It is interesting to observe that the full Schwarz method requires significantly fewer Schwarz iterations to converge when applied to our FOM-QOpInf and QOpInf-QOpInf couplings than the corresponding FOM-FOM coupling: whereas the FOM-FOM coupling fails to achieve the specified Schwarz tolerances in the maximum number of Schwarz iterations allowed (32), the hybrid couplings involving QOpInf ROMs converge in a mere 7--10 Schwarz iterations.  This enables O-SAM to achieve impressive speedups of up to $6.13\times$ and $106\times$ when performing FOM-QOpInf and QOpInf-QOpInf  couplings, all while maintaining good accuracy with respect to the corresponding FOM-FOM coupled solution (Tables \ref{tab:tension-specimen-metrics-repro}--\ref{tab:tension-specimen-metrics-predi} and Figures \ref{fig:tension-specimen-disp-solns}--\ref{fig:tension-specimen-stress-vm-solns}(a),(b)).  
As discussed earlier in Section \ref{sec:clamped}, we believe the reduction in the number of Schwarz iterations required for convergence when coupling in ROMs is due to the fact the ROM solutions tend to be smoother than their FOM analogs, making them easier to couple, together with the fact that the ROMs use problem-specific POD bases to represent the solution, instead of generic and problem-agnostic finite element shape functions.
%We believe there are two possible explanations for this behavior. First, since our OpInf models rely on POD modes to approximate the solution, these models give rise to solutions that are inherently smoother than their FEM analogues, which can aid convergence of the O-SAM coupling method. It is likely that convergence is also accelerated by the fact that the shape functions underlying the OpInf-OpInf couplings are problem-specific and data-driven, unlike the problem-agnostic and generic finite element shape functions underlying the FOM-FOM couplings.

\begin{table}[ht!]
\begin{small}
    \centering
    \caption{3D nonlinear hyperelastic tension specimen problem, reproductive regime: errors and performance metrics for various O-SAM-based couplings.  The best coupling in terms of a combination of the overall accuracy and efficiency is highlighted in green.  The errors highlighted in red are deemed unacceptable. %\adg{this table is too wide, we should shrink it.} 
    }
    \label{tab:tension-specimen-metrics-repro}
    \begin{tabular}{c|c|c|c|c|c|c}
    & \multirow{2}{*}{Field} &  FOM-FOM & FOM-QOpInf  & FOM-QOpInf &\cellcolor{green!20} QOpInf-QOpInf & QOpInf-QOpInf \\ 
     & & full SAM & full SAM & 1-iter SAM & \cellcolor{green!20}full SAM & 1-iter SAM \\
    \hline
    $\romDimOneArg{1}$ &$-$ & $-$ & $-$ & $-$ & \cellcolor{green!20} 2 & 2  \\
    \hline 
    $\romDimOneArg{2}$ & $-$ & $-$ & 2 & 2 & \cellcolor{green!20}2 & 2  \\
    \hline  
    $\lambda_1$ & $-$ &$-$ & $-$ & $-$ & \cellcolor{green!20}$1.0\times 10^{-6}$ & $1.0\times 10^{-6}$ \\
    \hline 
    $\lambda_2$ & $-$ & $-$ & $1.0\times 10^{-7}$ & $1.0\times 10^{-7}$ & \cellcolor{green!20}$1.0\times 10^{-11}$ & $1.0\times 10^{-11}$ \\
    \hline 
    \multirow{3}{*}{$\Omega_1$ rel errors}  & $~u$ & $-$ &  $1.01\times 10^{-4}$ & $1.88 \times 10^{-2}$  & \cellcolor{green!20}$5.68\times 10^{-4}$ &\cellcolor{red!20} $2.46 \times 10^1$  \\ 
    & $\dot{~u}$ & $-$ & $8.23\times 10^{-3}$  & $3.94\times 10^{-2}$ & \cellcolor{green!20}$1.83 \times 10^{-2}$ &\cellcolor{red!20} 4.49\\ 
    & $\sigma_{vm}$ &  $-$ & $1.04\times 10^{-4}$ & $1.93\times 10^{-2}$ & \cellcolor{green!20}$1.20\times 10^{-3}$ & \cellcolor{red!20}3.17\\ 
    \hline
     \multirow{3}{*}{$\Omega_2$ rel errors} &  $~u$& $-$ & $6.57\times 10^{-5}$ & $6.87\times 10^{-3}$ & \cellcolor{green!20}$2.38 \times 10^{-4}$ & \cellcolor{red!20}$1.03 \times 10^{-1}$\\
     & $\dot{~u}$&  $-$ & $8.91\times 10^{-3}$ & $4.54 \times 10^{-2}$ &\cellcolor{green!20}$7.60 \times 10^{-3}$ & \cellcolor{red!20}1.89 \\
     & $\sigma_{vm}$& $-$ & $8.51\times 10^{-4}$ & \cellcolor{red!20} $1.22\times 10^{-1}$ &\cellcolor{green!20}$8.94\times 10^{-4}$ & \cellcolor{red!20}$4.33\times 10^{-1}$  \\
     \hline
     CPU time & $-$&  8h 19m 29.5s&1h 34m 53.1s
& 12m 37.9s& \cellcolor{green!20}5m 20.3s & 1m 43.5  \\
     \hline 

     Mean/max \# & & & & & \cellcolor{green!20} & \\
     Schwarz iters  & \multirow{-2}{*}{$-$}&  \multirow{-2}{*}{32.0/32} & \multirow{-2}{*}{7.71/8} & \multirow{-2}{*}{1.0/1}  & \cellcolor{green!20}\multirow{-2}{*}{9.04/10}& \multirow{-2}{*}{1.0/1}
     
    \end{tabular}
    \end{small}
\end{table}

\begin{table}[ht!]
\begin{small}
    \centering
    \caption{3D nonlinear hyperelastic tension specimen problem, predictive regime: errors and performance metrics for various O-SAM-based couplings. The best coupling in terms of a combination of the overall accuracy and efficiency is highlighted in green.  The errors highlighted in red are deemed unacceptable.}
    \label{tab:tension-specimen-metrics-predi}
    \begin{tabular}{c|c|c|c|c|c|c}
    & \multirow{2}{*}{Field}  & FOM-FOM & FOM-QOpInf  & FOM-QOpInf &\cellcolor{green!20} QOpInf-QOpInf & QOpInf-QOpInf \\ 
    & &  full SAM & full SAM & 1-iter SAM & \cellcolor{green!20}full SAM & 1-iter SAM \\
    \hline
    $\romDimOneArg{1}$ & $-$& $-$ & $-$ & $-$ & \cellcolor{green!20}2 & 2  \\
    \hline 
    $\romDimOneArg{2}$ & $-$&  $-$ & 2 & 2 & \cellcolor{green!20}2 & 2  \\
    \hline 
    $\lambda_1$  & $-$ & $-$ & $-$ & $-$ & \cellcolor{green!20}$1.0 \times 10^{-11}$ &  $1.0 \times 10^{-11}$ \\ 
    \hline
    $\lambda_2 $ & $-$ & $-$ & $1.0\times 10^{-6}$ & $1.0\times 10^{-6}$ & \cellcolor{green!20}$1.0\times 10^{-6}$ & $1.0\times 10^{-6}$ \\
    \hline 
    \multirow{3}{*}{$\Omega_1$ rel errors} & $~u$&  $-$ & $3.44\times 10^{-4}$ & $5.62\times 10^{-2}$ &\cellcolor{green!20}$5.73 \times 10^{-4}$ & $2.65 \times 10^{-2}$   \\ 
    & $\dot{~u}$&  $-$ &$1.72 \times 10^{-2}$ & \cellcolor{red!20}1.42 & \cellcolor{green!20}$1.83 \times 10^{-2}$ & \cellcolor{red!20}$1.24 \times 10^{-1}$  \\ 
    & $\sigma_{vm}$& $-$ & $3.41\times 10^{-4}$ & $5.53 \times 10^{-2}$ & \cellcolor{green!20}$8.53\times 10^{-4}$ & $2.75\times 10^{-2}$\\ 
    \hline
     \multirow{3}{*}{$\Omega_2$ rel errors} & $~u$& $-$ & $2.50 \times 10^{-4}$ & $4.41 \times 10^{-2}$ & \cellcolor{green!20}$5.78 \times 10^{-4}$ & $1.92 \times 10^{-2}$\\
     & $\dot{~u}$& $-$ & $1.86 \times 10^{-2}$ &\cellcolor{red!20} 1.27 & \cellcolor{green!20}$1.96 \times 10^{-2}$ & \cellcolor{red!20}$1.91 \times 10^{-1}$\\
     & $\sigma_{vm} $& $-$ & $2.40\times 10^{-3}$ & \cellcolor{red!20}$2.04$ & \cellcolor{green!20}$6.00\times 10^{-3}$ &\cellcolor{red!20} $2.33 \times 10^{-1}$ \\
     \hline
     CPU time &$-$ &   8h 19m 29.5s & 1h 21m 25.9s&12m 43.5s& \cellcolor{green!20}4m 42.1s& 1m 49.3s

\\
     \hline 
     Mean/max \# & & & & & \cellcolor{green!20} & \\
     Schwarz iters  & \multirow{-2}{*}{$-$}& \multirow{-2}{*}{32.0/32} & \multirow{-2}{*}{7.03/8} & \multirow{-2}{*}{1.0/1}  & \cellcolor{green!20}\multirow{-2}{*}{7.74/8}& \multirow{-2}{*}{1.0/1}
    \end{tabular}
    \end{small}
\end{table}

%% file: 07-conclusion.tex
\section{Conclusions} \label{sec:conc}

In this work, we have presented a hybrid DD-based approach for coupling subdomain-local high-fidelity FOMs with non-intrusive OpInf ROMs using O-SAM. This method addresses significant challenges faced by analysts using traditional high-fidelity simulation codes, reducing both the long runtime requirements and the extensive mesh generation challenges that analysts often face. 
%The creation of high-quality meshes for complex multiscale components often represents a major bottleneck in the modeling and simulation workflow, sometimes taking weeks to complete. 
Our approach mitigates these issues by enabling the seamless ``gluing together" of arbitrary combinations of subdomain-local FOMs and non-intrusive OpInf ROMs in a plug-and-play fashion, thereby enhancing online efficiency, and  paving the way for more flexible and efficient simulation workflows in engineering applications.

Through a series of numerical experiments, we have demonstrated the efficacy of our approach on several complex 3D solid dynamics problems characterized by nonlinear behavior, implemented within the {\tt Norma.jl} \cite{Norma.jl} Julia code. The results indicate that our SAM-based coupling framework not only enhances computational efficiency by a factor as high as $106\times$ with respect to a comparable FOM-FOM coupling performed via O-SAM, but also maintains high accuracy across disparate models, meshes, and time integration schemes. Notably, the 
extension of our minimally-intrusive SAM-based coupling strategy to non-intrusive OpInf ROMs 
%the integration of OpInf allows for a non-intrusive 
%coupling strategy that 
significantly reduces the implementation burden associated with integrating data-driven models into existing mod/sim workflows.
%typically associated with traditional model order reduction techniques.  
Additionally, we demonstrate that subdomain-local quadratic OpInf (QOpInf) ROMs, when coupled with each other and with (fully nonlinear) FOMs, can deliver highly accurate solutions even when applied to problems with more generic (non-polynomial) nonlinearities.  We propose some innovative strategies to enhance the accuracy, robustness and efficiency of SAM when coupling subdomain-local OpInf ROMs, including a method for automatically optimizing subdomain-local OpInf regularization parameters and a mechanism for reducing the size of the learned boundary operator used in the Schwarz iteration process.  While attention is focused on solid dynamics, %we emphasize that 
our general O-SAM methodology can be applied to a wide range of PDEs.

The research summarized herein has informed several directions for future work.  First, as discussed in Section \ref{sec:analysis}, an extension of the theoretical analysis in \cite{Mota:2017, Mota:2022} to OpInf ROMs is currently lacking, and an ongoing research endeavor.  We are also currently in the process of developing a non-overlapping version of SAM that is capable of creating hybrid OpInf-OpInf and OpInf-FOM couplings \cite{Rodriguez:2025}.  This non-overlapping variant of SAM, termed NO-SAM, has several advantages, namely that it is more flexible  and can readily handle a wider range of problems than O-SAM, e.g., multi-material problems, multi-physics problems and problems with interfaces.  Additionally, we are exploring mechanisms for accelerating SAM through various strategies such as optimized  transmission conditions \cite{Rodriguez:2025}, Aitken \cite{deparis2006domain} and Anderson \cite{Walker:2011} acceleration, and the introduction of additional parallelism in the form of additive Schwarz \cite{Gander:2008}.
Towards enhancing the accuracy of SAM-based ROM-ROM and ROM-FOM-based couplings for highly nonlinear problems, e.g., problems with plasticity, we are extending SAM to work with non-intrusive structure-preserving NN-based OpInf ROMs \cite{erics_inprep_paper} and other classes of non-intrusive ROMs, such as kernel manifold ROMs \cite{diaz2025interpretableflexiblenonintrusivereducedorder, diaz2025kernelmanifoldsnonlinearaugmentationdimensionality}.
We are also developing an adaptive algorithm that will enable error indicator-driven online switching between ROMs and FOM in a way that manages both accuracy and efficiency.  
This adaptive algorithm is expected to enable the construction of robust hybrid models in which each subdomain-local ROM is trained on local (uncoupled) simulation data via a ``bottom-up" (rather than ``top-down") training approach.
Towards improving SAM's usability, we are developing a workflow for auto-tuning various SAM inputs  (e.g., interface location, subdomain count, overlap size) through multi-objective optimization (e.g., by simultaneously minimizing the CPU time and a stress recovery-based error indicator) using the {\tt GPTune} library \cite{gptune}, which performs a Bayesian, gradient-free optimization of black-box models using interpretable Gaussian process (GP) surrogates.
%\todo{am I missing any planned future work?} %Could mention different regularization parameters for different terms in the OpInf regularization problem, but I kind of think this is too in-the-weeds especially if we are moving in the direction of nonlinear OpInf.}  
Finally, we are beginning the implementation of our OpInf-FOM O-SAM-based couplings within Sandia's production code, {\tt SIERRA/SM} \cite{sierrasm}, towards making it accessible to analysts running a variety of mission-critical problems on complex geometries. 

%Future work will focus on refining the optimization of regularization parameters within the OpInf framework and exploring the applicability of our method to a broader range of physical phenomena. Additionally, we aim to investigate the potential for further enhancing the robustness and stability of the coupling process, particularly in the context of highly nonlinear systems.

%In conclusion, the proposed SAM-based coupling strategy represents a significant advancement in the field of multiscale modeling and simulation, offering a promising pathway for the development of efficient and accurate computational tools in engineering and scientific research.

%% file: 08-appendix.tex
\section*{Appendix A. Constitutive models}   \label{sec:appendix}

%\todo{Have Alejandro double check this section.}

\subsection*{A.1. Linear elastic material model}

For a simple linear elastic material, the Helmholtz free energy density %\todo{Is this same as strain energy?}
takes the form
\begin{equation}
    A(~F):= \frac{1}{2} \lambda (\text{tr}(~\epsilon))^2 + \mu \text{tr}(~\epsilon^2), 
\end{equation}
where $~\epsilon := \frac{1}{2} (\nabla ~u + (\nabla ~u)^{\intercal})$ and $\nabla ~u := ~F - ~I$, with $~F:= \nabla ~\varphi$ denoting the deformation gradient, $~I$ denoting the $3\times 3$ identity matrix, and $\lambda, \mu > 0$ denoting the Lam\'{e} coefficients.  

The governing PDEs now take the form \eqref{eq:dynamic_elasticity_pde}, where  
\begin{equation}
~\sigma := ~C : ~\epsilon,
\end{equation}
where $~C$ is the fourth-order stiffness tensor.

\subsection*{A.2.  Saint Venant--Kirchhoff material model}

In addition to describing the mathematical formulation of the Saint Venant--Kirchhoff material model, we also provide a derivation demonstrating that this constitutive model gives rise to PDEs with cubic nonlinearities. 

For the Saint Venant--Kirchhoff material model, it is well-known \cite{holzapfel2000} that the Helmholtz free-energy density $A(~F)$ implicit in the Piola-Kirchhoff stress tensor $~P$ takes the form 
\begin{equation} \label{eq:svk_W}
A(~F) = A(~E)=\frac{\lambda}{2} (\text{tr} ~E)^2 + \mu \text{tr}(~E^2),     
\end{equation}
where $~E := \frac{1}{2}~F^{\intercal}~F - ~I$, with $~F:= ~I + \nabla ~u $ denoting the deformation gradient, and where $\lambda, \mu > 0$ are the Lam\'{e} parameters.  
To derive the strong form of the dynamic solid mechanics PDEs, it is necessary to calculate the first Piola-Kirchhoff stress $~P := \frac{\partial A}{\partial ~F}$.  We first calculate the second Piola-Kirchhoff stress, $~S := \frac{\partial A}{\partial ~E}$.    The derivative of the first term in 
\eqref{eq:svk_W} is: 
% \begin{equation}
%     \frac{\partial (\text{tr}~E)^2 }{\partial ~E} = 2 \text{tr} ~E \frac{\partial (\text{tr}~E)} {\partial ~E}  = 2 \text{tr}~E \frac{\partial (E_{ij} \delta_{ij})}{\partial E_{ij}} = 2 \text{tr}~E \delta_{ij}= 2 \text{tr}(~E) ~I 
%  \end{equation}
\begin{equation}
    \frac{\partial (\text{tr}~E)^2 }{\partial ~E} = 2 \text{tr} ~E \frac{\partial (\text{tr}~E)} {\partial ~E}  = 2 \text{tr}~E \frac{\partial (~E:~I)}{\partial ~E} = 2 \text{tr}(~E) ~I .
 \end{equation}
 Since the second term in \eqref{eq:svk_W} is $E_{ij} E_{ij} = ~E : ~E$, its derivative with respect to $~E$ is $2~E$.  
 It follows that 
\begin{equation} \label{eq:svk_S}
     ~S := \frac{\partial A}{\partial ~E} = \lambda \text{tr}(~E) ~I + 2\mu ~E.
\end{equation}
  Remark that 
  \begin{equation}
  ~E = \frac{1}{2} (~F^{\intercal} ~F - ~I) = \frac{1}{2}\left[ (~I + \nabla ~u)^{\intercal} (~I + \nabla ~u) - ~I\right] = \frac{1}{2} \left[\nabla ~u + (\nabla ~u)^T + (\nabla ~u)^{\intercal} \nabla ~u\right].   
  \end{equation}
  Then,
  \eqref{eq:svk_S} can be rewritten as 
 \begin{equation} \label{eq:svk_S2}
     ~S = \frac{1}{2} \lambda \left[ 2 \text{tr} (\nabla ~u)  + \text{tr}((\nabla ~u)^{\intercal} \nabla ~u)\right] ~I + \mu  \left[\nabla ~u + (\nabla ~u)^{\textcolor{blue}{T}} + (\nabla ~u)^{\intercal} \nabla ~u\right].
 \end{equation}
% Now, the Hamiltonian of the PDE takes the form 
%\begin{equation}
%    H(~u, \dot{~u}) = \frac{1}{2} \int_{\Omega} \rho |\dot{~u}|^2 dV + \int_{\Omega} A(~u) dV, 
%\end{equation}
%with $~A(u)$ given by \eqref{eq:svk_W}.  

To calculate the 
first Piola-Kirchhoff stress $~P$, we compute
\begin{equation} \label{eq:svk_P}
\begin{array}{rl}
    ~P &:= \frac{\partial A}{\partial ~F}  \\
    &= \frac{\partial A}{\partial ~E} \frac{\partial ~E}{\partial ~F} \\
    &= ~F ~S \\
    &=  \frac{1}{2} \lambda ~F \left[ 2 \text{tr} (\nabla ~u)  + \text{tr}((\nabla ~u)^{\intercal} \nabla ~u)\right] ~I + \mu  \left[\nabla ~u + (\nabla ~u)^{\intercal} + (\nabla ~u)^{\intercal} \nabla ~u\right] \\
    &=  \frac{1}{2} \lambda ( ~I + \nabla ~u) \left[ 2 \text{tr} (\nabla ~u)  + \text{tr}((\nabla ~u)^{\intercal} \nabla ~u)\right] ~I + \mu   ( ~I + \nabla ~u)\left[\nabla ~u + (\nabla ~u)^{\intercal} + (\nabla ~u)^{\intercal} \nabla ~u\right]. \\
    \end{array}
\end{equation}
Now, the strong form of the PDE is \eqref{eq:dynamic_elasticity_pde},
%\begin{equation} \label{eq:svk_strong}
%\rho \ddot{~u} = \text{Div} ~P + \rho ~B,
%\end{equation}
where $~P$ is given by \eqref{eq:svk_P}.  In particular, from \eqref{eq:svk_P}, one can see that the governing PDE \eqref{eq:dynamic_elasticity_pde} has just cubic nonlinearities.

\subsection*{A.3.  Neohookean material model}  

For the Neohookean material model employed herein, the Helmholtz free-energy density is decomposed into a volumetric and deviatoric component 
\begin{equation} \label{eq:A_nh}
A(~F) := A^{\text{vol}}(~F) + A^{\text{dev}}(~F),
\end{equation}
where
\begin{equation}  \label{eq:Avol}
    A^{\text{vol}}(~F) :=  \frac{1}{4} \kappa (\det(~C) - \log(\det(~C)) - 1), 
\end{equation}
\begin{equation} \label{eq:Adev}
    A^{\text{dev}}(~F) := \frac{1}{2} \mu ( (\det(~C))^{-1/3} \text{tr}(~C) - 3), 
\end{equation}
where $~C:= ~F^{\intercal} ~F$ and $~F$ denotes the deformation gradient, defined above.  %\todo{Double check $~F$ definition with Alejandro.}  
In \eqref{eq:Avol}--\eqref{eq:Adev}, $\kappa$ denotes the bulk modulus and $\mu$ is a Lam\'{e} parameter. From \eqref{eq:A_nh}, one can shown that 
\begin{equation}
    ~S := ~S^{\text{vol}} + ~S^{\text{dev}},
\end{equation}
where 
\begin{equation}
   ~S^{\text{vol}}:= \frac{1}{2} \kappa (\det{~C} - 1) ~C^{-1},
\end{equation}
\begin{equation}
    ~S^{\text{dev}}:= \mu (\det(~C))^{-1/3} \left(~I - \frac{1}{2}~C^{-1} \text{tr}(~C)\right),
\end{equation}
where $~I$ is the $3 \times 3$ identity matrix.  Now, the first Piola-Kirchhoff stress takes the form $~P:= ~F~S$, as before.

%% file: reproducibility.tex
\section*{Code availability and reproducibility}\label{sec:reproducibility}

The {\tt Norma.jl} and {\tt norma-opinf} codes used to generate the results presented herein are available on GitHub at {\tt https://github.com/sandialabs/Norma.jl} and {\tt https://github.com/sandialabs/norma-opinf}, respectively.  In particular, we used the {\tt d11e20114a9aeb9d2427c03a946c89fee729b3d1} sha of {\tt Norma.jl} and the {\tt 763e37198e197919f4f78dcee4b26efb4809c1cd} sha of {\tt norma-opinf}.  Input files for variants of the clamped, bolted joint, torsion and tension-specimen problems can be found in the {\tt Norma.jl/examples/ahead} directory of {\tt Norma.jl}.

%% file: acknowledgments.tex
\section*{Acknowledgements} \label{sec:acknowl}

Support for this work was received through Sandia National Laboratories' Laboratory Directed Research and Development (LDRD) program and through the U.S. Department of Energy, Office of Science, Office of Advanced Scientific Computing Research, Mathematical Multifaceted Integrated Capability Centers (MMICCs) program, under Field Work Proposal 22025291 and the Multifaceted Mathematics for Predictive Digital Twins (M2dt) project. Additionally, the writing of this manuscript was funded in part by Irina Tezaur’s Presidential Early Career Award for Scientists and Engineers (PECASE).

Sandia National Laboratories is a multi-mission laboratory managed and operated by National Technology and Engineering Solutions of Sandia, LLC., a wholly owned subsidiary of Honeywell International, Inc., for the U.S. Department of Energy’s National Nuclear Security Administration under contract DE-NA0003525.

%The authors wish to thank Alejandro Mota for creating the {\tt Norma.jl} code in which our numerical examples are implemented, and for assisting with the setup of our numerical experiments within this code.  